\documentclass[reqno,11pt]{book}
\usepackage{psfig, amsmath, amstexnb, amsthm}
\usepackage{amssymb}  
\input dynsys.sty

\begin{document}

\thispagestyle{empty}

\begin{center}
{\Huge Geometrical Theory\\
 of Dynamical Systems\\}
\vspace{2cm}
{\Large
Nils Berglund\\
Department of Mathematics\\
ETH Z\"urich\\
8092 {\bf Z\"urich}\\
Switzerland\\
}
\vspace{2cm}
{\large Lecture Notes\\
Winter Semester 2000-2001\\}
\vspace{2cm}
{\Large Version:\\
November 14, 2001}
\end{center}

\chapter*{Preface}

This text is a slightly edited version of lecture notes for a course I gave
at ETH, during the Winter term 2000-2001, to undergraduate Mathematics and
Physics students. The choice of topics covered here is somewhat arbitrary,
and was partly imposed by time limitations. Evidently, these notes do not
intend to replace the many existing excellent textbooks on the subject, a
few of which are listed in the bibliography, but they might provide a
reasonably concise, albeit certainly biased introduction to this huge
domain. The approach used here has probably been influenced by my first
teacher in Dynamical Systems, Prof.\ Herv\'e Kunz. I also wish to
acknowledge my student's contribution in mercilessly tracking down a
substantial amount of typos. 

\bigskip
\noindent
Files available at 
{\tt http://www.math.ethz.ch/$\sim$berglund}

\noindent
Please send any comments to 
{\tt berglund@math.ethz.ch}

\bigskip
\noindent
\begin{flushright}
Z\"urich, November 2001
\end{flushright}

\addtolength\textheight{2mm}
\tableofcontents
\goodbreak
\addtolength\textheight{-2mm}
\newpage
\setcounter{page}{0}
\newpage

\chapter{Examples of Dynamical Systems}
\label{ch_ex}

The last 30 years have witnessed a renewed interest in dynamical systems,
partly due to the ``discovery'' of chaotic behaviour, and ongoing research
has brought many new insights in their behaviour.  What are dynamical
systems, and what is their geometrical theory? Dynamical systems can be
defined in a fairly abstract way, but we prefer to start with a few
examples of historical importance before giving general definitions. 
This will allow us to specify the class of systems that we want to study,
and to explain the differences between the geometrical approach and other
approaches.


\section{The Motion of the Moon}
\label{sec_moon}

The problem of the Moon's motion is a particular case of the $N$-body
problem, which gives a nice illustration of the historical evolution that
led to the development of the theory of dynamical systems. This section
follows mainly Gutzwiller's article \cite{Gutzwiller98}. 

Everyone knows that the phases of the Moon follow a cycle of a bit less
than 30 days. Other regularities in the Moon's motion were known to the
Babylonians as early as 1000 B.C. One can look, for instance, at the time
interval between Sunset and Moonrise at Full Moon. This interval is not
constant, but follows a cycle over 19 years, including 235 Full Moons (the
\defwd{Metonic Cycle}). Solar and Lunar Eclipses also follow a cycle with a
period of 18 years and 11 days, containing 223 Full Moons (the \defwd{Saros
Cycle}).

Greek astronomy started in the 5th century B.C.\ and initiated developments
culminating in the work of Ptolemy in the second century A.D. In contrast
with the Babylonians, who looked for regularities in long rows of numbers,
the Greeks introduced geometrical models for their astronomical
observations. To account for the various observed deviations from
periodicity, they invented the model of \defwd{epicycles}. In modern
notation, and assuming a planar motion with Cartesian coordinates
$(x,y)\in\R^2$, the complex number $z=x+\icx y\in\C$ evolves as a function
of time $t$ according to the law
\begin{equation}
\label{i1}
z = a \e^{\icx \w_1 t} (1+\eps \e^{\icx \w_2 t}),
\end{equation} 
where $a$, $\eps$, $\w_1$ and $\w_2$ are parameters which are fitted to
experimental data. 

The epicycle model was refined in subsequent centuries, with more terms being
included into the sum \eqref{i1} to explain the various ``inequalities''
(periodic deviations from the uniform motion of the Moon). Four inequalities
were discovered by Tycho Brahe alone in the 16th century. These terms could
be partly explained when Kepler discovered his three laws in 1609:
\begin{enum}
\item	the trajectory of a planet follows an ellipse admitting the Sun as a
focus,
\item	equal areas, measured with respect to the Sun, are swept in equal
time intervals,
\item	when several planets orbit the Sun, the period of the motion
squared is proportional to the third power of the semi-major axis of the
ellipse.
\end{enum}
Expanding the solution into Fourier series produces sums for which
\eqref{i1} is a first approximation.  However, while these laws describe
the motion of the planets quite accurately, they fail to fit the
observations for the Moon in a satisfactory way. 

A decisive new point of view was introduced by Newton when he postulated
his law of Universal Gravitation (published in his {\em Principia} in
1687). A system with $N$ planets is described by a set of \defwd{ordinary
differential equations}
\begin{equation}
\label{i2}
m_i \dtot{^2x_i}{t^2} = \sum_{\stackrel{j=1,\dots,N}{j\neq i}} 
\frac{Gm_im_j(x_j-x_i)}{\norm{x_i-x_j}^3}, 
\qquad i=1,\dots,N.
\end{equation}
Here the $x_i\in\R^3$ are vectors specifying the position of the planets,
the $m_i$ are positive scalars giving the masses of the particles, and $G$
is a universal constant. Newton proved that for two bodies ($N=2$), the
equation \eqref{i2} is equivalent to Kepler's first two laws. With three or
more bodies, however, there is no simple solution to the equations of
motion, and Kepler's third law is only valid approximately, when the
interaction between planets is neglected.

The three-body problem initiated a huge amount of research in the following
two hundred years. Newton himself invented several clever tricks allowing
him to compute corrections to Kepler's laws in the motion of the Moon. He
failed, however, to explain all the anomalies. Perturbation theory was
subsequently systematized by mathematicians such as Laplace,  Euler,
Lagrange, Poisson and Hamilton, who developed the methods of analytical
mechanics. As a first step, one can introduce the \defwd{Hamiltonian
function}
\begin{equation}
\label{i3}
\fctndef{H}{(\R^3)^N\times(\R^3)^N}
{\R}{(p,q)}{\displaystyle \sum_{i=1}^N \frac{p_i^2}{2m_i} -
\sum_{i<j}\frac{Gm_im_j}{\norm{q_i-q_j}},}
\end{equation}
where $p_i=m_iv_i\in\R^3$ are the momenta of the planets and
$q_i=x_i\in\R^3$ for $i=1,\dots,N$. The equation of motion \eqref{i2} is
then equivalent to the equations
\begin{equation}
\label{i4}
\dtot{q_i}t = \dpar{H}{p_i}, 
\qquad
\dtot{p_i}t = -\dpar{H}{q_i}.
\end{equation}
One advantage of this formulation is that all the information on the motion
is contained in the scalar function $H$. The main advantage, however, is
that the structure \eqref{i4} of the equations of motion is preserved under
special changes of variables, called \defwd{canonical transformations}.  In
the case of the two-body problem, a good set of coordinates is given by the
\defwd{Delaunay variables} $(I,\ph)\in\R^3\times\T^3$ (actually, there are
$6+6$ variables, but $6$ of them correspond to the trivial motion of the
center of mass of the system). The \defwd{action variables} $I_1, I_2, I_3$
are related to the semi-major axis, eccentricity and inclination of the
Kepler ellipse, while the \defwd{angle variables} $\ph_1, \ph_2, \ph_3$
describe the position of the planet and the spatial orientation of the
ellipse. The two-body Hamiltonian takes the form
\begin{equation}
\label{i5}
H(I,\ph) = -\frac{\mu}{2I_1^2},
\end{equation}
where $\mu=\frac{m_1m_2}{m_1+m_2}$. 
The equations of motion are 
\begin{align}
\nonumber
\dtot{\ph_1}t &= \dpar{H}{I_1} = \frac{\mu}{I_1^3}, &
\dtot{I_1}t &= -\dpar{H}{\ph_1} = 0 \\
\label{i6}
\dtot{\ph_2}t &= \dpar{H}{I_2} = 0, &
\dtot{I_2}t &= -\dpar{H}{\ph_2} = 0 \\
\nonumber
\dtot{\ph_3}t &= \dpar{H}{I_3} = 0, &
\dtot{I_3}t &= -\dpar{H}{\ph_3} = 0, 
\end{align}
describing the fact that the planet moves on an elliptical orbit with fixed
dimensions and orientation. 

In the case of the three-body problem Moon--Earth--Sun, one can use two sets
of Delaunay variables $(I,\ph)$ and $(J,\psi)$ describing, respectively, the
motion of the system Moon--Earth, and the motion around the Sun of the
center of mass of the system Moon--Earth. The Hamiltonian takes the form 
\begin{equation}
\label{i7}
H(I,J,\ph,\psi) = H_0(I,J) + H_1(I,J,\ph,\psi).
\end{equation}
The unperturbed part of the motion is governed by the Hamiltonian
\begin{equation}
\label{i8}
H_0(I,J) = -\frac{\mu}{2I_1^2}-\frac{\mu'}{2J_1^2},
\end{equation}
where  $\mu'=\frac{(m_1+m_2)m_3}{m_1+m_2+m_3}$. Due to the special initial
conditions of the system, the \defwd{perturbing function} $H_1$ has a small
amplitude. It depends on several small parameters: the initial
eccentricities $\eps\simeq 1/18$ of the Moon and $\eps'\simeq 1/60$ of the
Earth, their inclinations $i$ and $i'$, and the  ratio $a/a'\simeq 1/400$
of the semi-major axes of the two subsystems. All these quantities are
functions of the actions $I$ and $J$. The standard approach is to expand
$H_1$ in a trigonometric series
\begin{equation}
\label{i9}
H_1 = -\frac{Gm_1m_2m_3}{m_1+m_2} \frac{a^2}{{a'}^3}
\sum_{j\in\Z^6} C_j \e^{\icx(j_1\ph_1 + j_2\ph_2 + j_3\ph_3 + j_4\psi_1 +
j_5\psi_2 + j_6\psi_3)}. 
\end{equation}
The coefficients $C_j$ are in turn expanded into Taylor series of the small
parameters,
\begin{equation}
\label{i10}
C_j = \sum_{k\in\N^5} c_{jk}
\Bigpar{\frac{a}{a'}}^{k_1} \eps^{k_2} {\eps'}^{k_3}
\gamma^{k_4} {\gamma'}^{k_5},
\end{equation}
where $\gamma=\sin i/2$ and $\gamma'=\sin i'/2$. The solutions can then be
expanded into similar series, thus yielding a Fourier expansion of the form
\eqref{i1} (in fact, it is better to simplify the Hamiltonian by successive
canonical transformations, but the results are equivalent). The most
impressive achievement in this line of work is due to Delaunay, who
published in 1860 and 1867 two volumes of over 900 pages. They contain
expansions up to order 10, which are simplified with 505 transformations.
The main result for the trajectory of the Moon is a series containing 460
terms, filling 53 pages. 

At the turn of the century, these perturbative calculations were criticized
by Poincar\'e, who questioned the convergence of the expansions. Indeed,
although the magnitude of the first few orders decreases, he showed that
this magnitude may become  extremely large at sufficiently high order. This
phenomenon is related to the problem of small divisors appearing in the
expansion, which we will discuss in a simpler example in the next section.

Poincar\'e introduced a whole set of new methods to attack the problem from
a geometric point of view. Instead of trying to compute the solution for a
given initial condition, he wanted to understand the qualitative nature of
solutions for all initial conditions, or, as we would say nowadays, the
geometric structure of phase space. He thereby introduced concepts such as
invariant points, curves and manifolds. He also provided examples where the
solution cannot be written as a linear combination of periodic terms, a
first encounter with chaotic motion. 

The question of convergence of the perturbation series continued nonetheless
to be investigated, and was finally solved in a series of theorems by
Kolmogorov, Arnol'd and Moser (the so-called \defwd{KAM theory}) in the
1950s. They prove that the series converges for (very) small perturbations,
for initial conditions living on a Cantor set.

This did not solve the question of the motion of the Moon completely,
although fairly accurate ephemerides can be computed for relatively short
time spans of a few decades. Using a combination of analytical and
numerical methods, the existence of chaos in the Solar System was
demonstrated by Laskar in 1989 \cite{Laskar89}, implying that exact
positions of the planets cannot be predicted for times more than a few
hundred thousand years in the future.


\section{The Standard Map}
\label{sec_sm}

The standard map describes the motion of a ``rotator'' with one angular
degree of freedom $q\in\fS^1$ ($\fS^1$ denotes the circle $\R/2\pi\Z$),
which is periodically kicked by a pendulum-like force of intensity
proportional to $-\sin q$. If $q_n$ and $p_n$ denote the position and
momentum just before the \nth{n} kick,  one has
\begin{equation}
\label{sm1}
\begin{split}
q_{n+1} &= q_n + p_{n+1} 	\qquad\pmod{2\pi} \\
p_{n+1} &= p_n - \eps\sin q_n.
\end{split}
\end{equation}
For $\eps=0$, the dynamics is very simple and one has explicitly
\begin{equation}
\label{sm2}
\begin{split}
q_n &= q_0 + n p_0  	\qquad\pmod{2\pi} \\
p_n &= p_0.
\end{split}
\end{equation}
Let us now analyse the iterated map \eqref{sm1} according to the
perturbative method. The idea is to look for a change of variables
$(q,p)\mapsto(\ph,I)$ transforming the system into a similar one, but
without the term $\eps\sin q_n$. Let us write 
\begin{equation}
\label{sm3}
\begin{split}
q &= \ph + f(\ph,I) \\
p &= I + g(\ph,I),
\end{split}
\end{equation}
where $f$ and $g$ are unknown functions, which are $2\pi$-periodic in $\ph$.
We impose that this change of variables transforms the map \eqref{sm1} into
the map
\begin{equation}
\label{sm4}
\begin{split}
\ph_{n+1} &= \ph_n + I_{n+1} 	\qquad\pmod{2\pi} \\
I_{n+1} &= I_n \bydef \w.
\end{split}
\end{equation}
This is equivalent to requiring that $f$ and $g$ solve the functional
equations
\begin{equation}
\label{sm5}
\begin{split}
f(\ph+\w,\w) &= f(\ph,\w) + g(\ph+\w,\w) \\
g(\ph+\w,\w) &= g(\ph,\w) - \eps\sin(\ph+f(\ph,\w)).
\end{split}
\end{equation}
One can try to solve these equations by expanding $f$ and $g$ into Taylor
series in $\eps$ and Fourier series in $\ph$:
\begin{align}
\label{sm6a}
f(\ph,\w) &= \sum_{j=1}^\infty \eps^j f_j(\ph,\w) &
f_j(\ph,\w) &= \sum_{k=-\infty}^\infty a_{k,j}(\w) \e^{\icx k\ph} \\
\label{sm6b}
g(\ph,\w) &= \sum_{j=1}^\infty \eps^j g_j(\ph,\w) &
g_j(\ph,\w) &= \sum_{k=-\infty}^\infty b_{k,j}(\w) \e^{\icx k\ph}.
\end{align}
We will use the expansion
\begin{equation}
\label{sm7}
\sin(\ph+\eps f) = \sin\ph + \eps f_1 \cos\ph + \eps^2 \bigpar{f_2\cos\ph -
\tfrac12 f_1^2\sin\ph} + \Order{\eps^3}.
\end{equation}
At order $\eps$, we have to solve the relations
\begin{equation}
\label{sm8}
\begin{split}
f_1(\ph+\w,\w) &= f_1(\ph,\w) + g_1(\ph+\w,\w) \\
g_1(\ph+\w,\w) &= g_1(\ph,\w) - \sin\ph,
\end{split}
\end{equation}
which become, in Fourier components, 
\begin{equation}
\label{sm9}
\begin{split}
a_{k,1}\e^{\icx k\w} &= a_{k,1} + b_{k,1}\e^{\icx k\w} \\
b_{k,1}\e^{\icx k\w} &= b_{k,1} - c_{k,1},
\end{split}
\end{equation}
where $c_{k,1}$ are the Fourier components of $\sin\ph$, that is,
$c_{1,1}=-c_{-1,1}=1/(2\icx)$ and all other components vanish. We thus get 
\begin{equation}
\label{sm10}
\begin{split}
g_1(\ph,\w) &= \frac{\e^{\icx\ph}}{2\icx(1-\e^{i\w})} -
\frac{\e^{-\icx\ph}}{2\icx(1-\e^{-i\w})} 
= \frac{\cos(\ph-\w/2)}{2\sin(\w/2)}\\
f_1(\ph,\w) &= -\frac{\e^{\icx\w}\e^{\icx\ph}}{2\icx(1-\e^{i\w})^2} +
\frac{\e^{-\icx\w}\e^{-\icx\ph}}{2\icx(1-\e^{-i\w})^2}
= \frac{\sin\ph}{4\sin^2(\w/2)}. 
\end{split}
\end{equation}
Note that this is only possible for $\e^{\icx\w}\neq 1$, that is, $\w\neq
0\pmod{2\pi}$. At order $\eps^2$ we obtain similar relations as
\eqref{sm9}, but now $c_{k,2}$ denotes the Fourier coefficients of
$f_1(\ph,\w)\cos\ph$, which are nonzero for $\abs{k}=2$. Thus $g_2$ and
$f_2$ only exist if $\e^{2\icx\w}\neq 1$, or $\w\neq 0,\pi\pmod{2\pi}$.
Similarly, we will find that $g_j$ and $f_j$ only exist if
$\e^{j\icx\w}\neq 1$, so the equations \eqref{sm5} can only be solved for
irrational $\w/(2\pi)$. Even then, the expansions of $f$ and $g$ will
contain small terms of the form $1-\e^{\icx k\w}$ in the denominators, so
that the convergence of the series is not clear at all. In fact, the
convergence has been proved by Moser for certain irrational $\w$ called
\defwd{Diophantine numbers} \cite{Moser73}.

\begin{figure}
 \centerline{\psfig{figure=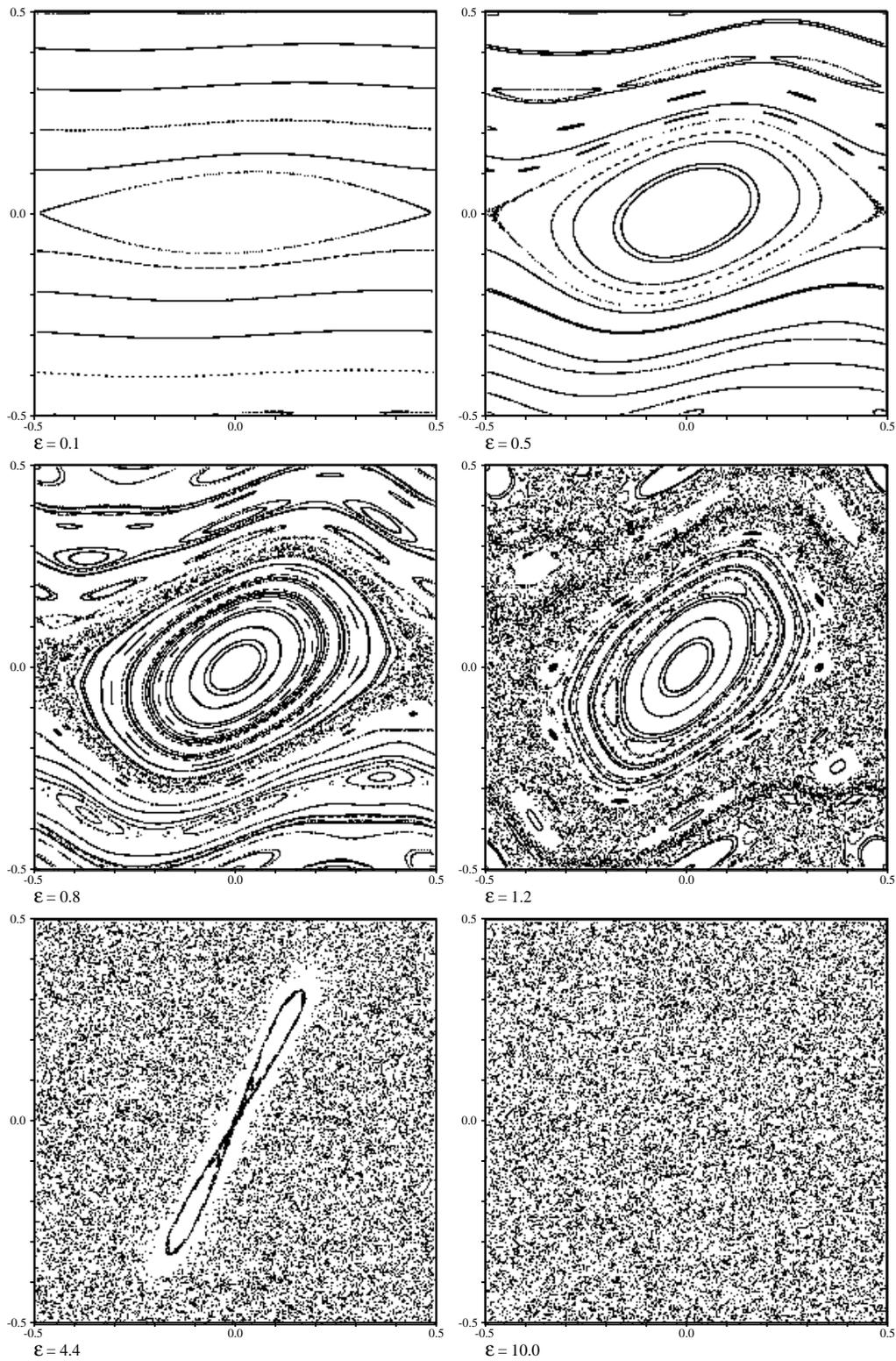,height=205mm,clip=t}}
 \vspace{3mm}
 \caption[]
 {Phase portraits of the standard map, obtained by representing several
 orbits with different initial conditions, for increasing values of the
 perturbation $\eps$. From left to right and top to bottom: $\eps=0.1$,
 $\eps=0.5$, $\eps=0.8$, $\eps=1.2$, $\eps=4.4$ and $\eps=10$.}
\label{fig_sm}
\end{figure}
Now let us turn to the geometric approach. We can consider $(q,p)$ as
coordinates in the plane (or on the cylinder because of the periodicity of
$q$). For given $(q_0,p_0)$, the set of points $\set{(q_n,p_n)}_{n\geqs0}$
is called the \defwd{orbit} with initial condition $(q_0,p_0)$. We would
like to know what the different orbits look like. The simplest case is the
\defwd{fixed point}: if
\begin{equation}
\label{sm11}
\begin{split}
q_{n+1} &= q_n \\
p_{n+1} &= p_n
\end{split}
\end{equation}
then the orbit will consist of a single point. The fixed points of the
standard map are $(0,k)$ and $(\pi,k)$ with $k\in\Z$. 
We can also have \defwd{periodic orbits}, consisting of $m$ points, if 
\begin{equation}
\label{sm12}
\begin{split}
q_{n+m} &= q_n \\
p_{n+m} &= p_n.
\end{split}
\end{equation}
Another possible orbit is the \defwd{invariant curve}. For instance if the
equations \eqref{sm5} admit a solution, we can write 
\begin{equation}
\label{sm13}
\begin{split}
q_n &= \ph_n + f(\ph_n,\w) \qquad\pmod{2\pi}\\
p_n &= \w + g(\ph_n,\w),
\end{split}
\end{equation}
where $\ph_n = \ph_0+n\w$. This is the parametric equation of a curve
winding around the cylinder. Since $\w$ is irrational, the points fill the
curve densely.
One can also analyse the dynamics in the vicinity of periodic orbits. It
turns out that for this kind of map, most periodic orbits are of one of two
types: elliptic orbits are surrounded by invariant curves, while hyperbolic
orbits attract other orbits from one direction and expel them into another
one.  There are, however, much more exotic types of orbits. Some live on
invariant Cantor sets, others densely fill regions of phase space with a
positive surface. 

The aim of the geometrical theory of dynamical systems is to classify the
possible behaviours and to find ways to determine the most important
qualitative features of the system. An important advantage is that large
classes of dynamical systems have a similar qualitative behaviour. This does
not mean that the perturbative approach is useless. But it is in general
preferable to start by analysing the system from a qualitative point of
view, and then, if necessary, use more sophisticated methods in order to
obtain more detailed information.

 
\section{The Lorenz Model}
\label{sec_lorenz}

Convection is an important mechanism in the dynamics of the atmosphere: warm
air has a lower density and therefore rises to higher altitudes, where it
cools down and falls again, giving rise to patterns in the atmospheric
currents. 

This mechanism can be modeled in the laboratory, an experiment known as
\defwd{Rayleigh-B\'enard convection}. A fluid is contained between two
horizontal plates, the upper one at temperature $T_0$ and the lower one at
temperature $T_1 = T_0+\Delta T>T_0$. The temperature difference $\Delta T$
is the \defwd{control parameter}, which can be modified. 

For small values of $\Delta T$, the fluid remains at rest, and the
temperature decreases linearly in the vertical direction. At slightly larger
$\Delta T$, convection rolls appear (their shape depends on the geometry of
the set-up). The flow is still stationary, that is, the fluid velocity at
any given point does not change in time. 

For still larger $\Delta T$, the spatial arrangement of the rolls remains
fixed, but their time dependence becomes more complex. Usually, it starts by
getting periodic. Then different scenarios are observed, depending on the
set-up. One of them is the \defwd{period doubling cascade}: the
time-dependence of the velocity field has period $P, 2P, 4P, \dots,
2^nP,\dots$, where the \nth{n} period doubling occurs for a temperature
difference $\Delta T_n$ satisfying
\begin{equation}
\label{lm1}
\lim_{n\to\infty} \frac{\Delta T_n - \Delta T_{n-1}}{\Delta T_{n+1}-\Delta
T_n} = \delta \simeq 4.4...
\end{equation}
These $\Delta T_n$ accumulate at some finite $\Delta T_\infty$ for which
the behaviour is no longer periodic, but displays temporal chaos. In this
situation, the direction of rotation of the rolls changes erratically in
time. 

For very large $\Delta T$, the behaviour can become \defwd{turbulent}: not
only is the time dependence nonperiodic, but the spatial arrangement of the
velocity field also changes. 

RB convection has been modeled in the following way. For simplicity, one
considers the two-dimensional case, with an infinite extension in the
horizontal $x_1$-direction, while the vertical $x_2$-direction is bounded
between $-\frac12$ and $\frac12$. Let $\cD=\R\times[-\frac12,\frac12]$. The
state of the system is described by three fields
\begin{align}
\nonumber
v &: \cD\times\R \to \R^2 && \text{velocity,} \\
\label{lm2}
T &: \cD\times\R \to \R   && \text{temperature,} \\
\nonumber
p &: \cD\times\R \to \R   && \text{pressure.} 
\end{align}
The deviation $\th(x,t)$ from the linear temperature profile is defined by
\begin{equation}
\label{lm3}
T(x,t) = T_0 + \Delta T \Bigpar{\frac12-x_2} + T_1 \th(x,t).
\end{equation}
The equations of hydrodynamics take the following form:
\begin{align}
\nonumber
\frac1\sigma \Bigbrak{\dpar vt + (v\cdot\nabla) v}
&= \Delta v - \nabla p + (0,\th)^T \\
\label{lm4}
\dpar \th t + (v\cdot\nabla) \th
&= \Delta \th + R v_2 \\
\nonumber
\nabla \cdot v &= 0.
\end{align}
Here $\sigma$, the \defwd{Prandtl number}, is a constant related to physical
properties of the fluid, while $R$, the \defwd{Reynolds number}, is
proportional to $\Delta T$. Furthermore,
\begin{align}
\nonumber
\nabla p &= \Bigpar{\dpar{p}{x_1}, \dpar{p}{x_2}}^T \\
\nonumber
\Delta \th &= \dpar{^2\th}{x_1^2} + \dpar{^2\th}{x_2^2} \\
\nonumber
\nabla \cdot v &= \dpar{v_1}{x_1} + \dpar{v_2}{x_2} \\
\nonumber
(v\cdot\nabla)\th &= v_1 \dpar{\th}{x_1} + v_2 \dpar{\th}{x_2}. 
\end{align}
The terms containing $(v\cdot\nabla)$ introduce the nonlinearity into the
system.  The boundary conditions require that $\th$, $v_2$ and
$\tdpar{v_1}{x_2}$ should vanish for $x_2 = \pm \frac12$. We thus have to
solve four coupled nonlinear partial differential equations for the four
fields $v_1, v_2, \th, p$. The continuity equation $\nabla\cdot v=0$ can be
satisfied by introducing the \defwd{vorticity} $\psi : \cD\times\R \to \R$,
such that
\begin{equation}
\label{lm5}
(v_1, v_2) = 
\Bigpar{-\dpar{\psi}{x_2}, \dpar{\psi}{x_1}}.
\end{equation}
It is also possible to eliminate the pressure from the two equations for
$\tdpar vt$. We are left with two equations for $\psi$ and $\th$. The
problem can be further simplified by assuming a periodic dependence on
$x_1$, of period $2\pi/q$. A possible approach (not the best one by modern
standards, but historically important) is to expand the two fields into
Fourier series (or ``modes''):
\begin{equation}
\label{lm6}
\begin{split}
\psi(x_1,x_2) &= \sum_{k\in\Z^2} a_k \e^{\icx k_1qx_1}\e^{\icx k_2\pi x_2}
\\
\th(x_1,x_2) &= \sum_{k\in\Z^2} b_k \e^{\icx k_1qx_1}\e^{\icx k_2\pi x_2}
\end{split}
\end{equation}
(where the boundary conditions impose some relations between Fourier
coefficients of the same $\abs{k_1}$ and $\abs{k_2}$). Note that the terms
of this sum are eigenfunctions of the linear operators in \eqref{lm4}.
Plugging these expansions into the equations, we obtain relations of the
form
\begin{equation}
\label{lm7}
\dtot{}{t} \Bigpar{\begin{matrix}a_k(t)\\b_k(t)\end{matrix}} 
= L_k \Bigpar{\begin{matrix}a_k(t)\\b_k(t)\end{matrix}} 
+ N\bigpar{\set{a_{k'},b_{k'}}_{k'\in\Z}},
\end{equation}
where $L_k$ are $2\times 2$ matrices and the term $N(\cdot)$ comes from the
nonlinear terms in $(v\cdot\nabla)$ and may depend on all other $k'$.
Without these nonlinear terms the problem would be easy to solve. 

\begin{figure}
 \centerline{\psfig{figure=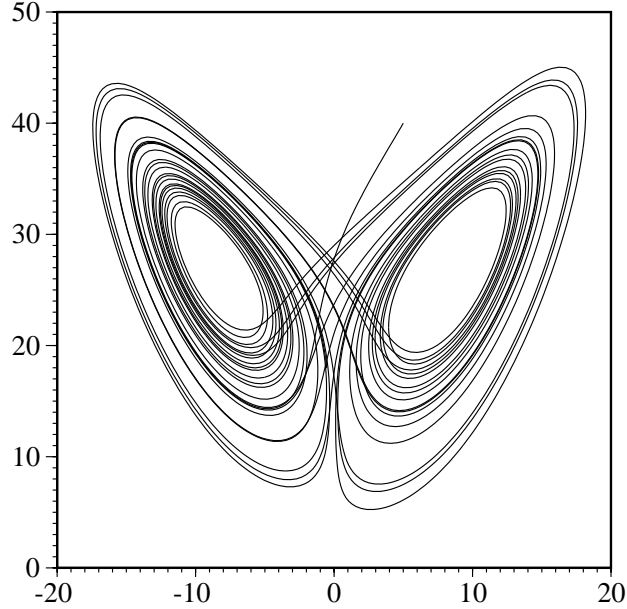,height=80mm,clip=t}}
 \vspace{3mm}
 \caption[]
 {One trajectory of the Lorenz equations \eqref{lm9} for $\sigma=10$,
 $b=8/3$ and $r=28$, projected on the $(X,Z)$-plane.}
\label{fig_lm}
\end{figure}
In 1962, Saltzmann considered approximations of the equations \eqref{lm7}
with finitely many terms, and observed that the dynamics seemed to be
dominated by three Fourier modes. In 1963, Lorenz decided to truncate the
equations to these modes \cite{Lorenz}, setting
\begin{equation}
\label{lm8}
\begin{split}
\psi(x_1,x_2) &= \alpha_1 X(t) \sin q x_1\cos\pi x_2 \\
\th(x_1,x_2) &= \alpha_2 Y(t) \cos q x_1\cos\pi x_2 + 
\alpha_3 Z(t) \sin2\pi x_2.
\end{split}
\end{equation} 
Here $\alpha_1 = \sqrt2(\pi^2+q^2)/(\pi q)$, $\alpha_2 = \sqrt2 \alpha_3$
and $\alpha_3 = (\pi^2+q^2)^3/(\pi q^2)$ are constants introduced only in
order to simplify the resulting equations. All other Fourier modes in the
expansion \eqref{lm7} are set to zero, a rather drastic approximation. After
scaling time by a factor $(\pi^2+q^2)$, one gets the equations
\begin{equation}
\label{lm9}
\begin{split}
\tdtot Xt &= \sigma (Y-X) \\
\tdtot Yt &= rX - Y - XZ \\
\tdtot Zt &= -bZ + XY,
\end{split}
\end{equation}
where $b=4\pi^2/(\pi^2+q^2)$, and $r=R q^2/(\pi^2+q^2)^3$ is proportional to
the control parameter $R$, and thus to $\Delta T$. These so-called
\defwd{Lorenz equations} are a very crude approximation of the equations
\eqref{lm4}, nevertheless they may exhibit very complicated dynamics. 

In fact, for $0\leqs r\leqs 1$, all solutions are attracted by the origin
$X=Y=Z=0$, corresponding to the fluid at rest. For $r>1$, a pair of
equilibria with $X\neq 0$ attracts the orbits, they correspond to
convection rolls with the two possible directions of rotation. Increasing
$r$ produces a very complicated sequence of bifurcations, including period
doubling cascades \cite{Sparrow}. For certain values of the parameters, a
\defwd{strange attractor} is formed, in which case the convection rolls
change their direction of rotation very erratically (\figref{fig_lm}), and
the dynamics is very sensitive to small changes in the initial conditions
(the \defwd{Butterfly effect}). The big surprise was that such a simple
approximation, containing only three modes, could capture such complex
behaviours. 


\section{The Logistic Map}
\label{sec_logistic}

Our last example is a famous map inspired by population dynamics. Consider a
population of animals that reproduce once a year. Let $P_n$ be the number of
individuals in the year number $n$. The offspring being usually proportional
to the number of adults, the simplest model for the evolution of the
population from one year to the next is the linear equation
\begin{equation}
\label{lg1}
P_{n+1} = \lambda P_n
\end{equation}
where $\lambda$ is the natality rate (minus the mortality rate). This law
leads to an exponential growth of the form
\begin{equation}
\label{lg2}
P_n = \lambda^n P_0 = \e^{n\ln\lambda} P_0
\end{equation}
(the \defwd{Malthus law}). This model becomes unrealistic when the number of
individuals is so large that the limitation of resources becomes apparent.
The simplest possibility to limit the growth is to introduce a quadratic
term $-\beta P_n^2$, leading to the law
\begin{equation}
\label{lg3}
P_{n+1} = \lambda P_n - \beta P_n^2.
\end{equation}
The rescaled variable $x = \beta P$ then obeys the equation
\begin{equation}
\label{lg4}
x_{n+1} = f_\lambda(x_n) \defby \lambda x_n (1-x_n).
\end{equation}
The map $f_\lambda$ is called the \defwd{logistic map}. 
Observe that for $0\leqs\lambda\leqs 4$, $f_\lambda$ maps the interval
$[0,1]$ into itself. The dynamics of the sequence $x_n$ depends drastically
on the value of $\lambda$.

For $0\leqs\lambda\leqs 1$, all orbits converge to $0$, which means that the
population becomes extinct. For $1<\lambda\leqs 3$, all orbits starting at
$x_0>0$ converge to $1 - 1/\lambda$, and thus the population reaches a
stable equilibrium. For $3 < \lambda \leqs 1+\sqrt6$, the orbits converge to
a cycle of period $2$, so that the population asymptotically jumps back and
forth between two values. 

For $\lambda > 1+\sqrt6$, the system goes through a whole seqence of period
doublings. Similarly as in RB convection, the values $\lambda_n$ of the
parameter for which the \nth{n} period doubling occurs obey the law
\begin{equation}
\label{lg5}
\lim_{n\to\infty} \frac{\lambda_n - \lambda_{n-1}}{\lambda_{n+1}-\lambda_n} 
= \delta = 4.669...
\end{equation}
where $\delta$ is called the \defwd{Feigenbaum constant}. In 1978,
Feigenbaum as well as Coullet and Tresser independently outlined an
argument showing that such period doubling cascades should be observable
for a large class of systems, and that the constant $\delta$ is universal.
For instance, it also appears in the two-dimensional \defwd{H\'enon map}
\begin{equation}
\label{lg6}
\begin{split}
x_{n+1} &= 1 - \lambda x_n^2 + y_n \\
y_{n+1} &= b x_n	
\qquad\qquad\qquad 0 < b < 1
\end{split}
\end{equation}
Rigorous proofs of these properties were later worked out by Collet,
Eckmann, Koch, Lanford and others \cite{CE}.

\begin{figure}
 \centerline{\psfig{figure=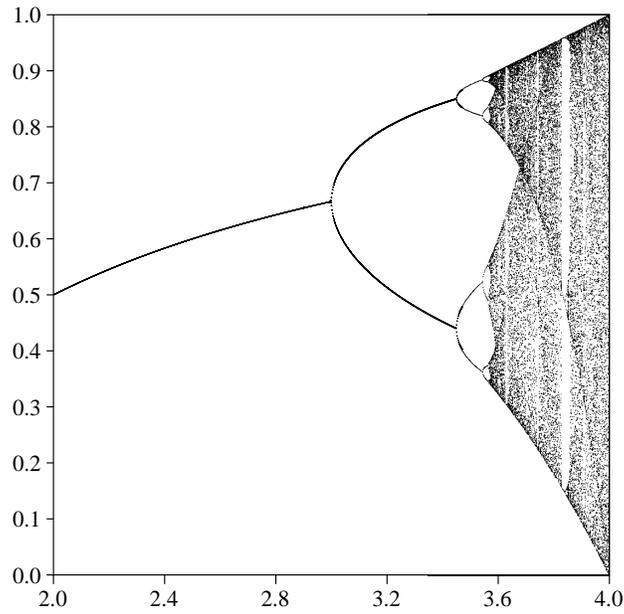,height=80mm,clip=t}}
 \vspace{3mm}
 \caption[]
 {Bifurcation diagram of the logistic map. For each value of $\lambda$ on
 the abscissa, the points $x_{1001}$ through $x_{1100}$ are represented on
 the ordinate, for an initial value $x_0=\frac12$.}
\label{fig_logi}
\end{figure}

For the logistic map, the period doublings accumulate at a value
$\lambda_\infty = 3.56...$, beyond which the orbits become chaotic. For
larger $\lambda$, there is a complicated interplay of regular and chaotic
motion (\figref{fig_logi}), containing other period doubling cascades.
Finally, for $\lambda=4$, one can prove that the dynamics is as random as
coin tossing.  


\chapter{Stationary and Periodic Solutions}
\label{ch_sp}

In Chapter \ref{ch_ex}, we have seen examples of two kinds of dynamical
systems: ordinary differential equations (ODEs) and iterated maps. There are
other types of dynamical systems, such as partial differential equations or
cellular automata. These are in general more difficult to analyse, although
some ideas developed for maps and ODEs can be carried over to their study.
Here we will concentrate on ODEs and maps, by starting with the simplest
kinds of dynamics: stationary and periodic. 


\section{Basic Concepts}
\label{sec_basic}


\subsection{Orbits and Flows}
\label{ssec_orbit}

Let $\cD\subset\R^n$ be an open domain. One type of dynamical systems we
will consider is given by a \defwd{map} $F: \cD\to\cD$. $\cD$ is called the
\defwd{phase space} of the system. It is possible to consider more general
differentiable manifolds as phase space, such as the circle, the cylinder
or the torus, but we will limit the discussion to Euclidean domains.
Generalizations to other manifolds are usually straightforward.

\begin{definition}
\label{def_orbit1}
The \defwd{(positive) orbit} of $F$ through a point $x_0\in\cD$ is the
sequence $(x_k)_{k\geqs 0}$ defined by $x_{k+1}=F(x_k)$ for all integers
$k\geqs 0$. We have thus
\begin{equation}
\label{orb1}
x_k = F^k(x)
\qquad
\text{where $F^k = \underbrace{F\circ F\circ \dots\circ F}_{\text{$k$ times}}$}
\end{equation}
In case $F$ is invertible, we can also define the \defwd{negative orbit} of
$x_0$ by the relations $F(x_k)=x_{k+1}$ for all $k<0$, which are equivalent
to \eqref{orb1} if we set $F^{-k} = (F^{-1})^k$ for all $k>0$. The
\defwd{orbit} of $x_0$ is then given by $(x_k)_{k\in\Z}$. 
\end{definition}

Note the trivial relation
\begin{equation}
\label{orb2}
F^{k+l}(x_0) = F^k(F^l(x_0))
\end{equation}
for all positive integers $k, l$ (and all integers if $F$ is invertible),
which will admit an analogue in the case of ODEs. 

As particular cases of maps $F$, we have \defwd{homeomorphisms}, which are
continuous maps admitting a continuous inverse, and \defwd{diffeomorphisms},
which are continuously differentiable maps admitting a continuously
differentiable inverse. Similarly, for all $r\geqs 1$, a
\defwd{$\cC^r$-diffeomorphism} is an invertible map $F$ such that both $F$
and $F^{-1}$ admit continuous derivatives up to order $r$.

The \defwd{ordinary differential equations} we are going to consider are of
the form
\begin{equation}
\label{orb3}
\dot{x} = f(x),
\end{equation}
where $f: \cD\to\R^n$, and $\dot{x}$ denotes $\dtot xt$. Equivalently, we
can write \eqref{orb3} as a system of equations for the components of $x$,
\begin{equation}
\label{orb4}
\dot{x}_i = f_i(x), \qquad
i = 1,\dots, n.
\end{equation}
$\cD\subset\R^n$ is again called \defwd{phase space} and $f$ is called a
\defwd{vector field}. 

To define the orbits of $f$, we have to treat the problem of existence and
uniqueness a bit more carefully. The following results are assumed to be
known from basic analysis (see for instance \cite{Hale}, \cite{Hartman} or
\cite{HS}).

\begin{theorem}[Peano-Cauchy]
\label{thm_PeanoCauchy}
Let $f$ be continuous. For every $x_0\in\cD$, there exists at least one
local solution of \eqref{orb3} through $x_0$, that is, there is an open
interval $I\ni 0$ and a function $x:I\to\cD$ such that $x(0)=x_0$ and
$\dot{x}(t) = f(x(t))$ for all $t\in I$. 
\end{theorem}

\begin{theorem}
\label{thm_maxinterval}
Every solution $x(t)$ with $x(0)=x_0$ can be continued to a maximal
interval of existence $(t_1,t_2)\ni 0$. If $t_2<\infty$ or $t_1 >
-\infty$,  then for any compact $\cK\subset\cD$, there exists a time
$t\in(t_1,t_2)$ with $x(t)\notin\cK$ (this means that solutions will
diverge or reach $\partial\cD$).
\end{theorem}

\begin{theorem}[Picard-Lindel\"of]
\label{thm_PicardLindelof}
Assume $f$ is continuous and \defwd{locally Lipschitzian}, that is, for
every compact $\cK\subset\cD$, there exists a constant $L_\cK$ such that
$\norm{f(x)-f(y)} \leqs L_\cK \norm{x-y}$ for all $x,y\in\cK$. Then there is
a unique solution $x(t)$ of \eqref{orb3} with $x(0)=x_0$ for every
$x_0\in\cD$. 
\end{theorem}

Note in particular that if $f$ is continuously differentiable, then it is
locally Lipschitzian. We will usually consider vector fields which are at
least once continuously differentiable. 

\begin{example}
\label{ex_orbit1}
It is easy to give counterexamples to global existence and uniqueness. For
instance,
\begin{equation}
\label{orb5}
\dot x = x^2 
\qquad\Rightarrow\qquad
x(t) = \frac1{\frac1{x_0}-t}
\end{equation}
has a solution diverging for $t=\frac1{x_0}$. 
A physically interesting counterexample to uniqueness is the \defwd{leaky
bucket equation}
\begin{equation}
\label{orb6}
\dot x = - \sqrt{\abs{x}}.
\end{equation}
Here $x$ is proportional to the height of water in a bucket with a hole in
the bottom, and \eqref{orb6} reflects the fact that the kinetic energy
(proportional to $\dot{x}^2$) of the water leaving the bucket is equal to
the potential energy of the water inside. For every $c$, \eqref{orb6}
admits the solution
\begin{equation}
\label{orb7}
x(t) = 
\begin{cases}
\frac14(t-c)^2 & \text{for $t<c$} \\
0 & \text{for $t\geqs c$.}
\end{cases}
\end{equation}
In particular, for any $c\leqs 0$, \eqref{orb7} is a solution of
\eqref{orb6} such that $x(0)=0$.  This reflects the fact that if the bucket
is empty at time $0$, we do not know at what time it was full.
\end{example}

For simplicity, we will henceforth assume that the ODE \eqref{orb3} admits a
unique global solution for all $x_0\in\cD$. This allows to introduce the
following definitions:\footnote{In the case of $x(t)$ existing for all
positive $t$ but not necessarily for all negative $t$, the definition
remains valid with orbit replaced by \defwd{positive orbit}, flow replaced
by \defwd{semi-flow} and group replaced by \defwd{semi-group}.}

\begin{definition}
\label{def_orbit2}
Let $x_0\in\cD$ and let $x(t)$ be the unique solution of \eqref{orb3} with
\defwd{initial condition} $x(0)=x_0$. 
\begin{itemiz}
\item	The \defwd{integral curve} through $x_0$ is the set
	$\setsuch{(x,t)\in\cD\times\R}{x=x(t)}$.
\item 	The \defwd{orbit} through $x_0$ is the set
	$\setsuch{x\in\cD}{x=x(t),t\in\R}$.
\item	The \defwd{flow} of the equation \eqref{orb3} is the map
\begin{equation}
\label{orb8}
\fctndef{\ph}{\cD\times\R}{\cD}
{(x_0,t)}{\ph_t(x_0)=x(t)}
\end{equation}
\end{itemiz}
\end{definition}

Geometrically speaking, the orbit is a curve in phase space containing $x_0$
such that the vector field $f(x)$ is tangent to the curve at any point $x$
of the curve. Uniqueness of the solution means that there is only one orbit
through any point in phase space.

By definition, we have $\ph_0(x_0)=x_0$ for all $x_0\in\cD$, and uniqueness
of solutions implies that $\ph_t(\ph_s(x_0)) = \ph_{t+s}(x_0)$. These
properties can be rewritten as 
\begin{equation}
\label{orb9}
\ph_0 = \identity
\qquad
\ph_t \circ \ph_s = \ph_{t+s}
\end{equation}
which means that the family $\set{\ph_t}_t$ forms a group. Note the
similarity between this relation and the relation \eqref{orb2} for iterated
maps. 

\begin{example}
\label{ex_orbit2}
In the case $f(x)=-x$, $x\in\R$, we have 
\begin{equation}
\label{orb10}
\ph_t(x_0) = x_0 \e^{-t}.
\end{equation}
The system admits three distinct orbits $(0,\infty)$, $(-\infty,0)$ and
$\set{0}$. 
\end{example}


\subsection{Evolution of Volumes}
\label{ssec_vol}

Let $\cM\subset\cD$ be a compact subset of phase space. We can define its
\defwd{volume} by a usual Riemann integral:
\begin{equation}
\label{vol1}
\volume(\cM) = \int_\cM \6x, 
\qquad \6x = \6x_1\dots\6x_n.
\end{equation}
The set $\cM$ will evolve under the influence of the dynamics: we can define
the sets $\cM_k = F^k(\cM)$ or $\cM(t) = \ph_t(\cM)$. How does their volume
evolve with time? The answer is actually quite simple. 

Consider first the case of a map $F$. We assume that $F$ is continuously
differentiable and denote by 
\begin{equation}
\label{vol2}
\dpar Fx(x)
\end{equation}
the \defwd{Jacobian matrix} of $F$ at $x$, which is the $n\times n$ matrix
$A$ with elements $a_{ij} = \dpar{F_i}{x_j}(x)$. 

\goodbreak
\begin{prop}
\label{prop_vol1}
Assume $F$ is a diffeomorphism and let $V_k = \volume(\cM_k)$.
Then 
\begin{equation}
\label{vol3}
V_{k+1} = \int_{\cM_k}\Bigabs{\det\Bigpar{\dpar Fx(x)}} \6x.
\end{equation}
\end{prop} 
\begin{proof}
This is a simple application of the formula for a change of variables in an
integral:
\[
V_{k+1} = \int_{\cM_{k+1}} \6y = \int_{\cM_k} \Bigabs{\det\Bigpar{\dpar
yx(x)}} \6x = \int_{\cM_k}\Bigabs{\det\Bigpar{\dpar Fx(x)}} \6x.
\]
\end{proof}

\begin{definition}
\label{def_vol1}
The map $F$ is called \defwd{conservative} if 
\begin{equation}
\label{vol4}
\Bigabs{\det\Bigpar{\dpar Fx(x)}} = 1 
\qquad \forall x\in\cD
\end{equation}
The map $F$ is called \defwd{dissipative} if 
\begin{equation}
\label{vol5}
\Bigabs{\det\Bigpar{\dpar Fx(x)}} < 1 
\qquad \forall x\in\cD.
\end{equation}
\end{definition}

Proposition \ref{prop_vol1} implies that $V_{k+1}=V_k$ if $F$ is
conservative and $V_{k+1}<V_k$ if $F$ is dissipative. More generally, if
$\abs{\det\dpar Fx(x)}\leqs \lambda$ for some constant $\lambda$ and all
$x\in\cD$, then $V_k\leqs\lambda^kV_0$. 

\begin{figure}
 \centerline{\psfig{figure=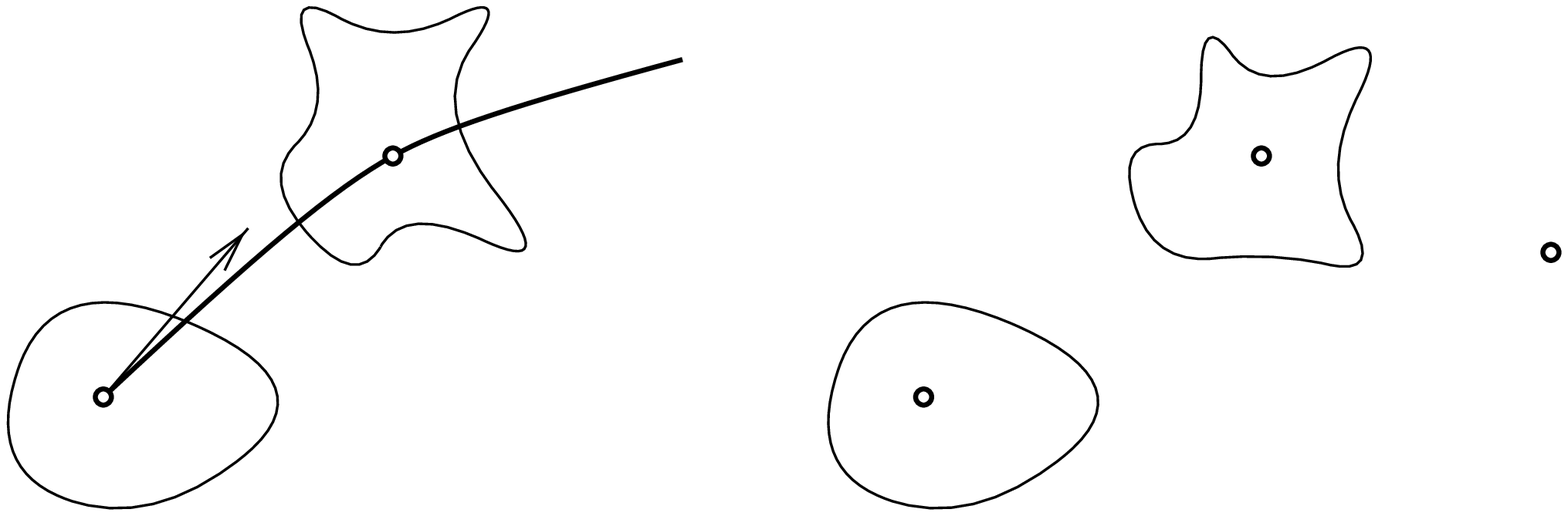,width=100mm,clip=t}}
 \figtext{
 	\writefig	1.9	3.5	a
 	\writefig	8.0	3.5	b
 	\writefig	3.2	1.0	$x(0)$
 	\writefig	2.5	2.2	$f(x(0))$
 	\writefig	4.2	1.0	$\cM(0)$
 	\writefig	5.7	2.8	$x(t)$
 	\writefig	5.8	2.5	$=\ph_{t}(x(0))$
 	\writefig	5.7	1.8	$\cM(t)$
 	\writefig	8.45	1.05	$x_1$
 	\writefig	9.5	1.1	$\cM_1$
 	\writefig	10.5	2.5	$x_2$
 	\writefig	10.9	1.75	$\cM_2$
 	\writefig	12.35	1.9	$x_3$
 }
 \captionspace
 \caption[Evolution of volumes in phase space]
 {(a) Evolution of a volume with the flow, (b) evolution with an
 iterated map.}
\label{fig_flowmap}
\end{figure}

For differential equations, the result is the following:

\begin{prop}
\label{prop_vol2}
Assume $f$ is continuously differentiable and let $V(t) = \volume(\cM(t))$.
Then 
\begin{equation}
\label{vol6}
\dtot{}{t} V(t) = \int_{\cM(t)}\nabla\cdot f(x) \6x,
\end{equation}
where $\nabla\cdot f = \sum_{i=1}^n \dpar{f_i}{x_i}$ is the divergence. 
\end{prop}
\begin{proof}
We have 
\[
V(t) = \int_{\cM(t)}\6y = \int_\cM \Bigabs{\det \dpar{}{x}\ph_t(x)} \6x.
\]
Let us fix $x\in\cM$, let $y(t)=\ph_t(x)$ and set
\[
J(t) \defby \dpar{}{x}\ph_t(x), 
\qquad
A(t) \defby \dpar fx(y(t)).
\]
Note that by definition of $\ph_t$, $J(0)=\one$ is the identity matrix. Now
we can compute
\[
\dtot{}{t} J(t) = \dpar{}{x}\dot{y}(t) = \dpar{}{x} f(\ph_t(x)) 
= \dpar fx(\ph_t(x)) \dpar{}{x}\ph_t(x),
\]
and thus
\[
\dtot{}t J(t) = A(t)J(t), 
\qquad J(0) = \one.
\]
This is a linear, time-dependent differential equation for $J(t)$, which is
known to admit a unique global solution. This implies in particular that
$\det J(t)\neq 0$ $\forall t$, since otherwise $J(t)$ would not be
surjective, contradicting uniqueness. Since $\det J(0)=1$, continuity
implies that $\det J(t)>0$ $\forall t$. Now let us determine the evolution
of $\det J(t)$. By Taylor's formula, there exists $\th\in[0,1]$ such that
\[
\begin{split}
J(t+\eps) &= J(t) + \eps\dtot{}{t} J(t+\th\eps) \\
&= J(t) \bigbrak{\one + \eps J(t)^{-1} A(t+\th\eps)J(t+\th\eps)}.
\end{split}
\]
From linear algebra, we know that for any $n\times n$ matrix $B$,
\[
\det(\one + \eps B) = 1 + \eps\Tr B + r(\eps)
\]
with $\lim_{\eps\to0}r(\eps)/\eps = 0$ (this is a consequence of the
definition of the determinant as a sum over permutations). Using
$\Tr(AB)=\Tr(BA)$, this leads to
\[
\det J(t+\eps) = \det J(t) \bigbrak{1 + \eps
\Tr\bigpar{A(t+\th\eps)J(t+\th\eps)J(t)^{-1}} + r(\eps)},
\]
and thus
\[
\dtot{}{t} \det J(t) 
= \lim_{\eps\to0} \frac{\det J(t+\eps) - \det J(t)}{\eps} 
= \Tr(A(t)) \det J(t). 
\]
Taking the derivative of $V(t)$ we get 
\[
\begin{split}
\dtot{}{t} V(t) 
&= \int_\cM \dtot{}{t} \det J(t) \6x \\
&= \int_\cM \Tr \Bigpar{\dpar fx(y(t))} \det J(t) \6x 
\vrule height 20pt depth 0pt width 0pt\\
&= \int_{\cM(t)} \Tr \Bigpar{\dpar fx(y)} \6y,
\vrule height 20pt depth 0pt width 0pt
\end{split}
\]
and the conclusion follows from the fact that $\Tr\dpar fx = \nabla\cdot f$.
\end{proof}

\begin{definition}
\label{def_vol2}
The vector field $f$ is called \defwd{conservative} if 
\begin{equation}
\label{vol6b}
\nabla \cdot f(x) = 0 
\qquad \forall x\in\cD
\end{equation}
The vector field $f$ is called \defwd{dissipative} if 
\begin{equation}
\label{vol7}
\nabla \cdot f(x) < 0  
\qquad \forall x\in\cD.
\end{equation}
\end{definition}

Proposition \ref{prop_vol2} implies that $V(t)$ is constant if $f$ is
conservative, and monotonously decreasing when $f$ is dissipative.  More
generally, if $\nabla \cdot f(x) \leqs c$ $\forall x\in\cD$, then $V(t)
\leqs V(0)\e^{ct}$. 

Of course, one can easily write down dynamical systems which are neither
conservative nor dissipative, but the conservative and dissipative
situations are very common in applications.

\begin{example}
\label{ex_vol}
Consider a Hamiltonian system, with Hamiltonian $H\in\cC^2(\R^{2m},\R)$.
Then $x=(q,p)\in \R^m\times\R^m$ and the equations \eqref{i4} take the form
\begin{equation}
\label{vol8}
f_i(x) = 
\begin{cases}
\displaystyle\dpar H{p_i} & i=1,\dots m \\
\vrule height 20pt depth 0pt width 0pt
-\displaystyle\dpar H{q_{i-m}} & i=m+1,\dots 2m.
\end{cases}
\end{equation}
This implies that
\begin{equation}
\label{vol9}
\nabla\cdot f = \sum_{i=1}^m \dpar{}{q_i} \Bigpar{\dpar H{p_i}} 
+ \sum_{i=1}^m \dpar{}{p_i} \Bigpar{-\dpar H{q_i}} = 0.
\end{equation}
Thus all (sufficiently smooth) Hamiltonian systems are conservative.
\end{example}

\begin{exercise}
\label{exo_vol}
Determine whether the following systems are conservative, dissipative, or
none of the above:
\begin{itemiz}
\item	the standard map \eqref{sm1};
\item	the Lorenz model \eqref{lm9};
\item	the logistic map \eqref{lg4};
\item	the H\'enon map \eqref{lg6}.
\end{itemiz}
\end{exercise}


\section{Stationary Solutions}
\label{sec_stat}

A \defwd{stationary solution} of a dynamical system is a solution that does
not change in time. We thus define

\begin{definition}\hfill
\label{def_stat}
\begin{itemiz}
\item	A \defwd{fixed point} of the map $F$ is a point $x^\star\in\cD$ such
that 
\begin{equation}
\label{stat1}
F(x^\star) = x^\star.
\end{equation}
\item	A \defwd{singular point} of the vector field $f$ is a point 
$x^\star\in\cD$ such that 
\begin{equation}
\label{stat2}
f(x^\star) = 0.
\end{equation}
\end{itemiz}
\end{definition}

In both cases, $x^\star$ is also called \defwd{equilibrium point}. Its orbit
is simply $\set{x^\star}$ and is called a \defwd{stationary orbit}. Note
that a singular point of $f$ is also a fixed point of the flow, and
therefore sometimes abusively called a ``fixed point of $f$''.

We are now interested in the behaviour near an equilibrium point. In this
section, we will always assume that $f$ and $F$ are twice continuously
differentiable. If $x^\star$ is a singular point of $f$, the change of
variables $x = x^\star + y$ leads to the equation
\begin{equation}
\label{stat3}
\begin{split}
\dot{y} &= f(x^\star+y)\\
&= Ay + g(y),
\end{split}
\end{equation}
where we have introduced the Jacobian matrix
\begin{equation}
\label{stat4}
A = \dpar fx(x^\star).
\end{equation}
Taylor's formula implies that there exists a neighbourhood $\cN$ of $0$ and
a constant $M>0$ such that
\begin{equation}
\label{stat5}
\norm{g(y)} \leqs M \norm{y}^2 
\qquad \forall y \in \cN.
\end{equation}
Similarly, the change of variables $x_k = x^\star + y_k$ transforms an
iterated map into
\begin{equation}
\label{stat6}
\begin{split}
y_{k+1} &= F(x^\star+y_k) - x^\star\\
&= By_k + G(y_k),
\end{split}
\end{equation}
where
\begin{equation}
\label{stat7}
B = \dpar Fx(x^\star), 
\qquad
\norm{G(y)} \leqs M \norm{y}^2 
\quad \forall y \in \cN.
\end{equation}


\subsection{Linear Case}
\label{ssec_slin}

Let us start by analysing the equations \eqref{stat3} and \eqref{stat6} in
the linear case, that is without the terms $g(x)$ and $G(x)$. Consider first
the ODE
\begin{equation}
\label{slin1}
\dot y = Ay.
\end{equation}
The solution can be written as 
\begin{equation}
\label{slin2}
y(t) = \e^{At} y(0),
\end{equation}
where the \defwd{exponential of $A$} is defined by the absolutely convergent
series
\begin{equation}
\label{slin3}
\e^{At} \equiv \exp(At) \defby \sum_{k=0}^\infty \frac{t^k}{k!} A^k.
\end{equation}
In order to understand the behaviour of $\e^{At}$, let us recall some facts
from linear algebra. We can write the characteristic polynomial of $A$ as
\begin{equation}
\label{slin4}
c_A(\lambda) = \det(\lambda\one-A) = \prod_{j=1}^m (\lambda-a_j)^{m_j},
\end{equation}
where $a_1,\dots,a_m\in\C$ are distinct eigenvalues of $A$, and $m_j$ are
their \defwd{algebraic multiplicities}. The \defwd{geometric multiplicity}
$g_j$ of $a_j$ is defined as the number of independent eigenvectors
associated with $a_j$, and satisfies $1\leqs g_j\leqs m_j$. 

The results on decomposition of matrices leading to the Jordan canonical
form can be formulated as follows \cite{HS}. The matrix $A$ can be
decomposed as 
\begin{equation}
\label{slin5}
A = S + N, 
\qquad
SN = NS.
\end{equation}
Here $S$, the \defwd{semisimple part}, can be written as 
\begin{equation}
\label{slin6}
S = \sum_{j=1}^m a_j P_j,
\end{equation}
where the $P_j$ are projectors on the eigenspaces of $A$, satisfying
$P_jP_k = \delta_{jk}P_j$, $\sum_j P_j=\one$ and $m_j=\dim(P_j\R^n)$. 
The \defwd{nilpotent part} $N$ can be written as 
\begin{equation}
\label{slin7}
N = \sum_{j=1}^m N_j,
\end{equation} 
where the $N_j$ satisfy the relations
\begin{equation}
\label{slin8}
N_j^{m_j} = 0, \qquad
N_jN_k = 0 \quad\text{for $j\neq k$,} \qquad
P_jN_k = N_kP_j = \delta_{jk}N_j.
\end{equation}
In an appropriate basis, each $N_j$ is block-diagonal, with $g_j$ blocks of
the form
\begin{equation}
\label{slin9}
\begin{pmatrix}
0 & 1 &  & 0 \\
  & \ddots & \ddots & \\
  & & \ddots & 1  \\
0 & & & 0
\end{pmatrix}.
\end{equation}
In fact, $N_j=0$ unless $g_j < m_j$. 

\begin{lemma}
\label{lem_slin1}
With the above notations
\begin{equation}
\label{slin10}
\e^{At} = \sum_{j=1}^m \e^{a_jt} P_j 
\Bigpar{\one + N_j t + \dots + \frac1{(m_j-1)!}N_j^{m_j-1} t^{m_j-1}}
\end{equation}
\end{lemma}
\begin{proof}
We use the fact that $e^{At}\e^{Bt}=\e^{(A+B)t}$ whenever $AB=BA$, which can
be checked by a direct calculation. Then $\e^{At}=\e^{St}\e^{Nt}$ with 
\[
\begin{split}
\e^{St} &= \prod_{j=1}^m\e^{a_jP_jt} 
= \prod_{j=1}^m\bigl(\one+(\e^{a_jt}-1)P_j\bigr) 
= \one + \sum_{j=1}^m(\e^{a_jt}-1)P_j = \sum_{j=1}^m\e^{a_jt}P_j,\\
\e^{Nt} &= \prod_{j=1}^m\e^{N_jt} = \one + \sum_{j=1}^m(\e^{N_jt}-\one).
\end{split}
\]
The result follows from the facts that $P_j(\e^{N_kt}-\one)=0$ for $j\neq
k$, and that $\e^{N_jt}$ contains only finitely many terms, being
nilpotent.
\end{proof}

The expression \eqref{slin10} shows that the long-time behaviour is
determined by the real parts of the eigenvalues $a_j$, while the nilpotent
terms, when present, influence the short time behaviour. This motivates the
following terminology:

\begin{definition}
\label{def_slin1}
The \defwd{unstable}, \defwd{stable} and \defwd{center subspace} of the
singular point $x^\star$ are defined, respectively, by
\begin{align}
\nonumber
E^+ &\defby P^+\R^n 
= \bigsetsuch{y}{\lim_{t\to-\infty}\e^{At}y = 0}, 
& P^+ &\defby \sum_{j:\re a_j>0}P_j,
\\
\label{slin11}
E^- &\defby P^-\R^n 
= \bigsetsuch{y}{\lim_{t\to+\infty}\e^{At}y = 0},  
& P^- &\defby \sum_{j:\re a_j<0}P_j,
\\
\nonumber
E^0 &\defby P^0\R^n,  
& P^0 &\defby \sum_{j:\re a_j=0}P_j.
\end{align}
The subspaces are invariant subspaces of $\e^{At}$, that is,
$\e^{At}E^+\subset E^+$, $\e^{At}E^-\subset E^-$ and $\e^{At}E^0\subset
E^0$. 
The fixed point is called
\begin{itemiz}
\item	a \defwd{sink} if $E^+ = E^0 = \set{0}$,
\item	a \defwd{source} if $E^- = E^0 = \set{0}$,
\item	a \defwd{hyperbolic point} if $E^0 = \set{0}$,
\item	an \defwd{elliptic point} if $E^+ = E^- = \set{0}$.
\end{itemiz}
\end{definition}

\begin{figure}
 \centerline{\psfig{figure=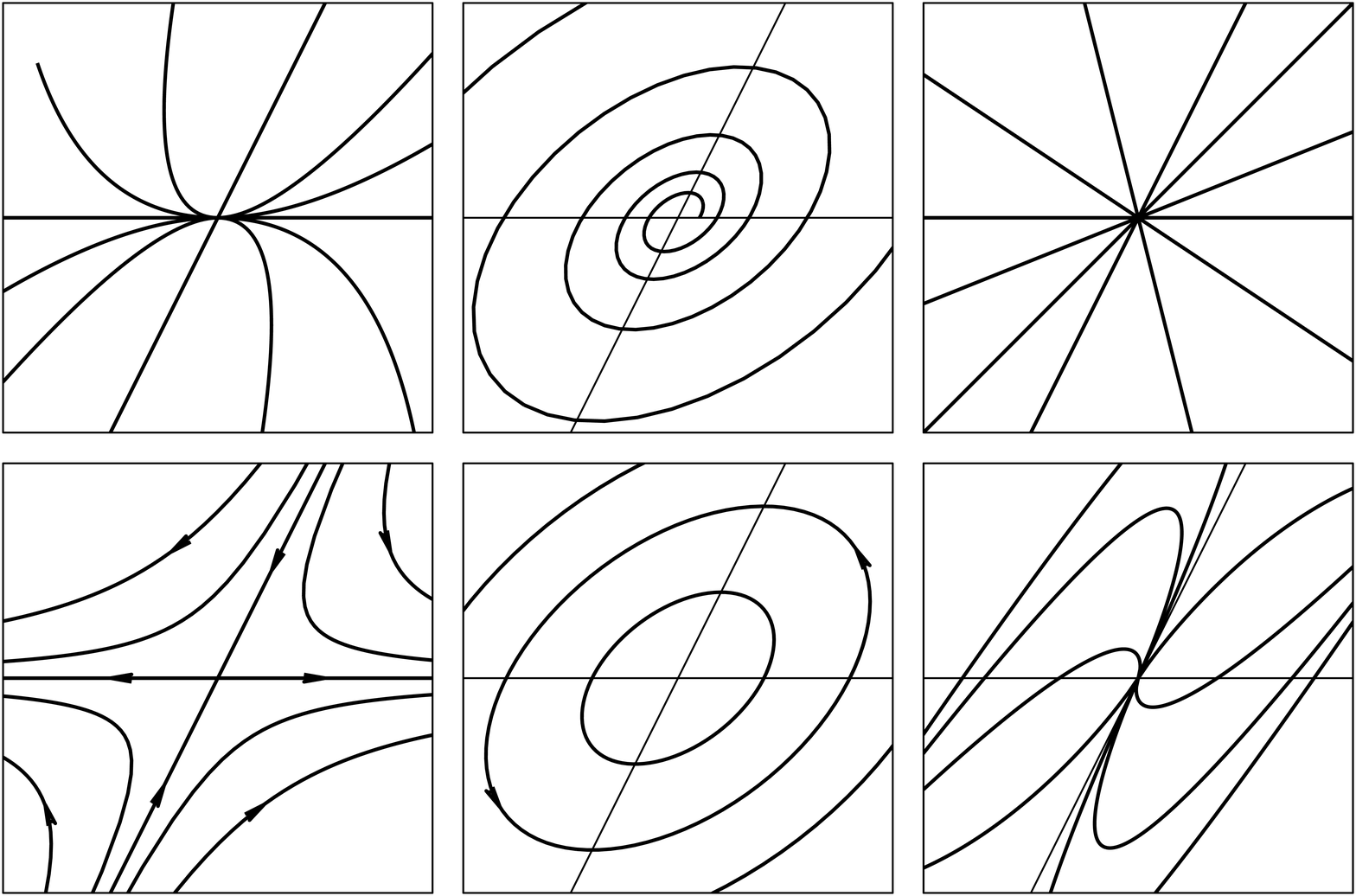,height=60mm,clip=t}}
 \figtext{
 	\writefig	3.0	6.2	a
 	\writefig	3.0	3.0	b
 	\writefig	6.1	6.2	c
 	\writefig	6.1	3.0	d
 	\writefig	9.2	6.2	e
 	\writefig	9.2	3.0	f
 }
 \captionspace
 \caption[Phase portraits of a linear two--dimensional system]
 {Phase portraits of a linear two--dimensional system: (a) node,
 (b) saddle, (c) focus, (d) center, (e) degenerate node, (f) improper node.}
\label{fig_slin}
\end{figure}

\begin{example}
\label{ex_slin}
Let $n=2$, and let $A$ be in Jordan canonical form, with $\det A\neq 0$. 
Then  we can distinguish between the following behaviours, depending on the
eigenvalues $a_1, a_2$ of $A$ (see \figref{fig_slin}).
\begin{enum}
\item	$a_1\neq a_2$
	\begin{enum}
	\item	If $a_1, a_2\in\R$, then $A = \bigl(
		\begin{smallmatrix} a_1 & 0 \\ 0 & a _2\end{smallmatrix} \bigr)$ 
		and
		\[
		\e^{At} = 
		\begin{pmatrix}
		\e^{a_1t} & 0 \\ 0 & \e^{a_2t}
		\end{pmatrix}
		\qquad
		\Rightarrow
		\qquad
		\begin{matrix}
		y_1(t) = \e^{a_1t} y_1(0) \\
		y_2(t) = \e^{a_2t} y_2(0)
		\end{matrix}
		\]
		The orbits are curves of the form $y_2 = c y_1^{a_2/a_1}$. 
		$x^\star$ is called a \defwd{node} if $a_1a_2>0$,
		and a \defwd{saddle} if $a_1a_2<0$.
	\item	If $a_1=\cc{a_2}=a+\icx\w\in\C$, then the real canonical
		form of $A$ is $A = \bigl(\begin{smallmatrix} 
		a & -\w \\ \w & a\end{smallmatrix} \bigr)$
		and 
		\[
		\e^{At} = \e^{at}
		\begin{pmatrix}
		\cos\w t & -\sin\w t \\ \sin\w t & \cos\w t
		\end{pmatrix}
		\qquad
		\Rightarrow
		\qquad
		\begin{matrix}
		y_1(t) = \e^{at} (y_1(0) \cos\w t - y_2(0) \sin\w t) \\
		y_2(t) = \e^{at} (y_1(0) \sin\w t + y_2(0) \cos\w t)
		\end{matrix}
		\]
		$x^\star$ is called a \defwd{focus} if $a\neq 0$, and a
		\defwd{center} if $a=0$. The orbits are spirals or ellipses.
	\end{enum}
\item	$a_1=a_2\bydef a$
	\begin{enum}
	\item	If $a$ has geometric multiplicity $2$, then $A = a\one$ and
		$\e^{At}=\e^{at}\one$; $x^\star$ is called a
		\defwd{degenerate node}.
	\item	If $a$ has geometric multiplicity $1$, then $A = \bigl(
		\begin{smallmatrix} a & 1 \\ 0 & a \end{smallmatrix} \bigr)$ 
		and 
		\[
		\e^{At} = \e^{at}
		\begin{pmatrix}
		1 & t \\ 0 & 1
		\end{pmatrix}
		\qquad
		\Rightarrow
		\qquad
		\begin{matrix}
		y_1(t) = \e^{at} (y_1(0) + y_2(0) t) \\
		y_2(t) = \e^{at}  y_2(0) \phantom{({}+ y_2(0))}
		\end{matrix}
		\]
		$x^\star$ is called an \defwd{improper node}.
	\end{enum}
\end{enum}
\end{example}

Let us now turn to the case of the linear iterated map
\begin{equation}
\label{slin12}
y_{k+1} = B y_k
\end{equation}
which admits the solution
\begin{equation}
\label{slin13}
y_k = B^k y_0.
\end{equation}
Using a similar decomposition $B=S+N$ into the semisimple and nilpotent
part, we arrive at
\begin{lemma}
\label{lem_slin2}
Let $b_i$ be the eigenvalues of $B$, and $P_i$, $N_i$ the associated
projectors and nilpotent matrices. Then
\begin{equation}
\label{slin14}
B^k = \sum_{i=1}^m P_i \sum_{j=0}^{\min\set{k,m_i-1}} \binom{k}{j}
b_i^{k-j} N_i^j.
\end{equation}
\end{lemma}
\begin{proof}
The main point is to observe that 
\[
B^k = \Bigpar{\sum_{i=1}^m (b_iP_i+N_i)}^k 
= \sum_{i=1}^m (b_iP_i+N_i)^k 
\]
because all cross-terms vanish. Then one applies the binomial formula.
\end{proof}

For large $k$, the behaviour of $B^k$ is dictated by the terms
$b_i^{k-m_i+1}$. This leads to the following equivalent of Definition
\ref{def_slin1}:

\begin{definition}
\label{def_slin2}
The \defwd{unstable}, \defwd{stable} and \defwd{center subspace} of the
fixed point $x^\star$ are defined, respectively, by
\begin{align}
\nonumber
E^+ &\defby P^+\R^n 
= \bigsetsuch{y}{\lim_{k\to-\infty}B^ky = 0}, 
& P^+ &\defby \sum_{j:\abs{b_j}>1}P_j,
\\
\label{slin15}
E^- &\defby P^-\R^n 
= \bigsetsuch{y}{\lim_{k\to+\infty}B^ky = 0},  
& P^- &\defby \sum_{j:\abs{b_j}<1}P_j,
\\
\nonumber
E^0 &\defby P^0\R^n,  
& P^0 &\defby \sum_{j:\abs{b_j}=1}P_j.
\end{align}
These subspaces are invariant under $B$. The remaining terminology on sinks,
sources, hyperbolic and elliptic points is unchanged. 
\end{definition}

\begin{exercise}
\label{exo_slin}
Find the equilibrium points of the standard map \eqref{sm1} and the Lorenz
equations \eqref{lm9}. Give the dimensions of their stable, unstable and
center subspaces.
{\em Hint:} To determine the sign of the real parts of the eigenvalues of a
$3\times3$ matrix, one can apply the Vi\`ete formula to its characteristic
polynomial.
\end{exercise}


\goodbreak
\subsection{Stability and Liapunov Functions}
\label{ssec_sstab}

\begin{definition}
\label{def_sstab1}
Let $x^\star$ be an equilibrium point of the system $\dot{x}=f(x)$. 
\begin{itemiz}
\item	$x^\star$ is called \defwd{stable} if for any $\eps>0$, one can find
	a $\delta = \delta(\eps)>0$ such that whenever
	$\norm{x_0-x^\star}<\delta$, one has
	$\norm{\ph_t(x_0)-x^\star}<\eps$ for all $t\geqs 0$. 
	
\item	$x^\star$ is called \defwd{asymptotically stable} if it is stable,
	and there is a $\delta_0>0$ such that
	$\lim_{t\to\infty}\ph_t(x_0)=x^\star$ for all $x_0$ such that 
	$\norm{x_0-x^\star}<\delta_0$. 
	
\item	The \defwd{basin of attraction} of an asymptotically stable
	equilibrium $x^\star$ is the set
	\begin{equation}
	\label{sstab1}
	\bigsetsuch{x\in\cD}{\lim_{t\to\infty}\ph_t(x)=x^\star}.
	\end{equation} 
	
\item	$x^\star$ is called \defwd{unstable} if it is not stable.
\end{itemiz}
If $x^\star$ is a fixed point of the map $F$, similar definitions hold with
$\ph_t(\cdot)$ replaced by $F^k(\cdot)$. 
\end{definition}

The linearization of the system $\dot{x}=f(x)$ around an equilibrium
$x^\star$ is the equation $\dot y=Ay$ with $A=\dpar fx(x^\star)$. $x^\star$
is called \defwd{linearly stable} if $y=0$ is a stable equilibrium of its
linearization, and similarly in the asymptotically stable and unstable
cases. Lemma \ref{lem_slin1} shows that 
\begin{itemiz}
\item	$x^\star$ is linearly asymptotically stable if and only
if all eigenvalues of $A$ have a strictly negative real part;
\item	$x^\star$ is linearly stable if and only if no
eigenvalue of $A$ has positive real part, and all purely imaginary
eigenvalues have equal algebraic and geometric multiplicities.
\end{itemiz}
The problem is now to determine relations between linear and nonlinear
stability. A useful method to do this is due to Liapunov. Here we will
limit the discussion to differential equations, although similar results
can be obtained for maps. 

\begin{theorem}[Liapunov]
\label{thm_Liapunov}
Let $x^\star$ be a singular point of $f$, let $\cU$ be a \nbh\ of $x^\star$
and set $\cU_0\defby\cU\setminus\set{x^\star}$. Assume there exists a
continuous function $V:\cU\to\R$, continuously differentiable on $\cU_0$,
such that
\begin{itemiz}
\item[1.] $V(x)>V(x^\star)$ for all $x\in\cU_0$;
\item[2.] the derivative of $V$ along orbits is negative in $\cU_0$, that
is, 
\begin{equation}
\label{sstab2}
\dot{V}(x) \defby \dtot{}{t} V(\ph_t(x))\Bigevalat{t=0} 
= \pscal{\nabla V(x)}{f(x)} \leqs 0 
\qquad \forall x\in\cU_0. 
\end{equation}	
\end{itemiz}
Then $x^\star$ is stable. If, furthermore, 
\begin{itemiz}
\item[3.] the derivative of $V$ along orbits is strictly negative, 
\begin{equation}
\label{sstab3}
\dot{V}(x) = \pscal{\nabla V(x)}{f(x)} < 0 
\qquad \forall x\in\cU_0, 
\end{equation}	
\end{itemiz}
then $x^\star$ is asymptotically stable. 
\end{theorem}
\begin{proof}
Pick $\eps>0$ small enough that the closed ball $\bcB(x^\star,\eps)$ with
center $x^\star$ and radius $\eps$, is contained in $\cU$. Let
$\cS=\partial\bcB(x^\star,\eps)$ be the sphere of radius $\eps$ centered
in $x^\star$. $\cS$ being compact, $V$ admits a minimum on $\cS$, that we
call $\beta$. Consider the open set
\[
\cW = \bigsetsuch{x\in\bcB(x^\star,\eps)}{V(x)<\beta}.
\]
$x^\star\in\cW$ by Condition 1., and thus there exists $\delta>0$ such that
the open ball $\cB(x^\star,\delta)$ is contained in $\cW$. For any
$x_0\in\cB(x^\star,\delta)$, we have $V(\ph_t(x_0))<\beta$ for all $t\geqs
0$, and thus $\ph_t(x_0)\in\cW$ for all $t\geqs0$ by Condition 2., which
proves that $x^\star$ is stable. 

Assume now that \eqref{sstab3} holds. Since the positive orbit of
$x_0\in\cW$ is bounded, it admits a convergent subsequence
$(x_n)_{n\geqs0} = (\ph_{t_n}(x_0))_{n\geqs0} \to y^\star\in\bcW$,
$t_n\to\infty$. Consider the function $t\mapsto V(\ph_t(x_0))$. It is
continuously differentiable, monotonously decreasing, and admits a
subsequence converging to $V(y^\star)$. Thus $V(\ph_t(x_0))$ must converge
to $V(y^\star)$ as $t\to\infty$. Let $\delta>0$ be a small constant and
define the compact set
\[
\cK = \bigsetsuch{x\in\bcW}{V(y^\star)\leqs V(x)\leqs V(y^\star)+\delta}.
\]
If $y^\star\neq x^\star$, then $x^\star\notin\cK$, and thus the maximum of
$\dot V$ on $\cK$ is a strictly negative constant $c$. Take $n$ large enough
that $x_n\in\cK$. Then $\ph_t(x_n)\in\cK$ for all $t\geqs 0$. But this
implies that $V(\ph_t(x_n))\leqs V(x_n)+ct$ which becomes smaller than
$V(y^\star)$ for $t$ large enough, which is impossible. Thus
$y^\star=x^\star$, and all orbits starting in $\cW$ converge to $x^\star$. 
\end{proof}

\begin{figure}
 \centerline{\psfig{figure=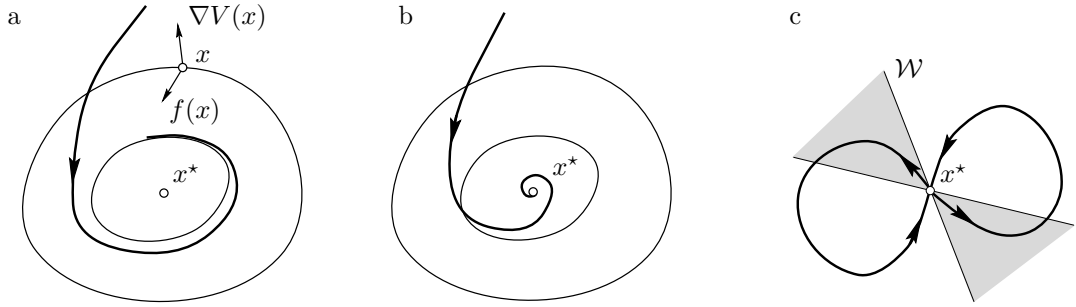,height=40mm,clip=t}}
 \figtext{
 	\writefig	0.2	4.2	a
 	\writefig	5.4	4.2	b
 	\writefig	10.6	4.2	c
 	\writefig	2.4	2.1	$x^\star$
 	\writefig	7.45	2.2	$x^\star$
 	\writefig	12.6	2.1	$x^\star$
 	\writefig	2.7	3.7	$x$
 	\writefig	2.6	4.2	$\nabla V(x)$
 	\writefig	2.35	2.95	$f(x)$
 	\writefig	12.0	3.5	$\cW$
 }
 \captionspace
 \caption[]
 {(a) Stable fixed point with level curves of a Liapunov function. (b)
 Asymptotically stable fixed point, here the vector field must cross all
 level curves in the same direction. (c) Example of an unstable fixed point
 with the set $\cW$ of \v Cetaev's Theorem.}
\label{fig_liap}
\end{figure}

The interpretation of \eqref{sstab2} is that the vector field crosses all 
level sets of $V$ in the same direction (\figref{fig_liap}). $V$ is called
a \defwd{Liapunov function} for $x^\star$, and a \defwd{strict Liapunov
function} if \eqref{sstab3} holds. In fact, the proof also shows that if
$V$ is a strict Liapunov function on $\cU$, and $\cW$ is a set of the form 
$\cW = \setsuch{x}{V(x)<\beta}$ contained in $\cU$, then $\cW$ is contained
in the basin of attraction of $x^\star$. Thus Liapunov functions can be
used to estimate such basins of attraction. 

\begin{exercise}
\label{exo_Liapunov}
Give sufficient conditions on the parameters of the Lorenz equations
\ref{lm9} for {\em all} orbits to converge to the point $(0,0,0)$ (the origin
is said to be \defwd{globally asymptotically stable}). {\em Hint:} Try
Liapunov functions of the form $V(X,Y,Z)=\alpha X^2+\beta Y^2+\gamma Z^2$,
where $\alpha,\beta,\gamma>0$ are constants.
\end{exercise}

The advantage of the Liapunov method is that one does not need to solve the
differential equation. However, the method is not constructive, and the
form of $V$ has to be guessed in each case. In the linearly asymptotically
stable case, such a $V$ can always be constructed, which leads to the
following result.

\begin{cor}
\label{cor_liap}
Assume $x^\star$ is a linearly asymptotically stable equilibrium, that is,
all eigenvalues of $A=\dpar fx(x^\star)$ have strictly negative real parts.
Then $x^\star$ is asymptotically stable. 
\end{cor}
\begin{proof}
There are many different constructions of strict Liapunov functions. We
will give one of them. In order to satisfy Condition 1.\ of the theorem, we
will look for a quadratic form
\[
V(x) = \pscal{(x-x^\star)}{Q(x-x^\star)}
\]
where $Q$ is a symmetric, positive definite matrix. By assumption, there is
a constant $a_0>0$ such that $\re a_i\leqs-a_0$ for all eigenvalues $a_i$ of
$A$. Thus Lemma \ref{lem_slin1} implies that 
\[
\norm{\e^{At}y} \leqs p(t)\e^{-a_0t}\norm{y}
\qquad\forall y,
\]
where $p$ is a polynomial of degree less than $n$. This implies that the
function
\[
V(x) = \int_0^\infty \norm{\e^{As}(x-x^\star)}^2\6 s
\]
exists. $V(x)$ is of the above form with 
\[
Q = \int_0^\infty \e^{\transpose{A}s}\e^{As}\6s.
\]
$Q$ is clearly symmetric, positive definite and bounded, thus there is a
$K>0$ such that $\norm{Qy}\leqs K\norm{y}$ for all $y$. Now we calculate
the following expression in two different ways:
\[
\begin{split}
\int_0^\infty \dtot{}s \bigpar{\e^{\transpose{A}s}\e^{As}} \6s 
&= \lim_{t\to\infty}\e^{\transpose{A}t}\e^{At} - \one = -\one \\
\vrule height 20pt depth 0pt width 0pt
\int_0^\infty \dtot{}s \bigpar{\e^{\transpose{A}s}\e^{As}} \6s 
&= \int_0^\infty \bigpar{\transpose{A} \e^{\transpose{A}s}\e^{As} +
\e^{\transpose{A}s}\e^{As} A} \6s = \transpose{A} Q + Q A.
\end{split}
\]
We have thus proved that
\[
\transpose{A} Q + Q A = -\one.
\]
Now if $y(t)=\ph_t(x)-x^\star$, we have 
\[
\dot V = \dtot{}{t}V(\ph_t(x)) = \dtot{}{t}\bigpar{\pscal{y(t)}{Qy(t)}} =
\pscal{\dot y(t)}{Qy(t)} + \pscal{y(t)}{Q\dot y(t)}. 
\]
Inserting $\dot{y}=Ay+g(y)$ produces two terms. The first one is
\[
\pscal{A y}{Qy} + \pscal{y}{QA y}
= \pscal{y}{(\transpose{A} Q + Q A)y} = - \norm{y}^2,
\]
and the second one gives
\[
\pscal{g(y)}{Qy} + \pscal{y}{Qg(y)} 
= 2 \pscal{g(y)}{Qy} 
\leqs 2\norm{g(y)} \norm{Qy} 
\leqs 2MK\norm{y}^3.
\]
Hence $\dot V \leqs - \norm{y}^2 + 2MK\norm{y}^3$, which shows that $V$
is a strict Liapunov function for $\norm{x-x^\star}<\frac1{2MK}$, and the
result is proved. 
\end{proof}

There exists a characterization of unstable equilibria based on
similar ideas:

\begin{theorem}[\v Cetaev]
\label{thm_Cetaev}
Let $x^\star$ be a singular point of $f$, $\cU$ a \nbh\ of $x^\star$ and
$\cU_0 = \cU\setminus\set{x^\star}$. Assume there exists an open set $\cW$,
containing $x^\star$ in its closure, and a continuous function $V:\cU\to\R$,
which is continuously differentiable on $\cU_0$ and satisfies
\begin{enum}
\item	$V(x)>0$ for all $x\in\cU_0\cap\cW$;
\item	$\dot V(x)>0$ for all $x\in\cU_0\cap\cW$;
\item	$V(x)=0$ for all $x\in\cU_0\cap\partial\cW$.
\end{enum}
Then $x^\star$ is unstable.
\end{theorem}
\begin{proof}
First observe that the definition of an unstable point can be stated as
follows: there exists $\eps>0$ such that, for any $\delta>0$, one can find
$x_0$ with $\norm{x_0-x^\star}<\delta$ and $T>0$ such that
$\norm{\ph_T(x_0)-x^\star}\geqs\eps$. 

Now take $\eps>0$ sufficiently small that $\cB(x^\star,\eps)\subset\cU$.
For any $\delta>0$, $ \cB(x^\star,\delta)\cap\cW\neq\emptyset$. We can thus
take an $x_0\in\cU_0\cap\cW$ such that $\norm{x_0-x^\star}<\delta$, and
by Condition 1.\ $V(x_0)>0$. Now assume by contradiction that
$\norm{\ph_t(x_0)-x^\star}<\eps$ for all $t\geqs0$. $\ph_t(x_0)$ must stay
in $\cU_0\cap\cW$ for all $t$, because it cannot reach the boundary of
$\cW$ where $V=0$. Thus there exists a sequence $x_n=\ph_{t_n}(x_0)$
converging to some $x_1\in\cU_0\cup\cW$. But this contradicts the fact that
$\dot V(x_1)>0$, as in the proof of Theorem \ref{thm_Liapunov}, and thus
$\norm{\ph_t(x_0)-x^\star}$ must become larger than $\eps$. 
\end{proof}

\begin{cor}
\label{cor_cetaev}
Assume $x^\star$ is an equilibrium point such that $A=\dpar fx(x^\star)$ has
at least one eigenvalue with positive real part. Then $x^\star$ is unstable.
\end{cor}
\begin{proof}
Consider first the case of $A$ having no purely imaginary eigenvalues. We
can choose a coordinate system along the unstable and stable subspaces of
$x^\star$, in which the dynamics is described by the equation
\[
\begin{split}
\dot y_+ &= A_+ y_+ + g_+(y) \\
\dot y_- &= A_- y_- + g_-(y),
\end{split}
\]
where all eigenvalues of $A_+$ have a strictly positive real part, all
eigenvalues of $A_-$ have a strictly negative real part, $y=(y_+,y_-)$ and
the terms $g_\pm$ are bounded in norm by a positive constant $M$ times
$\norm{y}^2$. We can define the matrices
\[
Q_- = \int_0^\infty \e^{\transpose{A_-}s}\e^{A_-s} \6s, 
\qquad
Q_+ = \int_0^\infty \e^{-\transpose{A_+}s}\e^{-A_+s} \6s,
\]
which are bounded, symmetric, positive definite, and satisfy
$\transpose{A_-}Q_-+Q_-A_-=-\one$ and $\transpose{A_+}Q_++Q_+A_+=\one$. 
Define the quadratic form
\[
V(y) = \pscal{y_+}{Q_+y_+} - \pscal{y_-}{Q_-y_-}.
\] 
The cone $\cW=\setsuch{y}{V(y)>0}$ is non-empty, because it contains an
eigenvector of $A$ corresponding to an eigenvalue with positive real part.
Proceeding similarly as in the proof of Corollary \ref{cor_liap}, we find
\[
\begin{split}
\dot V(y) &= \norm{y_+}^2 + 2 \pscal{g_+(y)}{Q_+y_+} 
+ \norm{y_-}^2 - 2 \pscal{g_-(y)}{Q_-y_-} \\
&\geqs \norm{y}^2 - 2 KM \norm{y}^3.
\end{split}
\]
Thus \v Cetaev's theorem can be applied to show that $x^\star$ is unstable. 
If $A$ also has purely imaginary eigenvalues, we obtain the additional
equation 
\[
\dot{y}_0 = A_0 y_0 + g_0(y),
\]
where all eigenvalues of $A_0$ are purely imaginary. Let $S$ be an
invertible complex matrix of the same dimension as $A_0$ and consider the
function
\[
V(y) = \pscal{y_+}{Q_+y_+} - \pscal{y_-}{Q_-y_-} - \norm{Sy_0}^2,
\]
where $\pscal{u}{v}=\sum \cc{u_i}v_i$ for complex vectors $u, v$. 
Proceeding as above, we obtain that
\[
\begin{split}
\dot V(y) \geqs{} &\norm{y_+}^2 + \norm{y_-}^2 
- 2 KM \norm{y}^2(\norm{y_+}+\norm{y_-})\\
&- 2\re (\pscal{Sy_0}{A_0Sy_0})
 - 2 \re (\pscal{Sg_0(y)}{Sy_0}).
\end{split}
\]
We shall prove below that for any $\eps>0$, one can construct a matrix
$S(\eps)$ such that $\re (\pscal{Sy_0}{A_0Sy_0}) \leqs \eps\norm{Sy_0}^2$. 
We now take $\cU=\setsuch{y}{\norm{y}<\delta}$, where $\delta>0$ has to be
determined, and $\cW=\setsuch{y}{V(y)>0}$. If $y\in\cW$, we have 
\[
\re (\pscal{Sy_0}{A_0Sy_0}) \leqs \eps\norm{Sy_0}^2 
< \eps \pscal{y_+}{Q_+Y_+} < \eps K \norm{y_+}^2.
\]
We introduce the constants
\[
C(\eps) = \sup_{y_0\neq0}\frac{\norm{Sy_0}}{\norm{y_0}},
\qquad
c(\eps) = \sup_{y_0\neq0}\frac{\norm{y_0}}{\norm{Sy_0}}.
\]
Then we have 
\begin{gather*}
\re (\pscal{Sg_0(y)}{Sy_0}) \leqs C \norm{g_0(y)} \norm{Sy_0} 
\leqs CMK^{1/2} \norm{y}^2 \norm{y_+} \\
\norm{y}^2 \leqs \norm{y_+}^2 + \norm{y_-}^2 + c^2\norm{Sy_0}^2 
< (1+c^2K)(\norm{y_+}^2 + \norm{y_-}^2).
\end{gather*}
Putting everything together, we obtain for all $y\in\cU\cap\cW$
\[
\begin{split}
\dot V(y) & > (1-\eps K) (\norm{y_+}^2 + \norm{y_-}^2) 
- 2KM \norm{y}^2 \bigbrak{(1+CK^{-1/2})\norm{y_+} + \norm{y_-}} \\
& > (\norm{y_+}^2 + \norm{y_-}^2) 
\bigbrak{1 - \eps K - 2KM (1+c(\eps)^2K) (2+C(\eps)K^{-1/2}) \delta}.
\end{split}
\]
Taking $\eps<1/K^2$ and then $\delta$ small enough, we can guarantee that
this quantity is positive for all $y\in\cU_0\cap\cW$, and the corollary is
proved. 
\end{proof}

In the proof we have used the following result.

\begin{lemma}
\label{lem_cetaev}
Assume all the eigenvalues of the $m\times m$ matrix $A_0$ are purely
imaginary. For every $\eps>0$, there exists a complex invertible matrix $S$
such that
\begin{equation}
\label{sstab4}
\bigabs{\re (\pscal{Sz}{SA_0z})} \leqs \eps\norm{Sz}^2
\qquad \forall z\in\C^m.
\end{equation}
\end{lemma}
\begin{proof}
Let $S_0$ be such that $B=S_0A_0S_0^{-1}$ is in complex Jordan canonical
form, that is, the diagonal elements $b_{jj}=\icx\lambda_j$ of $B$ are
purely imaginary, $b_{jj+1}=\sigma_j$ is either zero or one, and all other
elements of $B$ are zero. Let $S_1$ be the diagonal matrix with entries $1,
\eps^{-1}, \dots, \eps^{1-m}$. Then
\[
C \defby S_1BS_1^{-1} = 
\begin{pmatrix}
\icx\lambda_1 & \eps\sigma_1 & & 0 \\
& \ddots & \ddots & \\
& & \ddots & \eps\sigma_{m-1} \\
0 & & & \icx\lambda_m
\end{pmatrix}
\]
Let $S=S_1S_0$ and $u=Sz$. Then
\[
\biggabs{\re \frac{\pscal{Sz}{SA_0z}}{\norm{Sz}^2} }
= \biggabs{\re \frac{\pscal{u}{Cu}}{\norm{u}^2} }
= \biggabs{\re \frac{\sum \icx\lambda_j\abs{u_j}^2 
+ \eps \sum\sigma_j\cc{u_j}u_{j+1}}{\sum \abs{u_j}^2}}
\leqs \eps \frac{\sum\abs{u_j}\abs{u_{j+1}}}{\sum \abs{u_j}^2}
\leqs \eps,
\]
where we used the fact that the upper sum has $m-1$ terms and the lower one
$m$ terms.
\end{proof}

\begin{table}
\begin{center}
\begin{tabular}{|l|c|c|}
\hline
\vrule height 13pt depth 7pt width 0pt
& ODE $\dot x=f(x)$ & Map $x_{k+1}=F(x_k)$ \\
\hline
\vrule height 13pt depth 7pt width 0pt
Conservative & $\nabla\cdot f=0$ & $\bigabs{\det{\dpar Fx}}=1$ \\
\vrule height 0pt depth 7pt width 0pt
Dissipative & $\nabla\cdot f<0$ & $\bigabs{\det{\dpar Fx}}<1$ \\
\hline
\vrule height 13pt depth 7pt width 0pt
Equilibrium & $f(x^\star)=0$ & $F(x^\star)=x^\star$ \\
\hline
\vrule height 13pt depth 7pt width 0pt
Asympt.\ stable if & $\re a_i<0$ $\forall i$ & $\abs{b_i}<1$ $\forall i$ \\
\vrule height 0pt depth 7pt width 0pt
Unstable if & $\exists i:$ $\re a_i>0$ & $\exists i:$ $\abs{b_i}>1$ \\
\hline
\end{tabular}
\end{center}
\caption[]
{Comparison of some properties of ordinary differential equations and
iterated maps. Here $a_i$ and $b_i$ are eigenvalues of the matrices $A=\dpar
fx(x^\star)$ and $B=\dpar Fx(x^\star)$. }
\label{t_sstab}
\end{table}

The two corollaries of this section can be stated as follows: if
$x^\star$ is a hyperbolic equilibrium, then it has the same type of
stability as the linearized system. The same properties are valid for maps
(see \tabref{t_sstab}). This is in general not true for non-hyperbolic
equilibria. This situation will be studied in more detail in Chapter 3.


\subsection{Invariant Manifolds}
\label{ssec_smanif}

For hyperbolic equilibrium points, the analogies between nonlinear and
linear systems can be pushed further. One of them has to do with invariant
manifolds, which generalize the invariant subspaces of the linear case. 
We assume in this section that $f$ is of class $\cC^r$ with $r\geqs 2$. 

\begin{definition}
\label{def_stabmanif}
Let $x^\star$ be a singular point of the system $\dot{x}=f(x)$, and let
$\cU$ be a \nbh\ of $x^\star$. The \defwd{local stable and unstable
manifolds} of $x^\star$ in $\cU$ are defined, respectively, by 
\begin{equation}
\begin{split}
\label{Wloc}
\Wloc{s}(x^\star) &\defby \bigsetsuch{x\in \cU}{\lim_{t\to\infty}\ph_t(x)
= x^\star \text{ and } \ph_t(x)\in \cU \;\forall t\geqs 0} \\
\Wloc{u}(x^\star) &\defby \bigsetsuch{x\in \cU}{\lim_{t\to-\infty}\ph_t(x) 
= x^\star \text{ and } \ph_t(x)\in \cU \;\forall t\leqs 0}.
\end{split}
\end{equation}
The \defwd{global stable and unstable manifolds} of $x^\star$ are defined by 
\begin{equation}
\begin{split}
\label{Wglo}
\Wglo{s}(x^\star) &= \bigcup_{t\leqs 0}\ph_t(\Wloc{s}(x^\star)),\\
\Wglo{u}(x^\star) &= \bigcup_{t\geqs 0}\ph_t(\Wloc{u}(x^\star)).
\end{split}
\end{equation}
\end{definition}

Similar definitions can be made for maps. Global invariant manifolds can
have a very complicated structure, and may return infinitely often to a
\nbh\ of the equilibrium point. This is why one prefers to define separately
local and global invariant manifolds. 
The following theorem states that local invariant manifolds have a nice
structure.

\begin{theorem}[Stable manifold theorem]
\label{thm_stabmanif}
Let $x^\star$ be a hyperbolic equilibrium point of the system $\dot x =
f(x)$, such that the matrix $\dpar fx(x^\star)$ has $n_+$ eigenvalues with
positive real parts and $n_-$ eigenvalues with negative real parts, with
$n_+, n_-\geqs 1$. Then $x^\star$ admits, in a \nbh\ $\cU$,  
\begin{itemiz}
\item	a local stable manifold $\Wloc{s}(x^\star)$, which is a
differentiable manifold of class $\cC^r$ and dimension $n_-$, tangent to
the stable subspace $E^-$ at $x^\star$, and which can be represented as a
graph; 
\item	a local unstable manifold $\Wloc{u}(x^\star)$, which is a
differentiable manifold of class $\cC^r$ and dimension $n_+$, tangent to
the unstable subspace $E^+$ at $x^\star$, and which can be represented as a
graph.
\end{itemiz}
\end{theorem}

\begin{figure}
 \centerline{\psfig{figure=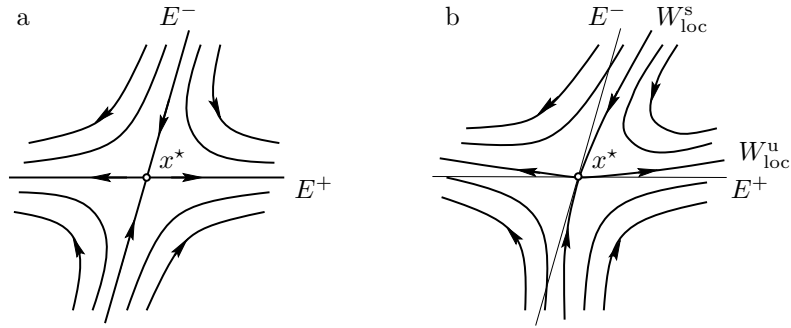,height=40mm,clip=t}}
 \figtext{
 	\writefig	2.7	4.5	a
 	\writefig	8.4	4.5	b
 	\writefig	4.6	2.6	$x^\star$
 	\writefig	10.35	2.6	$x^\star$
 	\writefig	6.4	2.2	$E^+$
 	\writefig	4.6	4.5	$E^-$
 	\writefig	12.2	2.2	$E^+$
 	\writefig	10.3	4.5	$E^-$
 	\writefig	12.3	2.7	$\Wloc{u}$
 	\writefig	11.2	4.5	$\Wloc{s}$
 }
 \captionspace
 \caption[]
 {Orbits near a hyperbolic fixed point: (a) orbits of the
 linearized system, (b) orbits of the nonlinear system with local stable
 and unstable manifolds.}
\label{fig_hyperbol}
\end{figure}

We will omit the proof of this result, but we will give in Chapter 3 the
proof of the center manifold theorem, which relies on similar ideas. Let us
now explain a bit more precisely what this result means. The geometric
interpretation is shown in \figref{fig_hyperbol}. To explain the meaning of
``a differentiable manifold of class $\cC^r$ tangent to $E^\pm$ and
representable as a graph'', let us introduce a coordinate system along the
invariant subspaces of the linearization. The vector field near $x^\star$
can be written as
\begin{equation}
\label{smanif1}
\begin{split}
\dot y_+ &= A_+y_+ + g_+(y_+,y_-) \\
\dot y_- &= A_-y_- + g_-(y_+,y_-),
\end{split}
\end{equation}
where $A_+$ is a $n_+\times n_+$ matrix, which has only eigenvalues with
positive real parts, and $A_-$ is a $n_-\times n_-$ matrix, which has only
eigenvalues with negative real parts. The terms $g_\pm$ are nonlinear and
satisfy $\norm{g_\pm(y_+,y_-)}\leqs M\norm{y}^2$ in $\cU$, where $M$ is a
positive constant. The theorem implies the existence of a function of class
$\cC^r$ 
\begin{equation}
\label{smanif2}
\hu: \cU_+ \to \R^{n_-}, 
\qquad \hu(0)=0,
\qquad \dpar \hu{y_+}(0) = 0,
\end{equation}
where $\cU_+$ is a \nbh\ of the origin in $\R^{n_+}$, such that the local
unstable manifold is given by the equation
\begin{equation}
\label{smanif3}
y_- = \hu(y_+).
\end{equation}
Similar relations hold for the stable manifold. 

In order to determine the function $\hu$, let us compute $\dot y_-$ in two
ways, for a given orbit on the unstable manifold:
\begin{equation}
\label{smanif4}
\begin{split}
\dot y_- &= A_- \hu(y_+) + g_-(y_+,\hu(y_+)) \\
\dot y_- &= \dpar{\hu}{y_+}(y_+)\dot y_+ 
= \dpar{\hu}{y_+}(y_+) \Bigbrak{A_+ y_+ + g_+(y_+,\hu(y_+))}.
\end{split}
\end{equation}
Since both expressions must be equal, we obtain that $\hu$ must satisfy the
partial differential equation
\begin{equation}
\label{smanif5}
A_- \hu(y_+) + g_-(y_+,\hu(y_+)) = \dpar{\hu}{y_+}(y_+) \Bigbrak{A_+ y_+ +
g_+(y_+,\hu(y_+))}.
\end{equation}
This equation is difficult to solve in general. However, since we know by
the theorem that $\hu$ is of class $\cC^r$, we can compute $\hu$
perturbatively, by inserting its Taylor expansion into \eqref{smanif5} and
solving order by order.

\begin{example}
\label{ex_smanif}
Consider, for $n=2$, the system
\begin{equation}
\label{smanif6}
\begin{split}
\dot y_1 &= y_1 \\
\dot y_2 &= -y_2 + y_1^2.
\end{split}
\end{equation}
Then the equation \eqref{smanif5} reduces to
\begin{equation}
\label{smanif7}
- \hu(y_1) + y_1^2 = {\hu}'(y_1) y_1,
\end{equation}
which admits the solution
\begin{equation}
\label{smanif8}
\hu(y_1) = \frac13 y_1^2.
\end{equation}
This is confirmed by the explicit solution of \eqref{smanif6},
\begin{equation}
\label{smanif9}
\begin{split}
y_1(t) &= y_1(0)\e^t \\
y_2(t) &= \Bigpar{y_2(0) - \frac13y_1(0)^2}\e^{-t} + \frac13y_1(0)^2
\e^{2t}.
\end{split}
\end{equation}
\end{example}


\subsection{Normal Forms}
\label{ssec_snform}

In this section we will examine some further connections between the flow
near an equilibrium point and its linearization. We will assume that
$x^\star$ is a singular point of the system $\dot x=f(x)$, and that $f$ is
of class $\cC^r$, $r\geqs 2$, in a \nbh\ of $x^\star$.

\begin{definition}
\label{def_equivalence}
Let $\cU$ and $\cW$ be open sets in $\R^n$, and let $f:\cU\to\R^n$ and
$g:\cW\to\R^n$ be two vector fields of class $\cC^r$. These vector fields
are called
\begin{itemiz}
\item	\defwd{topologically equivalent} if there exists a homeomorphism
$h:\cU\to\cW$ taking the orbits of $\dot x=f(x)$ to the orbits of $\dot
y=g(y)$ and preserving the sense of time;
\item	\defwd{differentiably equivalent} if there exists a diffeomorphism
$h:\cU\to\cW$ taking the orbits of $\dot x=f(x)$ to the orbits of $\dot
y=g(y)$ and preserving the sense of time.
\end{itemiz}
If, in addition, $h$ preserves parametrization of the orbits by time, the
vector fields are called \defwd{conjugate}.
\end{definition}

Equivalence means that if $\ph_t$ and $\psi_t$ are the flows of the two
systems, then 
\begin{equation}
\label{snorm1}
\psi_t\circ h = h\circ \ph_{\tau(t)}
\end{equation}
on $\cU$, where $\tau$ is a homeomorphism from $\R$ to $\R$. If
$\tau(t)=t$ for all $t$, the systems are conjugate.

\begin{theorem}[Hartman-Grobman]
\label{thm_Hartman}
Let $x^\star$ be a hyperbolic equilibrium point of $\dot x=f(x)$, that is,
the matrix $A=\dpar fx(x^\star)$ has no eigenvalue with zero real part.
Then, in a sufficiently small \nbh\ of $x^\star$, $f$ is topologically
conjugate to the linearization $\dot y=Ay$. 
\end{theorem}

Note, however, that topological equivalence is not a very strong property,
since $h$ need not be differentiable. In fact, one can show that all linear
systems with the same number of eigenvalues with positive and negative real
parts are topologically equivalent (see for instance \cite{HK}). So for
instance, the node and focus in \figref{fig_slin} are topologically
equivalent. On the other hand, differentiable equivalence is harder to
achieve as shows the following example.

\begin{example}
\label{ex_snorm}
Consider the following vector field and its linearization:
\begin{align}
\nonumber
\dot{y_1} &= 2y_1 + y_2^2 & 
\dot{z_1} &= 2z_1 \\
\dot{y_2} &= y_2 & 
\dot{z_2} &= z_2.
\label{snorm2}
\end{align}
The orbits can be found by solving the differential equations
\begin{equation}
\label{snorm3}
\dtot{y_1}{y_2} = 2\frac{y_1}{y_2} + y_2, \qquad\qquad
\dtot{z_1}{z_2} = 2\frac{z_1}{z_2},
\end{equation}
which admit the solutions
\begin{equation}
\label{snorm4}
y_1 = \bigbrak{c+\log\abs{y_2}} y_2^2, \qquad\qquad
z_1 = c z_2^2.
\end{equation}
Because of the logarithm, the two flows are $\cC^1$- but not
$\cC^2$-conjugate. 
\end{example}

The theory of normal forms allows to explain these phenomena, and to obtain
conditions for the existence of such $\cC^r$-conjugacies. Consider the system
for $y = x - x^\star$,
\begin{equation}
\label{snorm5}
\dot y = A y + g(y).
\end{equation}
We can try to simplify the nonlinear term by a change of coordinates $y = z
+ h(z)$, which leads to
\begin{equation}
\label{snorm6}
\dot z + \dpar hz(z) \dot z = Az + A h(z) + g(z+h(z)).
\end{equation}
Assume that $h(z)$ solves the partial differential equation
\begin{equation}
\label{snorm7}
\dpar hz(z)Az - A h(z) = g(z+h(z)).
\end{equation}
Then we obtain for $z$ the linear equation 
\begin{equation}
\label{snorm8}
\dot z = Az.
\end{equation}
Unfortunately, we do not know how to solve the equation \eqref{snorm7} in
general. One can, however, work with Taylor series. To this end, we rewrite
the system \eqref{snorm5} as
\begin{equation}
\label{snorm9}
\dot y = A y + g_2(y) + g_3(y) + \dots + g_{r-1}(y) + \Order{\norm{y}^r}.
\end{equation}
Here the last term is bounded in norm by a constant times $\norm{y}^r$, and
the terms $g_k(y)$ are homogeneous polynomial maps of degree $k$ from
$\R^n$ to $\R^n$ ($g_k(\lambda y)=\lambda^k g_k(y)$ $\forall\lambda\in\R$).
Let us denote by $\cH_k$ the set of all such maps. $\cH_k$ is a vector
space for the usual addition and multiplication by scalars. For instance,
when $n=2$, $\cH_2$ admits the basis vectors
\begin{equation}
\label{snorm10}
\begin{pmatrix} y_1^2  \\ 0 \end{pmatrix}, \quad
\begin{pmatrix} y_1y_2 \\ 0 \end{pmatrix}, \quad
\begin{pmatrix} y_2^2  \\ 0 \end{pmatrix}, \quad
\begin{pmatrix} 0 \\ y_1^2  \end{pmatrix}, \quad
\begin{pmatrix} 0 \\ y_1y_2 \end{pmatrix}, \quad
\begin{pmatrix} 0 \\ y_2^2  \end{pmatrix}.
\end{equation}
We now define a linear map from $\cH_k$ to itself given by
\begin{equation}
\label{snorm11}
\fctndef{\ad_k A}{\cH_k}{\cH_k \vrule height 0pt depth 10pt width 0pt}
{h(y)}{\vrule height 10pt depth 0pt width 0pt
\displaystyle\dpar hy(y)Ay - A h(y).}
\end{equation}
The fundamental result of normal form theory is the following:

\begin{prop}
\label{prop_snorm}
For each $k$, $2\leqs k<r$, choose a complementary space $\cG_k$ of the
image of $\ad_k A$, that is, such that $\cG_k\oplus\ad_k A(\cH_k) = \cH_k$.
Then there exists,  in a \nbh\ of the origin, an analytic (polynomial)
change of variables $y=z+h(z)$ transforming \eqref{snorm9} into 
\begin{equation}
\label{snorm12}
\dot z = A z + \gres_2(z) + \gres_3(z) + \dots + \gres_{r-1}(z) 
+ \Order{\norm{z}^r},
\end{equation}
where $\gres_k \in \cG_k$, $2\leqs k<r$. 
\end{prop}
\begin{proof}
The proof proceeds by induction. Assume that for some $k$, $2\leqs k\leqs
r-1$, we have obtained an equation of the form
\[
\dot y = Ay + \sum_{j=2}^{k-1} \gres_j(y) + g_k(y) + \Order{\norm{y}^{k+1}}.
\]
We decompose the term $g_k(y)$ into a resonant and a non-resonant part,
\[
g_k(y) = \gres_k(y) + g^0_k(y), \qquad
g^0_k \in \ad_k A(\cH_k), \quad
\gres_k \in \cG_k.
\]
There exists $h_k\in\cH_k$ satisfying
\[
\ad_k A (h_k(z)) \defby \dpar{h_k}z(z)Az - Ah_k(z) = g^0_k(z).
\]
Observe that $\gres_j(z+h_k(z)) = \gres_j(z)+\Order{\norm{z}^{k+1}}$ for all
$j$, and a similar relation holds for $g_k$. Thus the change of variables
$y=z+h(z)$ yields the equation
\[
\begin{split}
\dot y = \Bigbrak{\one+\dpar{h_k}z(z)} \dot z 
&= Az + Ah_k(z) + \sum_{j=2}^{k-1} \gres_j(z) + g_k(z) +
\Order{\norm{z}^{k+1}}\\
&= \Bigbrak{\one+\dpar{h_k}z(z)}Az + \sum_{j=2}^{k} \gres_j(z) + 
\Order{\norm{z}^{k+1}},
\end{split}
\]
where we have used the definition of $h$ to get the second line. Now for
sufficiently small $z$, the matrix $\one+\dpar{h_k}z$ admits an inverse 
$\one+\Order{\norm{z}^{k-1}}$. Multiplying the above identity on the left by
this inverse, we have proved the induction step.
\end{proof}

The terms $\gres_k(z)$ are called \defwd{resonant} and \eqref{snorm12} is
called the \defwd{normal form} of \eqref{snorm9}. The equation $\ad_k A
(h_k(z)) = g^0_k(z)$ that has to be satisfied to eliminate the non-resonant
terms of order $k$ is called the \defwd{homological equation}. Whether a
term is resonant or not is a problem of linear algebra, which depends only
on the matrix $A$. While it can be difficult to determine the coefficients
of the resonant terms in a particular case, it is in general quite easy to
find which terms can be eliminated. In particular, the following result
holds:

\begin{lemma}
\label{lem_snorm}
Let $(a_1, \dots, a_n)$ be the eigenvalues of $A$, counting multiplicity.
Assume that for each $j$, $1\leqs j\leqs n$, we have 
\begin{equation}
\label{snorm13}
p_1 a_1 + \dots + p_n a_n \neq a_j
\end{equation}
for all $n$-tuples of non-negative integers $(p_1, \dots, p_n)$ satisfying
$p_1+\dots+p_n = k$. Then $\ad_k A$ is invertible, and thus there are no
resonant terms of order $k$.
\end{lemma}
\begin{proof}
We can assume that $A$ is in Jordan canonical form. Consider first the case
of a diagonal $A$. Let $(e_1,\dots,e_n)$ be the canonical basis of $\R^n$.
For $\cH_k$ we choose the basis vectors $z_1^{p_1}\dots z_n^{p_n} e_j$,
where $p_1+\dots+p_n=k$. Then an explicit calculation shows that
\[
\ad_k A (z_1^{p_1}\dots z_n^{p_n} e_j) = 
(p_1 a_1 + \dots + p_n a_n - a_j) z_1^{p_1}\dots z_n^{p_n} e_j.
\]
By assumption, the term in brackets is different from zero. Thus the linear
operator $\ad_k A$ is diagonal in the chosen basis, with nonzero elements on
the diagonal, which shows that it is invertible.

Consider now the case of a matrix $A$ that is not diagonal, but has
elements of the form $a_{jj+1}=1$. Then $\ad_k A$ applied to a basis vector
will contain additional, off-diagonal terms. One of them is proportional to
$z_1^{p_1}\dots z_n^{p_n} e_{j-1}$, while the others are of the form
$z_{k+1} \sdpar{}{z_k}(z_1^{p_1}\dots z_n^{p_n}) e_j$. One can show that
the basis vector of $\cH_k$ can be ordered in such a way that $\ad_k A$ is
represented by a triangular matrix, with the same diagonal elements as in
the case of a diagonal $A$, thus the conclusion is unchanged.
\end{proof}

The non-resonance condition \eqref{snorm13} is called a \defwd{Diophantine
condition}, since it involves integer coefficients. Thus resonant terms can
exist only when the eigenvalues of $A$ satisfy a relation of the form $p_1
a_1 + \dots + p_n a_n = a_j$, which is called a resonance of order
$p_1+\dots+p_n$. In Example \ref{ex_snorm}, the relation $2 a_2 = a_1$
induces a resonance of order $2$, which makes it impossible to eliminate
the term $y_2^2$ by a polynomial change of coordinates. 

In order to solve the question of differentiable equivalence, the really
difficult problem is to eliminate the remainder $\Order{\norm{z}^r}$ in the
normal form \eqref{snorm12}. This problem was solved by Poincar\'e for
sources and sinks, and by Sternberg and Chen for general hyperbolic
equilibria (see for instance \cite{Hartman}). We state here the main result
without proof.

\begin{theorem}[Poincar\'e-Sternberg-Chen]
\label{thm_snorm}
Let $A$ be a $n\times n$ matrix with no eigenvalues on the imaginary axis.
Consider the two equations
\begin{equation}
\label{snorm14}
\begin{split}
\dot{x} &= Ax + b_1(x) \\
\dot{y} &= Ay + b_2(y). 
\end{split}
\end{equation}
We assume that $b_1$ and $b_2$ are of class $\cC^r$, $r\geqs 2$, in a \nbh\
of the origin, that $b_1(x)=\Order{\norm{x}^2}$,
$b_2(x)=\Order{\norm{x}^2}$, and $b_1(x)-b_2(x)=\Order{\norm{x}^r}$. Then
for every $k\geqs 2$, there is an integer $N=N(k,n,A)\geqs k$ such that,
if $r\geqs N$, there exists a map $h$ of class $\cC^k$ such that the two
systems can be transformed into one another by the transformation
$x=y+h(y)$. 
\end{theorem}

This result implies that for the system $\dot y=Ay+g(y)$ to be
$\cC^k$-conjugate to its linearization $\dot z=Az$, it must be sufficiently
smooth and satisfy non-resonance conditions up to sufficiently high order
$N$, where this order and the conditions depend only on $k$, $n$ and the
eigenvalues of $A$. In the special case of all eigenvalues of $A$ having
real parts with the same sign (source or sink), Poincar\'e showed that
$N=k$. For general hyperbolic equilibria, $N$ can be much larger than $k$.

Another, even more important consequence, is that vector fields near
singular points can be classified by their normal forms, each normal form
being a representative of an equivalence class (with respect to
$\cC^k$-conjugacy). This property plays an important role in bifurcation
theory, as we shall see in Chapter 3.

Let us finally remark that in the much more difficult case of
non-hyperbolic equilibria, certain results on $\cC^k$-conjugacy have been
obtained by Siegel, Moser, Takens and others.


\goodbreak
\section{Periodic Solutions}
\label{sec_per}


\subsection{Periodic Orbits of Maps}
\label{sec_pmap}

\begin{definition}
\label{def_pmap}
Let $p\geqs 1$ be an integer. A \defwd{periodic orbit of period $p$} of the
map $F$ is a set of points $\set{x^\star_1,\dots,x^\star_p}$ such that
\begin{equation}
\label{pmap0}
F(x^\star_1)=x^\star_2, \quad \dots \quad F(x^\star_{p-1})=x^\star_p,
\quad F(x^\star_p)=x^\star_1.
\end{equation} 
Each point of the orbit is called a \defwd{periodic point of period $p$}.
Thus a periodic point $x^\star$ of period $p$ is also a fixed point of
$F^p$. $p$ is called the \defwd{least period} of $x^\star$ if
$F^j(x^\star)\neq x^\star$ for $1\leqs j<p$. 
\end{definition}

To find periodic orbits of an iterated map, it is thus sufficient to find
the fixed points of $F^p$, $p=1,2,\dots$. Unfortunately, this becomes
usually extremely difficult with increasing $p$. Moreover, the number of
periodic orbits of period $p$ often grows very quickly with $p$. Methods
that simplify the search for periodic orbits are known for special classes
of maps. For instance, for two-dimensional conservative maps, there exists
a variational method: periodic orbits of period $p$ correspond to
stationary points of some function of $\R^p$ to $\R$. 

Once a periodic orbit has been found, the problem of its linear stability
is rather easily solved. Indeed, it is sufficient to find the eigenvalues of
the matrix
\begin{equation}
\label{pmap1}
\begin{split}
\dpar {F^p}x(x^\star_1) &= 
\dpar Fx(F^{p-1}\bigpar{x^\star_1)} \dpar Fx(F^{p-2}\bigpar{x^\star_1)} 
\dots \dpar Fx(x^\star_1) \\
&= \dpar Fx(x^\star_p) \dpar Fx(x^\star_{p-1})\dots \dpar Fx(x^\star_1).
\end{split}
\end{equation}
Note that the result has to be invariant under cyclic permutations
of the matrices. 

The dynamics near any point of the periodic orbit can be inferred from the
dynamics near one of them, considered as a fixed point of $F^p$. Thus
periodic orbits can also be classified into sinks, sources, hyperbolic and
elliptic orbits, and the concepts of nonlinear stability, invariant
manifolds and normal forms can be carried over from fixed points to
periodic orbits.


\subsection{Periodic Orbits of Flows and Poincar\'e Sections}
\label{sec_pflow}

\begin{definition}
\label{def_pflow}
Let $f$ be a vector field, $\ph_t$ its flow and $T>0$ a constant. A
\defwd{periodic solution of period $T$} of $f$ is a function $\gamma(t)$
such that 
\begin{equation}
\label{pflow0}
\dot\gamma(t)=f(\gamma(t))\text{ and }\gamma(t+T)=\gamma(t) 
\qquad \forall t.
\end{equation} 
The corresponding closed curve $\Gamma=\setsuch{\gamma(t)}{0\leqs t\leqs
T}$ is called a \defwd{periodic orbit of period $T$}. Thus each point $x$
of this orbit is a fixed point of $\ph_T$. $T$ is called the \defwd{least
period} of the orbit if $\ph_t(x)\neq x$ for $0<t<T$. 
\end{definition}

Finding periodic orbits of differential equations is even more difficult
than for maps. There exist methods which help to find periodic orbits in a
number of particular cases, such as two-dimensional flows, or systems
admitting constants of the motion and small perturbations of them. 

Let us now assume that we have found a periodic solution $\gamma(t)$. We
would like to discuss its stability. The difference $y(t) = x(t)-\gamma(t)$
between an arbitrary solution and the periodic solution satisfies the
equation
\begin{equation}
\label{pflow1}
\dot y = f(\gamma(t)+y) - f(\gamma(t)).
\end{equation} 
If $f$ is twice continuously differentiable and $y$ is small, we may expand
$f$ into Taylor series, which yields
\begin{equation}
\label{pflow2}
\dot y = A(t) y + g(y,t), 
\qquad A(t) = \dpar fx(\gamma(t)), 
\quad \norm{g(y,t)}\leqs M\norm{y}^2.
\end{equation}
Let us examine the linearization of this equation, given by
\begin{equation}
\label{pflow3}
\dot y = A(t) y.
\end{equation}
Note that a similar equation already appeared in the proof of Proposition
\ref{prop_vol2}. This equation admits a unique global solution, which can
be represented, because of linearity, as 
\begin{equation}
\label{pflow4}
y(t) = U(t) y(0),
\end{equation}
where $U(t)$ is an $n\times n$ matrix-valued function solving the equation 
\begin{equation}
\label{pflow5}
\dot U(t) = A(t) U(t), \qquad U(0)=\one. 
\end{equation}
The function $U(t)$ is called the \defwd{principal solution} of the equation
\eqref{pflow3}. It should be clear that the linear stability of $\Gamma$ is
related to the asymptotic behaviour of the eigenvalues of $U(t)$.
Unfortunately, there is no general method to determine these eigenvalues.
Note, however, that $A(t)$ is periodic in $t$, and in this case we can say
more.

\begin{theorem}[Floquet]
\label{thm_floquet}
Let $A(t)=A(t+T)$ for all $t$. Then the principal solution of $\dot y = A(t)
y$ can be written as 
\begin{equation}
\label{pflow6}
U(t) = P(t) \e^{Bt},
\end{equation}
where $P(t+T)=P(t)$ for all $t$, $P(0)=\one$, and $B$ is a constant matrix.
\end{theorem}
\begin{proof}
The matrix $V(t)=U(t+T)$ satisfies the equation
\[
\dot V(t) = \dot U(t+T) = A(t+T) U(t+T) = A(t) V(t),
\]
which is the same as \eqref{pflow5}, except for the initial value
$V(0)=U(T)$. We already saw in Proposition \ref{prop_vol2} that $\det
U(t)\neq0$ for all $t$. Thus the matrix $V(t)U(T)^{-1}$ exists and satisfies
\eqref{pflow5}, including the initial condition. By uniqueness of the
solution, it must be equal to $U(t)$:
\[
V(t) U(T)^{-1} = U(t) 
\qquad \Rightarrow \qquad
U(t+T) = U(t) U(T).
\]
We claim that there exists a matrix $B$ such that $U(T)=\e^{BT}$. To see
this, let $\lambda_i\neq 0$ and $m_i$, $i=1,\dots,m$ be the eigenvalues of
$U(T)$ and their algebraic multiplicities. Let $U(T)=\sum_{i=1}^m
(\lambda_i P_i+N_i)$ be the decomposition of $U(T)$ into its semisimple and
nilpotent parts. Then, using Lemma \ref{lem_slin1}, it is easy to check
that 
\[
B = \frac 1T \sum_{i=1}^m \biggpar{\log(\lambda_i)P_i - \sum_{j=1}^{m_i}
\frac{(-N_i)^j}{j\lambda_i^j}}
\]
satisfies $\e^{BT}=U(T)$. Here the sum over $j$ is simply the Taylor
expansion of $\log(\one+N_i/\lambda_i)$.  $B$ is unique up to the
determination of the logarithms. We now define 
\[
P(t) = U(t) \e^{-Bt}.
\]
Then we have for all $t$
\[
P(t+T) = U(t+T) \e^{-B(t+T)} = U(t) \e^{BT} \e^{-B(t+T)} = P(t).
\]
Finally, $P(0)=U(0)=\one$, which completes the proof.
\end{proof}

Floquet's theorem shows that the solution of \eqref{pflow3} can be written
as
\begin{equation}
\label{pflow7}
y(t) = P(t) \e^{Bt} y(0).
\end{equation}
Since $P(t)$ is periodic, the long-time behaviour depends only on $B$. The
eigenvalues of $B$ are called the \defwd{characteristic exponents} of the
equation. $\e^{BT}$ is called the \defwd{monodromy matrix}, and its
eigenvalues, called the \defwd{characteristic multipliers}, are exponentials
of the characteristic exponents times $T$. Computing the characteristic
exponents is difficult in general, but the existence of the representation
\eqref{pflow7} is already useful to classify the possible behaviours near a
periodic orbit.

\begin{exercise}
\label{exo_Floquet}
Consider the Hill equation
\[
\ddot x = -\w(t)^2 x, \qquad
\w(t) = 
\begin{cases}
\Omega & \text{for $nT < t < (n+\frac12)T$} \\
1 & \text{for $(n+\frac12)T < t < (n+1)T$} 
\end{cases}
\quad \forall n\in\Z,
\]
where $T,\Omega$ are positive parameters. Write this system in the form
\eqref{pflow3}, and compute the monodromy matrix and the characteristic
exponents. Plot the exponents as a function of $\Omega$ for $T=\pi$. 
{\em Hint:} The eigenvalues of a $2\times 2$ matrix can be expressed as a
function of its determinant and its trace.
\end{exercise}

Once we have determined the linear stability of the periodic orbit, we
could proceed in a similar way as in the case of a stationary point, in
order to determine the nonlinear stability, the existence of invariant
manifolds, and similar properties. However, Poincar\'e invented a
remarkable method, which allows to shortcut all these steps by reducing the
problem to a simpler one, which has already been studied. Appropriately
enough, this method is called the \defwd{Poincar\'e section}. 

\begin{figure}
 \centerline{\psfig{figure=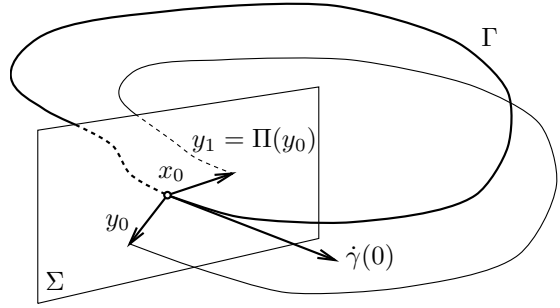,height=40mm,clip=t}}
 \figtext{
 	\writefig	4.2	0.65	$\Sigma$
 	\writefig	10.0	3.9	$\Gamma$
 	\writefig	5.7	2.15	$x_0$
 	\writefig	8.2	1.0	$\dot{\gamma}(0)$
 	\writefig	5.0	1.5	$y_0$
 	\writefig	6.15	2.55	$y_1=\Pi(y_0)$
 }
 \captionspace
 \caption[Poincar\'e map]
 {Definition of the Poincar\'e map associated with the periodic
 orbit $\gamma(t)$.}
\label{fig_Poincare}
\end{figure}

\begin{definition}
\label{def_Poincare}
Let $\gamma(t)$ be a periodic solution of period $T$, $x_0=\gamma(0)$, and
let $\Sigma$ be a hyperplane transverse to the orbit at $x_0$ (see
\figref{fig_Poincare}). By continuity of the flow, there is a \nbh\ $\cU$
of $x_0$ in $\Sigma$ such that for all $x=x_0+y\in \cU$, we can define a
continuous map $\tau(y)$, $\tau(0)=T$, such that $\ph_t(x)$ returns for the
first time to $\Sigma$ in a vicinity of $x_0$ at $t=\tau(y)$. The
\defwd{Poincar\'e map} $\Pi$ associated with the periodic orbit is defined
by
\begin{equation}
\label{Poincare}
x_0 + \Pi(y) \defby \ph_{\tau(y)}(x_0 + y).
\end{equation}
\end{definition}

\begin{prop}
\label{prop_Poincare}
The Poincar\'e map is as smooth as the vector field in a \nbh\ of
the origin. The characteristic multipliers of the periodic orbit are given
by $1$ and the $n-1$ eigenvalues of the Jacobian matrix $\dpar{\Pi}{y}(0)$.
\end{prop}
\begin{proof}
The smoothness of $\Pi$ follows directly from the smoothness of the flow and
the implicit function theorem. Let us now observe that 
\[
\dtot{}{t} \dot\gamma(t) = \dtot{}{t} f(\gamma(t)) 
= \dpar fx(\gamma(t)) \dot\gamma(t) = A(t)\dot\gamma(t). 
\]
Thus by Floquet's theorem, we can write
\[
\dot\gamma(t) = P(t) \e^{Bt}\dot\gamma(0),
\]
and, in particular, 
\[
\dot\gamma(0) = \dot\gamma(T) = P(T)\e^{BT}\dot\gamma(0) =
\e^{BT}\dot\gamma(0).
\]
This shows that $\dot\gamma(0)$ is an eigenvector of $\e^{BT}$ with
eigenvalue $1$. Let $(e_1,\dots,e_{n-1},\dot{\gamma}(0))$,
$e_1,\dots,e_{n-1}\in\Sigma$, be a basis of $\R^n$. In this basis, 
\[
\e^{BT} = \begin{pmatrix}\e^{B_{\Sigma}T} & 0 \\
\dots & 1 \end{pmatrix},
\]
where $B_{\Sigma}$ is the restriction of $B$ to $\Sigma$, and the dots
denote arbitrary entries. Now, if we consider momentarily $y$ as a vector
in $\R^n$ instead of $\Sigma$, linearization of \eqref{Poincare} gives
\[
\dpar\Pi y(0) 
= \dpar{}y \Bigpar{\ph_{\tau(y)}(x_0 + y)} \Bigevalat{y=0}
= \dpar{\ph_T}t(x_0) \dpar\tau y(0) + \dpar{\ph_T}x(x_0) 
= \dot\gamma(0) \dpar\tau y(0) + \e^{BT}.
\] 
The first term is a matrix with zero entries except on the last line, so
that $\sdpar{\Pi}{y}(0)$ has the same representation as $\e^{BT}$, save for
the entries marked by dots. In particular, when $y$ is restricted to
$\Sigma$, $\dpar{\Pi}{y}(0)$ has the same eigenvalues as
$\e^{B_{\Sigma}T}$.
\end{proof}

The consequence of this result is that, by studying the Poincar\'e map, we
obtain a complete characterization of the dynamics in a \nbh\ of the
periodic orbit. In particular, if $y=0$, considered as a fixed point of the
Poincar\'e map, admits invariant manifolds, they can be interpreted as the
intersection of $\Sigma$ and invariant manifolds of the periodic orbit. 


\chapter{Local Bifurcations}
\label{ch_bif}


Up to now, we have obtained quite a precise picture of the dynamics near
hyperbolic equilibrium solutions. One might wonder whether it is of any
interest to examine the case of nonhyperbolic equilibria, since a matrix
chosen at random will have eigenvalues on the imaginary axis with
probability zero. This argument no longer works, however, if the dynamical
system depends on a parameter:
\begin{equation}
\label{bif1}
\dot x = f(x,\lambda) 
\qquad\text{or}\qquad
x_{k+1} = F(x_k,\lambda),
\end{equation}
where $\lambda\in\R$ (or $\R^p$). By changing $\lambda$, it is quite
possible to encounter nonhyperbolic equilibria. 

An important result in this connection is the implicit function theorem:
\setcounter{theorem}{0}
\begin{theorem}
\label{thm_ift}
Let $\cN$ be a \nbh\ of $(x^\star,y^\star)$ in $\R^n\times\R^m$. Let
$f:\cN\to\R^n$ be of class $\cC^r$, $r\geqs 1$, and satisfy
\begin{align}
\label{bif2}
f(x^\star,y^\star) &= 0, \\
\det \dpar fx (x^\star,y^\star) &\neq 0.
\label{bif3}
\end{align}
Then there exists a \nbh\ $\cU$ of $y^\star$ in $\R^m$ and a unique
function $\ph:\cU\to\R^n$ of class $\cC^r$ such that 
\begin{align}
\label{bif4}
\ph(y^\star) &= x^\star, \\
f(\ph(y),y) &= 0 \qquad\text{for all $y\in\cU$}.
\label{bif5}
\end{align}
\end{theorem}

This result tells us under which conditions the equation $f(x,y)=0$ ``can be
solved with respect to $x$''. 
Assume $x^\star$ is an equilibrium point of $f(x,\lambda_0)$ and let $A$ be
the linearization $\dpar fx(x^\star,\lambda_0)$. Then the following
situations can occur:
\begin{itemiz}
\item	If $A$ has no eigenvalues with zero real part, then $f$ will admit
equilibrium points $x^\star(\lambda)$ for all $\lambda$ in a \nbh\ of
$\lambda_0$. By continuity of the eigenvalues of a matrix-valued function,
$x^\star(\lambda)$ will be hyperbolic near $\lambda_0$. The curve
$x^\star(\lambda)$ is usually called an \defwd{equilibrium branch} of $f$. 

\item	If $A$ has one or several eigenvalues equal to zero, then the
implicit function theorem can no longer be applied, and various interesting
phenomena can occur. For instance, the number of equilibrium points of $f$
may change at $\lambda=\lambda_0$. Such a situation is called a
\defwd{bifurcation}, and $(x^\star,\lambda_0)$ is called a
\defwd{bifurcation point} of $f$. 

\item	If $A$ has purely imaginary, nonzero eigenvalues, the implicit
function theorem can still be applied to show the existence of an
equilibrium branch $x^\star(\lambda)$, but its stability may change at
$\lambda=\lambda_0$. This situation is also called a bifurcation. 
\end{itemiz}

Let us point out that we consider here {\em local} bifurcations, that is,
changes of the orbit structure in a small \nbh\ of equilibria. We will not
discuss in any depth {\em global} bifurcations, which involve simultaneous
changes in a larger region of phase space.


\section{Center Manifolds}
\label{sec_cm}


\subsection{Existence of Center Manifolds}
\label{ssec_cmex}

One of the most useful methods to study the flow near a bifurcation point is
the center manifold theorem, which generalizes the stable manifold theorem
(Theorem \ref{thm_stabmanif}) to nonhyperbolic equilibrium points.

\begin{definition}
\label{def_cm1}
Let $\cU\subset\R^n$ be an open set. Let $\cS\subset\R^n$ have the
structure of a differentiable manifold. For $x\in\cU$, let
$(t^x_1,t^x_2)\ni0$ be the maximal interval such that $\ph_t(x)\in\cU$ for
all $t\in(t^x_1,t^x_2)$. $\cS$ is called a \defwd{local invariant manifold}
if $\ph_t(x)\in\cS$ for all $x\in\cS\cap\cU$ and all  $t\in(t^x_1,t^x_2)$. 
\end{definition}

\begin{theorem}
\label{thm_centermanifold}
Let $x^\star$ be a singular point of $f$, where $f$ is of class $\cC^r$,
$r\geqs 2$, in a \nbh\ of $x^\star$. Let $A=\dpar fx(x^\star)$ have,
respectively, $n_+$, $n_0$ and $n_-$ eigenvalues with positive, zero and
negative real parts, where $n_0>0$. Then there exist, in a \nbh\ of
$x^\star$, local invariant $\cC^r$ manifolds $\Wloc{u}$, $\Wloc{c}$ and
$\Wloc{s}$, of respective dimension $n_+$, $n_0$ and $n_-$, and such that 
\begin{itemiz}
\item	$\Wloc{u}$ is the unique local invariant manifold tangent to $E_+$
at $x^\star$, and $\ph_t(x)\to x^\star$ as $t\to-\infty$ for all
$x\in\Wloc{u}$. 
\item	$\Wloc{s}$ is the unique local invariant manifold tangent to $E_-$
at $x^\star$, and $\ph_t(x)\to x^\star$ as $t\to\infty$ for all
$x\in\Wloc{s}$. 
\item	$\Wloc{c}$ is tangent to $E_0$, but not necessarily unique.
\end{itemiz} 
\end{theorem}

Before giving a (partial) proof of this result, we shall introduce a useful
lemma from the theory of differential inequalities. 

\begin{lemma}[Gronwall's inequality]
\label{lem_Gronwall}
Let $\ph$, $\alpha$ and $\beta$ be continuous and real-valued on $[a,b]$,
with $\beta$ non-negative, and assume that
\begin{equation}
\label{cm1}
\ph(t) \leqs \alpha(t) + \int_a^t \beta(s)\ph(s) \6s 
\qquad \forall t\in[a,b].
\end{equation}
Then
\begin{equation}
\label{cm2}
\ph(t) \leqs \alpha(t) + \int_a^t \beta(s)\alpha(s) 
\e^{\int_s^t\beta(u)\6u}\6s 
\qquad \forall t\in[a,b].
\end{equation}
\end{lemma} 
\begin{proof}
Let 
\[
R(t) = \int_a^t \beta(s)\ph(s) \6s.
\]
Then $\ph(t)\leqs\alpha(t)+R(t)$ for all $t\in[a,b]$ and, since
$\beta(t)\geqs0$,  
\[
\dtot Rt(t) = \beta(t)\ph(t) \leqs \beta(t)\alpha(t) + \beta(t)R(t).
\]
Let $K(s)=\int_a^s\beta(u)\6u$. Then 
\[
\dtot {}{s} \e^{-K(s)}R(s) 
= \bigbrak{R'(s)-\beta(s)R(s)}\e^{-K(s)}
\leqs \beta(s)\alpha(s)\e^{-K(s)}
\]
and thus, integrating from $a$ to $t$,
\[
\e^{-K(t)}R(t) \leqs \int_a^t \beta(s)\alpha(s)\e^{-K(s)}\6s.
\]
We obtain the conclusion by multiplying this expression by $\e^{K(t)}$ and
inserting the result into \eqref{cm1}.
\end{proof}

There exist various generalizations of this result, for instance to
functions $\beta(t)$ that are only integrable.

Let us now proceed to the proof of the center manifold theorem. There exist
more or less sophisticated proofs. We will give here a rather
straightforward one taken from \cite{Carr}. For simplicity, we consider the
case $n_+=0$, and we will only prove the existence of a Lipschitz
continuous center manifold. 

\begin{proof}[{\sc Proof of Theorem \ref{thm_centermanifold}}]
We write the system near the equilibrium point as 
\[
\begin{split}
\dot{y} &= By + g_-(y,z) \\
\dot{z} &= Cz + g_0(y,z),
\end{split}
\]
where all eigenvalues of $B$ have strictly negative real parts, all
eigenvalues of $C$ have zero real parts, and $\norm{g_-(y,z)},
\norm{g_0(y,z)}\leqs M(\norm{y}^2+\norm{z}^2)$ in a \nbh\ of the origin. We
shall prove the existence of a local center manifold for a modified
equation, that agrees with the present equation in a small \nbh\ of the
equilibrium. Let $\psi:\R^{n_0}\to[0,1]$ be a $\cC^\infty$ function such
that $\psi(x)=1$ when $\norm{x}\leqs 1$ and $\psi(x)=0$ when $\norm{x}\geqs
2$. We introduce the functions
\[
G_-(y,z) = g_-\bigpar{y,z\psi(\frac z\eps)}, 
\qquad
G_0(y,z) = g_0\bigpar{y,z\psi(\frac z\eps)},
\]
Then the system 
\[
\begin{split}
\dot{y} &= By + G_-(y,z) \\
\dot{z} &= Cz + G_0(y,z),
\end{split}
\]
agrees with the original system for $\norm{z}\leqs\eps$.  We look for a
center manifold with equation $y=h(z)$, where $h$ is in a well-chosen
function space, in which we want to apply Banach's fixed point theorem. Let
$\rho>0$ and $\kappa>0$ be constants, and let $\cX$ be the set of Lipschitz
continuous functions $h:\R^{n_0}\to\R^{n_-}$ with Lipschitz constant
$\kappa$, $\norm{h(z)}\leqs \rho$ for all $z\in\R^{n_0}$ and $h(0)=0$. $\cX$
is a complete space with the supremum norm $\abs{\cdot}$.

For $h\in\cX$, we denote by $\ph_t(\cdot,h)$ the flow of the differential
equation
\[
\dot{z} = Cz + G_0(h(z),z).
\]
For any solution $(y(s),z(s))_{t_0\leqs s\leqs t}$, the relation
\[
y(t) = \e^{B(t-t_0)} y(t_0) + \int_{t_0}^t \e^{B(t-s)} G_-(y(s),z(s))\6s
\]
is satisfied, as is easily checked by differentiation. Among all possible
solutions, we want to select a class of solutions that are tangent to the
center subspace $y=0$ at the origin. In the linear case, the above solutions
satisfy this property only if $\e^{-Bt_0}y(t_0)=0$. We will thus consider
the class of particular solutions satisfying 
\[
y(t) = \int_{-\infty}^t \e^{B(t-s)} G_-(y(s),z(s))\6s
\]
If we set $t=0$ and require that $y(s)=h(z(s))$ for all $s$, we arrive at
the equality
\[
h(z(0)) = \int_{-\infty}^0 \e^{-Bs} G_-\bigpar{h(z(s)),z(s)}\6s, 
\]
where $z(s)=\ph_s(z(0),h)$. Thus we conclude that if $h$ is a fixed point
of the operator $T:\cX\to\cX$, defined by
\[
(Th)(z) = \int_{-\infty}^0 \e^{-Bs} G_-\bigpar{h(\ph_s(z,h)),\ph_s(z,h)}\6s,
\]
then $h$ is a center manifold of the equation. Note that there may be center
manifolds that do not satisfy this equation, and thus we will not be proving
uniqueness. 

Now we want to show that $T$ is a contraction on $\cX$ for an appropriate
choice of $\eps$, $\kappa$ and $\rho$. Observe first that since all
eigenvalues of $B$ have a strictly negative real part, Lemma \ref{lem_slin1}
implies the existence of positive constants $\beta$ and $K$ such that 
\[
\norm{\e^{-Bs}y} \leqs K\e^{\beta s}\norm{y}
\qquad\forall y\in\R^{n_-},\;\forall s\leqs 0.
\]
Since the eigenvalues of $C$ have zero real parts, the same lemma implies
that $\norm{\e^{Cs}z}$ is a polynomial in $s$. Hence, for every $\nu>0$,
there exists a constant $Q(\nu)$ such that 
\[
\norm{\e^{Cs}z} \leqs Q(\nu)\e^{\nu\abs{s}}\norm{z}
\qquad\forall z\in\R^{n_0},\;\forall s\in\R.
\]
It is possible that $Q(\nu)\to\infty$ as $\nu\to0$.
 
We first need to show that $T\cX\subset\cX$. 
Assume from now on that $\rho\leqs\eps$, so that $\norm{h(z)}\leqs\eps$. 
The definition of $G_-$ implies
that
\[
\norm{G_-(h(z),z)}\leqs M'\eps^2
\]
for a constant $M'>0$ depending on $M$. Moreover, the derivatives of $G_0$
are of order $\eps$, so that there is a constant $M''>0$ depending on $M$
with 
\[
\begin{split}
\norm{G_0(h(z_1),z_1) - G_0(h(z_2),z_2)} 
&\leqs M'' \eps\bigbrak{\norm{h(z_1)-h(z_2)} + \norm{z_1-z_2}} \\
&\leqs M'' \eps(1+\kappa) \norm{z_1-z_2}.
\end{split}
\]
A similar relation holds for $G_-$. 
The bound on $\norm{G_-}$ implies that 
\[
\norm{Th(z)} \leqs \int_{-\infty}^0 K\e^{\beta s} M'\eps^2 \6s
\leqs \eps^2\frac{KM'}{\beta},
\]
and hence $\norm{Th(z)} \leqs \rho$ provided $\eps \leqs
(\beta/KM')(\rho/\eps)$.
Next we want to estimate the Lipschitz constant of $T$. Let
$z_1,z_2\in\R^{n_0}$. By definition of the flow $\ph_t(\cdot,h)$, 
\[
\begin{split}
\dtot{}{t} \bigpar{\ph_t(z_1,h)-\ph_t(z_2,h)} 
={}& C \bigpar{\ph_t(z_1,h)-\ph_t(z_2,h)} \\
&+ G_0\bigpar{h(\ph_t(z_1,h)), \ph_t(z_1,h)} 
- G_0\bigpar{h(\ph_t(z_2,h)), \ph_t(z_2,h)}.
\end{split}
\]
Taking into account the fact that $\ph_0(z,h)=z$, we get
\begin{multline*}
\ph_t(z_1,h)-\ph_t(z_2,h) 
={} \e^{Ct} \bigpar{z_1-z_2} \\
+ \int_0^t \e^{C(t-s)} \Bigbrak{G_0\bigpar{h(\ph_s(z_1,h)), \ph_s(z_1,h)} 
- G_0\bigpar{h(\ph_s(z_2,h)), \ph_s(z_2,h)}}\6s.
\end{multline*}
For $t\leqs 0$, we obtain with the properties of $G_0$
\[
\begin{split}
\norm{\ph_t(z_1,h)-\ph_t(z_2,h)} 
\leqs{}& Q(\nu)\e^{-\nu t} \norm{z_1-z_2} \\
&+ \int_t^0 Q(\nu)\e^{-\nu(t-s)} M''\eps(1+\kappa) 
\norm{\ph_s(z_1,h)-\ph_s(z_2,h)}\6s.
\end{split}
\]
We can now apply Gronwall's inequality to 
$\psi(t) = \e^{-\nu t} \norm{\ph_{-t}(z_1,h)-\ph_{-t}(z_2,h)}$, with the
result
\[
\norm{\ph_t(z_1,h)-\ph_t(z_2,h)} \leqs Q(\nu)\e^{-\gamma t} \norm{z_1-z_2},
\]
where $\gamma=\nu+Q(\nu)M''\eps(1+\kappa)$. We can arrange that
$\gamma<\beta$, taking for instance $\nu=\beta/2$ and $\eps$ small enough.
We thus obtain
\[
\begin{split}
\norm{Th(z_1)-Th(z_2)} 
&\leqs \int_{-\infty}^0 K\e^{\beta s}M''\eps(1+\kappa) 
\norm{\ph_s(z_1,h)-\ph_s(z_2,h)} \6s \\
&\leqs \frac{KM''\eps(1+\kappa)Q(\nu)}{\beta-\gamma} \norm{z_1-z_2}.
\end{split}
\]
For any given $\kappa$ and $\nu$, we can find $\eps$ small enough that
$\norm{Th(z_1)-Th(z_2)}\leqs\kappa\norm{z_1-z_2}$. This completes the proof
that $T\cX\subset\cX$.

Finally, we want to show that $T$ is a contraction. For $h_1,h_2\in\cX$,
using 
\begin{multline*}
\norm{G_0\bigpar{h_1(\ph_s(z,h_1)), \ph_s(z,h_1)} 
- G_0\bigpar{h_2(\ph_s(z,h_2)), \ph_s(z,h_2)}} \\
\leqs M''\eps \bigbrak{(1+\kappa)\norm{\ph_s(z,h_1)-\ph_s(z,h_2)} +
\abs{h_2-h_1}},
\end{multline*}
we obtain in a similar way by Gronwall's inequality that 
\[
\norm{\ph_t(z,h_1)-\ph_t(z,h_2)} \leqs 2 \frac{Q(\nu)}{\nu}\eps
M''\abs{h_2-h_1}
\]
(this is a rough estimate, where we have thrown away some $t$-dependent
terms). This leads to the bound
\[
\norm{Th_1(z)-Th_2(z)} \leqs \frac K\beta M''\eps
\Bigbrak{1+2M''\eps(1+\kappa)\frac{Q(\nu)}\nu} \abs{h_2-h_1}.
\]
Again, taking $\eps$ small enough, we can achieve that
$\norm{Th_1(z)-Th_2(z)} \leqs \lambda\abs{h_2-h_1}$ for some $\lambda<1$
and for all $z\in\R^{n_0}$. This shows that $T$ is a contraction, and we
have proved the existence of a Lipschitz continuous center manifold by
Banach's fixed point theorem. One can proceed in a similar way to show that
$T$ is a contraction in a space of Lipschitz differentiable functions.
\end{proof}

One can also prove that if $f$ is of class $\cC^r$, $r\geqs 2$, then $h$ is
also of class $\cC^r$. If $f$ is analytic or $\cC^\infty$, however, then $h$
will be $\cC^r$ for all $r\geqs 1$, but in general it will not be analytic,
and not even $\cC^\infty$. In fact, the size of the domain in which $h$ is
$\cC^r$ may become smaller and smaller as $r$ goes to infinity, see
\cite{Carr} for examples.  

As pointed out in the proof of the theorem, the center manifold is not
necessarily unique. It is easy to give examples of systems admitting a
continuous family of center manifolds. However, as we shall see, these
manifolds have to approach each other extremely fast near the equilibrium
point, and thus the dynamics will be qualitatively the same on all center
manifolds.

\begin{example}
\label{ex_cm1}
The system
\begin{equation}
\label{cm3}
\begin{split}
\dot{y} &= -y \\
\dot{z} &= -z^3
\end{split}
\end{equation}
admits a two-parameter family of center manifolds
\begin{equation}
\label{cm4}
y = h(z,c_1,c_2) = 
\begin{cases}
c_1 \e^{-1/2z^2}
& \text{for $z>0$} \\
0
& \text{for $z=0$} \\
c_2 \e^{-1/2z^2}
& \text{for $z<0$.}
\end{cases}
\end{equation}
The operator $T$ in the proof of Theorem \ref{thm_centermanifold} admits a
unique fixed point $h(z)\equiv0$, but there exist other center manifolds
which are not fixed points of $T$. Note, however, that all functions
$h(z,c_1,c_2)$ have identically zero Taylor expansions at $z=0$. 
\end{example}


\subsection{Properties of Center Manifolds}
\label{ssec_cmprop}

We assume in this section that $x^\star$ is a non-hyperbolic equilibrium
point of $f$, such that $A=\dpar fx(x^\star)$ has $n_0\geqs1$ eigenvalues
with zero real parts and $n_-\geqs1$ eigenvalues with negative real parts.
In appropriate coordinates, we can write
\begin{equation}
\label{cm5}
\begin{split}
\dot{y} &= By + g_-(y,z) \\
\dot{z} &= Cz + g_0(y,z),
\end{split}
\end{equation}
where all eigenvalues of $B$ have strictly negative real parts, all
eigenvalues of $C$ have zero real parts, and $g_-(y,z)$, $g_0(y,z)$ are
nonlinear terms. Theorem \ref{thm_centermanifold} shows the existence of a
local center manifold with parametric equation $y=h(z)$. The dynamics on
this locally invariant manifold is governed by the equation
\begin{equation}
\label{cm6}
\dot u = Cu + g_0(h(u),u).
\end{equation}
Equation \eqref{cm6} has the advantage to be of lower dimension than
\eqref{cm5}, and thus easier to analyse. The following result shows that
\eqref{cm6} is a good approximation of \eqref{cm5}.

\begin{theorem}
\label{thm_cmstab}
If the origin of \eqref{cm6} is stable (asymptotically stable, unstable),
then the origin of \eqref{cm5} is stable (asymptotically stable,
unstable). 

Assume that the origin of \eqref{cm6} is stable. For any solution
$(y(t),z(t))$ of \eqref{cm5} with  $(y(0),z(0))$ sufficiently small, there
exist a solution $u(t)$ of \eqref{cm6} and a constant $\gamma>0$ such that 
\begin{equation}
\label{cm7}
\begin{split}
y(t) &= h(u(t)) + \Order{\e^{-\gamma t}} \\
z(t) &= u(t) + \Order{\e^{-\gamma t}}
\end{split}
\end{equation}
as $t\to\infty$.
\end{theorem}

For a proof, see \cite{Carr}. The center manifold is said to be
\defwd{locally attractive}. Note that if the center manifold is not unique,
then \eqref{cm7} holds for any center manifold.

Let us now discuss how to compute center manifolds. In the proof of Theorem
\ref{thm_centermanifold}, we used the fact that $h$ is a fixed point of a
functional operator $T$. While being useful to prove existence of a center
manifold, $T$ is not very helpful for the computation of $h$, but there
exists another operator for this purpose. Replacing $y$ by $h(z)$ in
\eqref{cm5}, we obtain
\begin{equation}
\label{cm8}
\dpar hz(z) \bigbrak{Cz + g_0(h(z),z)} 
= B h(z) + g_-(h(z),z).
\end{equation}
For functions $\phi:\R^{n_0}\to\R^{n_-}$ which are continuously
differentiable in a \nbh\ of the origin, let us define
\begin{equation}
\label{cm9}
(L\phi)(z) = \dpar \phi z(z) \bigbrak{Cz + g_0(\phi(z),z)} 
- B \phi(z) - g_-(\phi(z),z).
\end{equation}
Then \eqref{cm8} implies that $(Lh)(z)=0$ for any center manifold of
\eqref{cm5}. This equation is impossible to solve in general.
However, its solutions can be computed perturbatively, and the approximation
procedure is justified by the following result, which is also proved in
\cite{Carr}. 

\begin{theorem}
\label{thm_cmapprox}
Let $\cU$ be a \nbh\ of the origin in $\R^{n_0}$ and let
$\phi\in\cC^1(\cU,\R^{n_-})$ satisfy $\phi(0)=0$ and $\dpar\phi z(0)=0$. If
there is a $q>1$ such that $(L\phi)(z) = \Order{\norm{z}^q}$ as $z\to0$,
then $\norm{h(z)-\phi(z)} = \Order{\norm{z}^q}$ as $z\to0$ for any center
manifold $h$.
\end{theorem}

An important consequence of this result is that if $h_1$ and $h_2$ are two
different center manifolds of $x^\star$, then one must have $h_1(z) - h_2(z)
=  \Order{\norm{z}^q}$ for all $q>1$, i.e., all center manifolds have the
same Taylor expansion at $z=0$. 

\begin{example}
\label{ex_cm}
Consider the two-dimensional system
\begin{equation}
\label{cmx1}
\begin{split}
\dot y &= -y + c z^2 \\
\dot z &= yz - z^3,
\end{split}
\end{equation}
where $c$ is a real parameter. We want to determine the stability of the
origin. Naively, one might think that since the first equation suggests
that $y(t)$ converges to zero as $t\to\infty$, the dynamics can be
approximated by projecting on the line $y=0$. This would lead to the
conclusion that the origin is asymptotically stable, because $\dot{z} =
-z^3$ when $y=0$. We will now compute the center manifold in order to find
the correct answer to the question of stability. The operator \eqref{cm9}
has the form
\begin{equation}
\label{cmx2}
(L\phi)(z) = \phi'(z) \bigbrak{z\phi(z)-z^3} + \phi(z) - cz^2. 
\end{equation}
Theorem \ref{thm_cmapprox} allows us to solve the equation $(Lh)(z)=0$
perturbatively, by an Ansatz of the form
\begin{equation}
\label{cmx3}
h(z) = h_2 z^2 + h_3 z^3 + h_4 z^4 + \Order{z^5}. 
\end{equation}
The equation $(Lh)(z)=0$ becomes
\begin{equation}
\label{cmx4}
(Lh)(z) = (h_2-c) z^2 + h_3 z^3 + \bigbrak{h_4 + 2 h_2(h_2-1)} z^4 +
\Order{z^5} = 0,
\end{equation}
which requires $h_2=c$, $h_3=0$ and $h_4=-2c(c-1)$. Hence the center
manifold has a Taylor expansion of the form
\begin{equation}
\label{cmx5}
h(z) = c z^2 - 2c(c-1) z^4 + \Order{z^5},
\end{equation}
and the motion on the center manifold is governed by the equation 
\begin{equation}
\label{cmx6}
\dot{u} = u \mskip1.5mu h(u) - u^3 = (c-1) u^3 - 2c(c-1) u^5 + \Order{u^6}.
\end{equation}
It is easy to show (using, for instance, $u^2$ as a Liapunov function) that
the equilibrium point $u=0$ is asymptotically stable if $c<1$ and unstable
if $c>1$. By Theorem \ref{thm_cmstab}, we conclude that the origin of the
system \eqref{cmx1} is asymptotically stable if $c<1$ and unstable if
$c>1$, which contradicts the naive approach when $c>1$. 

The case $c=1$ is special. In this case, the function $h(z)=z^2$ is an
exact solution of the equation $(Lh)(z)=0$, and the curve $y=z^2$ is the
unique center manifold of \eqref{cmx1}, which has the particularity to
consist only of equilibrium points. The origin is thus stable. 
\end{example}


\section{Bifurcations of Differential Equations}
\label{sec_bifode}

We consider in this section parameter-dependent ordinary differential
equations of the form
\begin{equation}
\label{bo1}
\dot x = f(x,\lambda),
\end{equation}
where $x\in\cD\subset\R^n$, $\lambda\in\Lambda\subset\R^p$ and $f$ is of
class $\cC^r$ for some $r\geqs 2$. We assume that $(x^\star,0)$ is a
bifurcation point of \eqref{bo1}, which means that
\begin{equation}
\label{bo2}
\begin{split}
f(x^\star,0) &= 0 \\
\dpar fx(x^\star,0) &= A,
\end{split}
\end{equation}
where the matrix $A$ has $n_0\geqs1$ eigenvalues on the imaginary axis.
When $\lambda=0$, the equilibrium point $x^\star$ admits a center manifold.
We would like, however, to examine the dynamics of \eqref{bo1} for all
$\lambda$ in a \nbh\ of $0$, where the center manifold theorem cannot be
applied directly. There is, however, an elegant trick to solve this
problem. In the enlarged phase space $\cD\times\Lambda$, consider the
system
\begin{equation}
\label{bo3}
\begin{split}
\dot{x} &= f(x,\lambda) \\
\dot{\lambda} &= 0.
\end{split}
\end{equation}
It admits $(x^\star,0)$ as a non-hyperbolic equilibrium point. The
linearization of \eqref{bo3} around this point is a matrix of the form 
\begin{equation}
\label{bo4}
\begin{pmatrix}
\vrule height 10pt depth 10pt width 0pt
A & \dpar f\lambda(x^\star,0) \\
0 & 0
\end{pmatrix},
\end{equation}
which has $n_0+p$ eigenvalues on the imaginary axis, including the
(possibly multiple) eigenvalue zero. This matrix can be made block-diagonal
by a linear change of variables, where one of the blocks contains all
eigenvalues with zero real part. In these variables, the system \eqref{bo3}
becomes
\begin{equation}
\label{bo5}
\begin{split}
\dot y &= B y + g_-(y,z,\lambda) \\
\dot z &= C z + D\lambda + g_0(y,z,\lambda) \\
\dot \lambda &= 0.
\end{split}
\end{equation}
Here the $(n-n_0)\times(n-n_0)$ matrix $B$ has only eigenvalues with
nonzero real parts, the $n_0\times n_0$ matrix $C$ has all eigenvalues on
the imaginary axis, $D$ is a matrix of size $n_0\times p$, and 
\begin{equation}
\label{bo6}
\norm{g_-(y,z,\lambda)}, \;\norm{g_0(y,z,\lambda)} 
\leqs M\bigpar{\norm{y}^2+\norm{z}^2+\norm{\lambda}^2} 
\end{equation}
in a \nbh\ of the origin, for some positive constant $M$. We can thus apply
Theorem \ref{thm_centermanifold}, which shows the existence of a local
invariant center manifold of the form $y=h(z,\lambda)$. The dynamics on this
manifold is governed by the $n_0$-dimensional equation
\begin{equation}
\label{bo7}
\dot u = Cu + D\lambda + g_0(h(u,\lambda),u,\lambda).
\end{equation}
If $A$ has no eigenvalues with positive real part, Theorem \ref{thm_cmstab}
shows that this equation gives a good approximation to the dynamics of
\eqref{bo1} for small $\lambda$ and near $x=x^\star$. The big advantage is
that generically, the number of eigenvalues on the imaginary axis at the
bifurcation point is small, and thus the reduced equation \eqref{bo7} is of
low dimension. 

With these preliminaries, it becomes possible to investigate bifurcations in
a systematic way. Recall that $A$ is a real matrix, and thus its eigenvalues
are either real or appear in complex conjugate pairs. Thus the two simplest
bifurcations, which we will consider below, involve either a single zero
eigenvalue, or a pair of conjugate imaginary eigenvalues. More complicated
bifurcations correspond to a double zero eigenvalue, a zero eigenvalue and
two conjugate imaginary ones, and so on. These cases are, however, less
\lq\lq generic\rq\rq, and we will not discuss them here. 

\begin{exercise}
\label{exo_cmanif}
Show that the Lorenz equations \eqref{lm9} admit $(X,Y,Z)=(0,0,0)$, $r=1$ as
a bifurcation point. Compute an approximation to second order of the center
manifold in the extended phase space, and deduce the equation governing the
dynamics on the center manifold. {\em Hint:} let $R^2=(\sigma-1)^2 +
4r\sigma$. To put the system into the form \eqref{bo5}, use the
transformation $X=\sigma(z-y_2)$, $Y=\frac12(1-\sigma)(y_2-z) +
\frac12R(y_2+z)$, $Z=y_1$, and $r=1+\lambda$.
\end{exercise}


\subsection{One-Dimensional Center Manifold}
\label{ssec_bifode1}

We first discuss bifurcations involving a single eigenvalue equal to zero,
i.e.\ $n_0=1$ and $C=0$. The dynamics on the one-dimensional
center manifold is governed by an equation of the form
\begin{equation}
\label{bos1}
\dot u = F(u,\lambda), \qquad u\in\R.
\end{equation}
We will also assume that $\lambda\in\R$. The origin $(0,0)$ is a bifurcation
point of \eqref{bos1}, meaning
\begin{equation}
\label{bos2}
F(0,0) = 0, \qquad
\dpar Fu(0,0) = 0.
\end{equation}
In order to understand the dynamics for small $u$ and $\lambda$, we need in
particular to determine the singular points of $F$, that is, we have to
solve the equation $F(u,\lambda)=0$ in a \nbh\ of the origin. Note that the
second condition in \eqref{bos2} implies that we cannot apply the implicit
function theorem. 

Let us start by expanding $F$ in Taylor series,
\begin{equation}
\label{bos3}
F(u,\lambda) = \sum_{\substack{p+q\leqs r \\ p,q\geqs0}} c_{pq} u^p\lambda^q +
\sum_{\substack{p+q=r \\ p,q\geqs0}} u^p\lambda^q R_{pq} (u,\lambda),
\end{equation}
where the functions $R_{pq}$ are continuous near the origin,
$R_{pq}(0,0)=0$ and 
\begin{equation}
\label{bos4}
c_{pq} = \frac1{p!q!}\dpar{^{p+q}F}{u^p\partial\lambda^q} (0,0).
\end{equation}
The bifurcation conditions \eqref{bos2} amount to $c_{00}=c_{10}=0$. 
An elegant way to describe the solutions of $F(u,\lambda)=0$ is based on
\defwd{Newton's polygon}. 

\begin{figure}
 \centerline{\psfig{figure=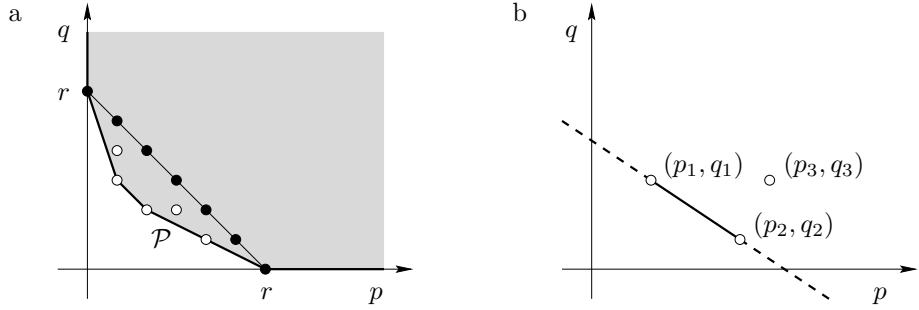,height=40mm,clip=t}}
 \figtext{
 	\writefig	1.0	4.2	a
 	\writefig	1.65	4.0	$q$
 	\writefig	1.65	3.15	$r$
 	\writefig	4.35	0.5	$r$
 	\writefig	2.9	1.2	$\cP$
 	\writefig	5.8	0.5	$p$
 	\writefig	7.7	4.2	b
 	\writefig	8.4	4.0	$q$
 	\writefig	9.7	2.2	$(p_1,q_1)$
 	\writefig	10.9	1.4	$(p_2,q_2)$
 	\writefig	11.3	2.2	$(p_3,q_3)$
 	\writefig	12.5	0.5	$p$
 }
 \captionspace
 \caption[Newton polygon]
 {(a) Definition of the Newton polygon. White circles correspond
 to points $(p,q)$ such that $c_{pq} \neq 0$, black circles to points 
 $(p,q)$ such that $p + q = r$. The slopes of full lines correspond to
 possible exponents of equilibrium branches. (b) The setting of the
 proof of Proposition \ref{prop_Newtonpolygon}.}
\label{fig_Newtonpol}
\end{figure}

\begin{definition}
\label{def_Newtonpolygon}
Consider the set 
\begin{equation}
\label{bos5}
\cA = \bigsetsuch{(p,q)\in\N^2}{\text{$p+q\leqs r$ and $c_{pq}\neq0$}}.
\end{equation}
(In case $F$ is analytic, we simply drop the condition $p+q\leqs r$). 
For each $(p,q)\in\cA$, we construct the sector
$\setsuch{(x,y)\in\R^2}{\text{$x\geqs p$ and $y\geqs q$}}$. The
\defwd{Newton polygon} $\cP$ of \eqref{bos3} is the broken line in $\R^2$
defined by the convex envelope of the union of all these sectors, see
\figref{fig_Newtonpol}a.
\end{definition}

\begin{prop}
\label{prop_Newtonpolygon}
Assume for simplicity\footnote{This assumption allows to neglect all
remainders $R_{pq}(u,\lambda)$, and to avoid pathological situations such
as $F(u,\lambda) = u^2 - \lambda^{5/2}$.} that $c_{pq}\neq0$ whenever
$p+q=r$. Assume further that the equation $F(u,\lambda)=0$ admits a
solution of the form $u=C\abs{\lambda}^\mu(1+\rho(\lambda))$ for small
$\lambda$, where $C\neq0$ and $\rho(\lambda)\to0$ continuously as
$\lambda\to0$. Then Newton's polygon must have a segment of slope $-\mu$. 
\end{prop}
\begin{proof}
It is sufficient to consider a function $F$ of the form
\[
F(u,\lambda) = \sum_{i=1}^3 c_i u^{p_i}\lambda^{q_i}.
\]
Indeed, if the expansion \eqref{bos3} contains only two terms, the result is
immediate, and if it has more than three terms, one can proceed by
induction. The hypothesis implies
\[
\sum_{i=1}^3 \sigma_i c_i C^{p_i}\abs{\lambda}^{\mu p_i + q_i}
(1+\rho(\lambda))^{p_i} = 0,
\]
where $\sigma_i=\pm1$. Assume for definiteness that $p_1\mu+q_1\leqs
p_2\mu+q_2 \leqs p_3\mu+q_3$. Consider first the case $p_1\mu+q_1<
p_2\mu+q_2$. Then division be $\abs{\lambda}^{\mu p_1+q_1}$ gives
\[
\sigma_1 c_1 C^{p_1} (1+\rho(\lambda))^{p_1} + 
\sum_{i=2}^3 \sigma_i c_i C^{p_i}\abs{\lambda}^{\mu p_i + q_i - 
\mu p_1 - q_1}
(1+\rho(\lambda))^{p_i} = 0.
\]
The exponent of $\abs{\lambda}$ is strictly positive. Thus taking the limit
$\lambda\to0$, we obtain $C^{p_1}=0$, a contradiction. We conclude
that we must have $p_1\mu+q_1=p_2\mu+q_2$. 
Graphically, the relation $p_1\mu+q_1 = p_2\mu+q_2 \leqs p_3\mu+q_3$ means
that $\mu = \frac{q_2-q_1}{p_1-p_2}$ is minus the slope of the segment from
$(p_1,q_1)$ to $(p_2,q_2)$, and that $(p_3,q_3)$ lies above, see
\figref{fig_Newtonpol}b. 
\end{proof}

This result does not prove the existence of equilibrium branches
$u=u^\star(\lambda)$, but it tells us where to look. By drawing Newton's
polygon, we obtain the possible values of $\mu$. By inserting the Ansatz
$u^\star(\lambda) = C\abs{\lambda}^\mu(1+\rho)$ into the equation
$F(u,\lambda)=0$, we can determine whether of not such a branch exists.
This is mainly a matter of signs of the coefficients $c_{pq}$. For
instance, the equation $u^2+\lambda^2=0$ admits no solution other than
$(0,0)$, while the equation $u^2-\lambda^2=0$ does admits solutions
$u=\pm\lambda$. 

\subsubsection*{Saddle-Node Bifurcation}

We now illustrate the procedure of computing equilibrium branches in the
generic case $r=2$, $c_{20}\neq0$, $c_{01}\neq0$. Then 
\begin{equation}
\label{bos6}
F(u,\lambda) = c_{01}\lambda + c_{20}u^2 + c_{11}u\lambda + c_{02}\lambda^2
+\sum_{p+q=2} u^p\lambda^q R_{pq}(u,\lambda),
\end{equation}
and Newton's polygon has two vertices $(0,1)$ and $(2,0)$, connected by a
segment with slope $-1/2$. Proposition \ref{prop_Newtonpolygon} tells us
that if there is an equilibrium branch, then it must be of the form
$u=C\abs{\lambda}^{1/2}(1+\rho(\lambda))$, where
$\lim_{\lambda\to0}\rho(\lambda)=0$. In fact, it turns out to be
easier to express $\lambda$ as a function of $u$:

\begin{lemma}
\label{lem_sn}
In a \nbh\ of the origin, there exists a continuous function $\bar\rho(u)$
with $\bar\rho(0)=0$ such that $F(u,\lambda)=0$ if and only if
\begin{equation}
\label{bos7}
\lambda = -\frac{c_{20}}{c_{01}} u^2 (1+\bar\rho(u)).
\end{equation}
\end{lemma}
\begin{proof}
Proposition \ref{prop_Newtonpolygon} indicates that any equilibrium branch
must be of the form $\lambda = \overline C u^2 (1+\bar\rho(u)))$. Define
the function
\[
G(\bar\rho,u) = \frac1{u^2} F(u,\overline C u^2 (1+\bar\rho)).
\]
Using the expansion \eqref{bos6} of $F$, it is easy to see that
\[
\lim_{u\to0} G(0,u) = c_{01}\overline C + c_{20},
\]
and thus $G(0,0)=0$ if and only if $\overline C = -c_{20}/c_{01}$. Moreover,
using a Taylor expansion to first order of $\dpar F\lambda$, one obtains
\[
\lim_{u\to0} \dpar G{\bar\rho}(0,u) = c_{01}\overline C \neq 0.
\]
Thus, by the implicit function theorem, there exists, for small $u$, a
unique function $\bar\rho(u)$ such that $\bar\rho(0)=0$ and
$G(\bar\rho(u),u)=0$.
\end{proof}

\begin{figure}
 \centerline{\psfig{figure=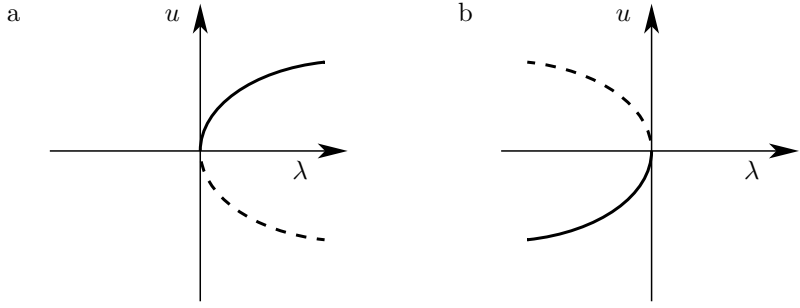,height=40mm}}
 \figtext{
 	\writefig	1.8	4.2	a
 	\writefig	3.9	4.2	$u$
 	\writefig	5.6	2.1	$\lambda$
 	\writefig	7.8	4.2	b
 	\writefig	9.9	4.2	$u$
 	\writefig	11.6	2.1	$\lambda$
 }
 \captionspace
 \caption[]{Saddle-node bifurcation, (a) in the case $c_{01}>0$, $c_{20}<0$
 (direct bifurcation), and (b) in the case $c_{01}>0$, $c_{20}>0$ (indirect
 bifurcation). The other cases are similar, with stable and unstable
 branches interchanged. Full curves indicate stable equilibrium branches,
 while dashed curves indicate unstable equilibrium branches.}
 \label{fig_saddlenode}
\end{figure}

Expressing $u$ as a function of $\lambda$, we find the existence of two
equilibrium branches
\begin{equation}
\label{bos8}
u = u^\star_\pm(\lambda) 
= \pm \sqrt{-\frac{c_{01}}{c_{20}}\lambda}\; \bigbrak{1 + \rho(\lambda)},
\end{equation}
which exist only for $\sign\lambda = -\sign(c_{01}/c_{20})$. Their stability
can be determined by using the Taylor expansion 
\begin{equation}
\label{bos9}
\dpar Fu(u,\lambda) = 2 c_{20}u + c_{11}\lambda 
+\sum_{p+q=1} u^p\lambda^q \overline R_{pq}(u,\lambda),
\end{equation}
where $\overline R_{pq}$ are some continuous functions vanishing at the
origin. Inserting \eqref{bos8}, we get 
\begin{equation}
\label{bos10}
\dpar Fu(u^\star_\pm(\lambda),\lambda) = \pm 2
c_{20}\sqrt{-\frac{c_{01}}{c_{20}}\lambda}
+ \bigorder{\sqrt{\abs{\lambda}}}.
\end{equation}
We thus obtain the following cases, depending on the signs of the
coefficients:
\begin{enum}
\item	If $c_{01}>0$ and $c_{20}<0$, the branches exist for $\lambda>0$,
$u^\star_+$ is stable and $u^\star_-$ is unstable;
\item	if $c_{01}<0$ and $c_{20}>0$, the branches exist for $\lambda>0$,
$u^\star_+$ is unstable and $u^\star_-$ is stable;
\item	if $c_{01}<0$ and $c_{20}<0$, the branches exist for $\lambda<0$,
$u^\star_+$ is stable and $u^\star_-$ is unstable;
\item	if $c_{01}>0$ and $c_{20}>0$, the branches exist for $\lambda<0$,
$u^\star_+$ is unstable and $u^\star_-$ is stable.
\end{enum}
These bifurcations are called \defwd{saddle-node bifurcations}, because
when considering the full system instead of its restriction to the center
manifold, they involve a saddle and a node. If the branches exist for
$\lambda>0$, the bifurcation is called \defwd{direct}, and if they exist
for $\lambda<0$, it is called \defwd{indirect} (\figref{fig_saddlenode}).

It is important to observe that the qualitative behaviour depends only on
those coefficients in the Taylor series which correspond to vertices of
Newton's polygon. Thus we could have thrown away all other terms, to
consider only the truncated equation, or \defwd{normal form}, 
\begin{equation}
\label{bos11}
\dot u = c_{20} u^2 + c_{01} \lambda.
\end{equation}

\subsubsection*{Transcritical Bifurcation}

Let us consider next the slightly less generic case where $c_{01}=0$, but 
$c_{20}\neq0$, $c_{11}\neq0$, $c_{02}\neq0$. Then the Taylor expansion of
$F$ takes the form
\begin{equation}
\label{bos12}
F(u,\lambda) = c_{20}u^2 + c_{11}u\lambda + c_{02}\lambda^2
+\sum_{p+q=2} u^p\lambda^q R_{pq}(u,\lambda),
\end{equation}
and Newton's polygon has three vertices $(0,2)$, $(1,1)$ and $(2,0)$,
connected by segments of slope $-1$. Proposition \ref{prop_Newtonpolygon}
tells us to look for equilibrium branches of the form
$u=C\lambda(1+\rho(\lambda))$. Proceeding in a similar way as in Lemma
\ref{lem_sn}, we obtain the conditions
\begin{equation}
\label{bos13}
\begin{split}
c_{20}C^2 + c_{11}C + c_{02} &= 0 \\
2 c_{20}C^2 + c_{11}C &\neq 0
\end{split}
\end{equation}
for the existence of a unique equilibrium branch of this form. We thus
conclude that if $c_{11}^2 - 4 c_{20} c_{02} > 0$, there are two
intersecting equilibrium branches. It is easy to see that the linearization
of $f$ around such a branch is $(2c_{20}C+c_{11})\lambda +
\order{\lambda}$, and thus one of the branches is stable, the other is
unstable, and they exchange stability at the bifurcation point. 
This bifurcation is called \defwd{transcritical}.

If $c_{11}^2 - 4 c_{20} c_{02} < 0$, there are no equilibrium branches near
the origin. Finally, if $c_{11}^2 - 4 c_{20} c_{02} = 0$ there may be
several branches with the same slope through the origin.

\begin{figure}
 \centerline{\psfig{figure=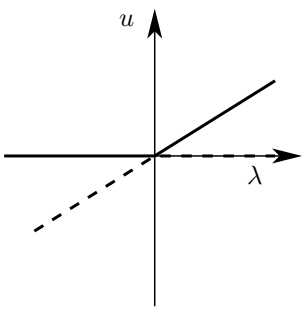,height=40mm}}
 \figtext{
 	\writefig	6.9	4.2	$u$
 	\writefig	8.6	2.1	$\lambda$
 }
 \captionspace
 \caption[]{Transcritical bifurcation, in the case of equation \eqref{bos14}
 with $c_{20}<0$ and $c_{11}>0$.}
 \label{fig_transcritical}
\end{figure}

Let us point out that the condition $c_{02}\neq 0$ is not essential. In
fact, one often carries out a change of variables taking one of the
equilibrium branches to the $\lambda$-axis. The resulting normal form is
\begin{equation}
\label{bos14}
\dot{u} = c_{20} u^2 + c_{11} u \lambda.
\end{equation}
One of the equilibrium branches is $u\equiv0$. The slope of the other branch
and the stability depend on the signs of $c_{20}$ and $c_{11}$
(\figref{fig_transcritical}). 

\subsubsection*{Pitchfork Bifurcation}

One can go on like that for ever, considering cases with more coefficients
in the Taylor series equal to zero, which is not especially interesting
unless one has to do with a concrete problem. However, sometimes symmetries
of the differential equation may cause many terms in the Taylor expansion to
vanish. Consider the case $F\in\cC^3$ satisfying
\begin{equation}
\label{bos15}
F(-u,\lambda) = -F(u,\lambda)
\qquad \forall (u,\lambda).
\end{equation}
Then $c_{pq}=0$ for even $p$. In fact, $F(u,\lambda)/u$ is $\cC^2$ near the
origin and thus we can write
\begin{equation}
\label{bos16}
F(u,\lambda) = u\Bigbrak{c_{11}\lambda + c_{30}u^2 + c_{12}\lambda^2
+\sum_{p+q=2} u^p\lambda^q \overline R_{pq}(u,\lambda)}. 
\end{equation}
The term in brackets is similar to the expansion for the saddle-node
bifurcation, and thus we obtain similar equilibrium branches. In addition,
there is the equilibrium branch $u\equiv 0$. Depending on the signs of the
coefficients, we have the following cases:
\begin{enum}
\item	If $c_{11}>0$ and $c_{30}<0$, the branch $u=0$ is stable for
$\lambda<0$ and unstable for $\lambda>0$, and two additional stable
branches exist for $\lambda>0$;
\item	if $c_{11}<0$ and $c_{30}>0$, the branch $u=0$ is unstable for
$\lambda<0$ and stable for $\lambda>0$, and two additional unstable
branches exist for $\lambda>0$;
\item	if $c_{11}>0$ and $c_{30}>0$, the branch $u=0$ is stable for
$\lambda<0$ and unstable for $\lambda>0$, and two additional unstable
branches exist for $\lambda<0$;
\item	if $c_{11}<0$ and $c_{30}<0$, the branch $u=0$ is unstable for
$\lambda<0$ and stable for $\lambda>0$, and two additional stable branches
exist for $\lambda<0$.
\end{enum}

\begin{figure}
 \centerline{\psfig{figure=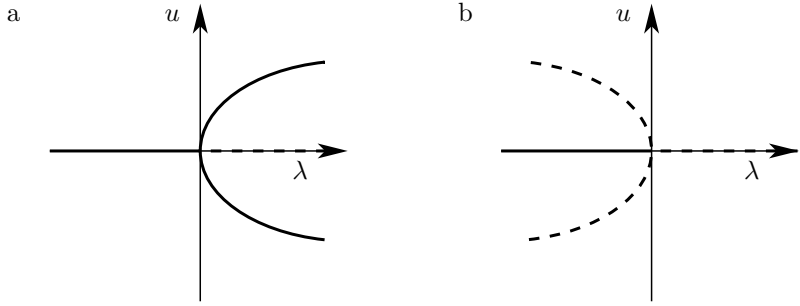,height=40mm}}
 \figtext{
 	\writefig	1.8	4.2	a
 	\writefig	3.9	4.2	$u$
 	\writefig	5.6	2.1	$\lambda$
 	\writefig	7.8	4.2	b
 	\writefig	9.9	4.2	$u$
 	\writefig	11.6	2.1	$\lambda$
 }
 \captionspace
 \caption[]{Pitchfork bifurcation, (a) in the case $c_{11}>0$, $c_{30}<0$
 (supercritical bifurcation), and (b) in the case $c_{11}>0$, $c_{30}>0$
 (subcritical bifurcation).}
 \label{fig_pitchfork}
\end{figure}

This situation in called a \defwd{pitchfork bifurcation}, which is said to
be \defwd{supercritical} if stable equilibrium branches are created, and
\defwd{subcritical} is unstable equilibrium branches are destroyed
(\figref{fig_pitchfork}). The normal form of the pitchfork bifurcation is 
\begin{equation}
\label{bos17}
\dot u = c_{11}\lambda u + c_{30}u^3.
\end{equation}
In Physics, this equation is often written in the form
\begin{equation}
\label{bos18}
\dot{u} = -\dpar Vu(u,\lambda), \qquad
V(u,\lambda) = -\frac{c_{11}}2\lambda u^2 - \frac{c_{30}}4 u^4.
\end{equation}
If $c_{30}<0$, the function $V(u,\lambda)$ has one or two minima, and in the
latter case it is called a \defwd{double-well potential}. 

\begin{exercise}
\label{exo_pitchfork}
Consider the case $F\in\cC^3$ with $c_{20}=c_{11}=0$ and
$c_{30},c_{11},c_{02}\neq 0$, without the symmetry assumption \eqref{bos15}
and discuss the shape and stability of equilibrium branches.
\end{exercise}


\subsection{Two-Dimensional Center Manifold: Hopf Bifurcation}
\label{ssec_bifode2}

We consider now the case $n_0=2$, with $C$ having eigenvalues
$\pm\icx\w_0$, where $\w_0\neq0$. The dynamics on the center manifold is
governed by a two-dimensional system of the form
\begin{equation}
\label{hopf1}
\dot u = F(u,\lambda),
\end{equation}
where we shall assume that $\lambda\in\R$ and $F\in\cC^3$. Since $\dpar
Fu(0,0) = C$ is invertible, the implicit function theorem
(Theorem~\ref{thm_ift}) shows the existence, near $\lambda=0$, of a unique
equilibrium branch $u^\star(\lambda)$, with $u(0)=0$ and
$F(u^\star(\lambda),\lambda)=0$. By continuity of the eigenvalues of a
matrix-valued function, the linearization of $F$ around this branch has
eigenvalues $a(\lambda)\pm\icx\w(\lambda)$, where $\w(0)=\w_0$ and
$a(0)=0$. A translation of $-u^\star(\lambda)$, followed by a linear change
of variables, puts the system \eqref{hopf1} into the form
\begin{equation}
\label{hopf2}
\begin{pmatrix}
\dot u_1 \\ \dot u_2 
\end{pmatrix} 
= 
\begin{pmatrix}
a(\lambda) & -\w(\lambda) \\ 
\w(\lambda) & a(\lambda) 
\end{pmatrix}
\begin{pmatrix}
u_1 \\ u_2 
\end{pmatrix} 
+ 
\begin{pmatrix}
g_1(u_1,u_2,\lambda) \\
g_2(u_1,u_2,\lambda)
\end{pmatrix},
\end{equation}
where the $g_i$ are nonlinear terms satisfying
$\norm{g_i(u_1,u_2,\lambda)}\leqs M(u_1^2+u_2^2)$ for small $u_1$ and $u_2$
and some constant $M>0$. 

\begin{figure}
 \centerline{\psfig{figure=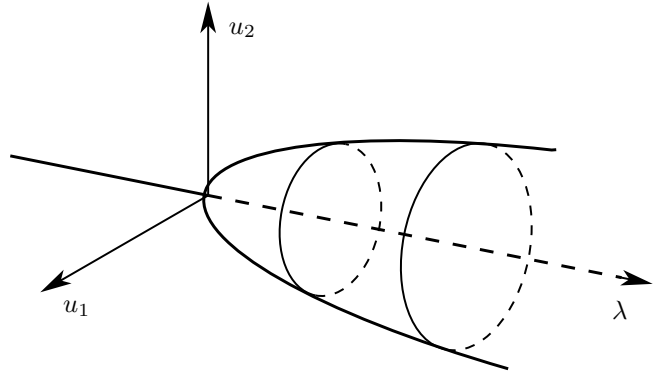,height=50mm,clip=t}}
 \figtext{
 	\writefig	3.8	1.3	$u_1$
 	\writefig	6.0	5.0	$u_2$
 	\writefig	11.1	1.2	$\lambda$
 }
 \captionspace
 \caption[Hopf bifurcation]{Supercritical Hopf bifurcation. The stationary
 solution $(u_1,u_2)=(0,0)$ is stable for $\lambda<0$ and unstable for
 $\lambda>0$. For positive $\lambda$, a stable periodic orbit close to a
 circle of radius $\sqrt\lambda$ appears.}
\label{fig_Hopf}
\end{figure}

Our strategy is now going to be to simplify the nonlinear terms as much as
possible, following the theory of normal forms developed in
Subsection~\ref{ssec_snform}. It turns out to be useful to introduce the
complex variable $z=u_1 + \icx u_2$ (an idea going back to Poincar\'e), 
which satisfies an equation of the form
\begin{equation}
\label{hopf3}
\dot z = \bigbrak{a(\lambda)+\icx\w(\lambda)}z + g(z,\cc z,\lambda),
\end{equation}
where $\cc z$ is the complex conjugate of $z$. This should actually be
considered as a two-dimensional system for the independent variables $z$
and $\cc z$: 
\begin{equation}
\label{hopf4}
\begin{split}
\dot z &= \bigbrak{a(\lambda)+\icx\w(\lambda)}z + g(z,\cc z,\lambda),\\
\dot{\cc{z}} &= \bigbrak{a(\lambda)-\icx\w(\lambda)}\cc z + 
\cc g(z,\cc z,\lambda).
\end{split}
\end{equation}
Lemma~\ref{lem_snorm} shows that monomials of the form $c_{p}z^{p_1}{\cc
z}^{p_2}$ in the nonlinear term $g$ can be eliminated by a nonlinear change
of variables, provided the non-resonance condition \eqref{snorm13} is
satisfied. In the present case, this condition has the form
\begin{equation}
\label{hopf5}
(p_1+p_2-1) a(\lambda) + (p_1-p_2\mp1)\icx\w(\lambda)\neq 0,
\end{equation}
where the signs $\mp$ refer, respectively, to the first and second
equation in \eqref{hopf4}. This condition can also be checked directly by
carrying out the transformation $z=\z+h_p\z^{p_1}{\cc\z}^{p_2}$ in
\eqref{hopf4}. Condition \eqref{hopf5} is hardest to satisfy for
$\lambda=0$, where it becomes
\begin{equation}
\label{hopf6}
(p_1-p_2\mp1)\icx\w_0\neq 0.
\end{equation}
Since $\w_0\neq0$ by assumption, this relation always holds for
$p_1+p_2=2$, so quadratic terms can always be eliminated. The only resonant
term of order $3$ in $g$ is $z^2\cc z = \abs{z}^2 z$, corresponding to
$(p_1,p_2)=(2,1)$. (Likewise, the term $z\cc z^2 = \abs{z}^2 \cc z$ is
resonant in $\cc g$.) We conclude from Proposition~\ref{prop_snorm} that there
exists a polynomial change of variables $z=\z+h(\z,\cc\z)$, transforming
\eqref{hopf3} into
\begin{equation}
\label{hopf7}
\dot \z = \bigbrak{a(\lambda)+\icx\w(\lambda)}\z + c(\lambda)\abs{\z}^2\z +
R(\z,\cc\z,\lambda),
\end{equation}
where $c(\lambda)\in\C$ and $R(\z,\cc\z,\lambda) = \order{\abs{\z}^3}$
(meaning that $\lim_{\z\to 0} \abs{\z}^{-3} R(\z,\cc\z,\lambda) = 0$).
Equation \eqref{hopf7} is the normal form of our bifurcation. To analyse it
further, we introduce polar coordinates $\z=r\e^{\icx\ph}$, in which the
system becomes
\begin{equation}
\label{hopf8}
\begin{split}
\dot r &= a(\lambda)r + \re c(\lambda) r^3 + R_1(r,\ph,\lambda) \\
\dot \ph &= \w(\lambda) + \im c(\lambda) r^2 + R_2(r,\ph,\lambda),
\end{split}
\end{equation}
where $R_1=\order{r^3}$ and $R_2=\order{r^2}$. We henceforth assume that
$a'(0)\neq 0$, and changing $\lambda$ into $-\lambda$ if necessary, we may
assume that $a'(0)>0$. If we discard the remainder $R_1$, the first
equation in \eqref{hopf8} describes a pitchfork bifurcation for $r$, which
is supercritical if $\re c(0)<0$, and subcritical if $\re c(0)>0$. The
first case corresponds to the appearance of a stable periodic orbit, of
amplitude $\sqrt{-a(\lambda)/\re c(\lambda)}$ (\figref{fig_Hopf}), the
second to the destruction of an unstable periodic orbit. The rotation
frequency on this orbit is $\w_0+\Order{\lambda}$. 

It remains to show that this picture is not destroyed by the remainders
$R_1$ and $R_2$. Note that in the present case we cannot apply the
Poincar\'e-Sternberg-Chen theorem (Theorem~\ref{thm_snorm}) because the
linear part is not hyperbolic.

\begin{theorem}[Andronov-Hopf]
\label{thm_hopf}
Assume that the system $\dot x = f(x,\lambda)$ admits an equilibrium branch
$x^\star(\lambda)$ such that the linearization of $f$ at $x^\star(\lambda)$
has two eigenvalues $a(\lambda)\pm\icx\w(\lambda)$, with $a(0)=0$,
$a'(0)>0$ and $\w(0)\neq 0$, and all other eigenvalues have strictly
negative real parts. If the coefficient $c(\lambda)$ in the normal form
\eqref{hopf7} satisfies $\re c(0)\neq 0$, then
\begin{itemiz}
\item	if $\re c(0)<0$, the system admits a stable isolated periodic orbit
for small positive $\lambda$, close to a circle with radius proportional to
$\sqrt{\lambda}$ (supercritical case);
\item	if $\re c(0)>0$, the system admits an unstable isolated periodic orbit
for small negative $\lambda$, close to a circle with radius proportional to
$\sqrt{-\lambda}$ (subcritical case).
\end{itemiz}
\end{theorem}
\begin{proof}
There exist various proofs of this result. One of them is based on the
method of averaging, another one on the Poincar\'e-Bendixson theorem. Since
we did not introduce these methods, we will give a straightforward
geometrical proof. The main idea is to consider the set
$\setsuch{(r,\ph)}{\ph=0,r>0}$ as a Poincar\'e section, and to examine the
associated Poincar\'e map. Consider the case $\re c(0)<0$. If we assume
that $0<\lambda\ll1$, set $r=\sqrt{\lambda}\mskip1.5mu\rho$ in
\eqref{hopf8}, and expand the $\lambda$-dependent terms, we arrive at the
equivalent system
\[
\begin{split}
\dot\rho &= \rho \bigbrak{a'(0)\lambda + \re c(0)\lambda\rho^2 +
\order{\lambda}} \\
\dot\ph &= \w_0 + \lambda\bigbrak{\im c(0)\rho^2 + \w'(0)} +
\order{\lambda}.
\end{split}
\]
Let $\rho_0=\sqrt{-a'(0)/\re c(0)}$. We shall consider this system in the
disc $\rho\leqs 2\rho_0$. The second equation tells us that $\dot\ph\neq 0$
(say, $\dot\ph>0$) for $\lambda$ small enough. Under this condition, we can
use $\ph$ instead of $t$ as new time variable, and write the system in the
form
\[
\begin{split}
\dot\rho &= \frac{a'(0)}{\w_0}\rho \Bigbrak{\lambda -
\frac{\rho^2}{\rho_0^2}\lambda + \order{\lambda}} \\
\dot\ph &= 1.
\end{split}
\]
We define a Poincar\'e map $P$ by the fact that the solution starting at
$(\rho,0)$ passes, after one revolution, through the point
$(P(\rho),2\pi)$. $P(\rho)$ is monotonous by uniqueness of solutions. Fix a
constant $\delta\in(0,1)$. The term $\order{\lambda}$ may depend on $\ph$,
but by taking $\lambda$ small enough, it can be made smaller in absolute
value than $(\delta/2)\lambda$ (uniformly in $\rho$ for $\rho\leqs
2\rho_0$). It follows that
\begin{align*}
&\text{for $0<\rho^2\leqs(1-\delta)\rho_0^2$:} &
\dot\rho &\geqs \frac{a'(0)}{\w_0}\rho\lambda
\Bigpar{1-(1-\delta)-\frac\delta2} > 0 \\
&\text{for $(1+\delta)\rho_0^2\leqs\rho^2\leqs(2\rho_0)^2$:} &
\dot\rho &\leqs \frac{a'(0)}{\w_0}\rho\lambda
\Bigpar{1-(1+\delta)+\frac\delta2} < 0.
\end{align*}
This means that $\rho(t)$ is monotonous outside the annulus
$(1-\delta)\rho_0^2\leqs\rho^2\leqs(1+\delta)\rho_0^2$, and thus $P$ cannot
have fixed points outside the interval
$I=[\sqrt{1-\delta}\,\rho_0,\sqrt{1+\delta}\,\rho_0]$. Furthermore, since
the vector field enters the annulus, $P$ maps the interval $I$ into itself,
and must admit a fixed point $\rho^\star$. Returning to the original
variables, we see that this fixed point corresponds to the desired periodic
orbit. 

Finally, it is also possible to show that $\rho^\star$ is the only fixed
point of $P$ in $(0,2\rho_0)$. We know that there can be no fixed points
outside $I$ (except $0$). By examining the error term $\order{\lambda}$ a
bit more carefully, one finds that $\dot\rho$ is a decreasing function of
$\rho$ near $\rho=\rho_0$. Thus, if $\rho_1<\rho_2$ belong to $I$, with
$\delta$ small enough, the orbits starting in $\rho_1$ and $\rho_2$ must
approach each other as time increases. This means that $P$ is contracting
in $I$, and thus its fixed point is unique.
\end{proof}

This bifurcation is called \defwd{Poincar\'e-Andronov-Hopf bifurcation}.
One should note that it is the first time we prove the existence of
periodic orbits in any generality. Bifurcations with eigenvalues crossing
the imaginary axis are common in parameter-dependent differential equations,
and thus periodic orbits are frequent in these systems.

The nature of the bifurcation depends crucially on the sign of $\re c(0)$.
This quantity can be determined by a straightforward, though rather tedious
computation, and is given, for instance, in \cite[p.\ 152]{GH}.

\begin{exercise}
\label{exo_vdpol}
Show that the van der Pol oscillator
\begin{equation}
\label{vdPol}
\begin{split}
\dot x_1 &= x_2 + \lambda x_1 - \frac13 x_1^3 \\
\dot x_2 &= -x_1
\end{split}
\end{equation}
displays a Hopf bifurcation at the origin, and determine whether it is
subcritical or supercritical.
\end{exercise}


\goodbreak
\section{Bifurcations of Maps}
\label{sec_bifmap}

We turn now to parameter-dependent iterated maps of the form
\begin{equation}
\label{bmap1}
x_{k+1} = F(x_k,\lambda),
\end{equation}
with $x\in\cD\subset\R^n$, $\lambda\in\Lambda\subset\R^p$ and $F\in\cC^r$
for some $r\geqs2$. We assume again that $(x^\star,0)$ is a bifurcation
point of \eqref{bmap1}, which means that 
\begin{equation}
\label{bmap2}
\begin{split}
F(x^\star,0) &= x^\star \\
\dpar Fx(x^\star,0) &= A,
\end{split}
\end{equation}
where $A$ has $n_0\geqs0$ eigenvalues of module $1$. For simplicity, we
assume that all other eigenvalues of $A$ have a module strictly smaller
than $1$. In appropriate coordinates, we can thus write this system as
\begin{equation}
\label{bmap3}
\begin{split}
y_{k+1} &= B y_k + g_-(y_k,z_k,\lambda_k) \\
z_{k+1} &= C z_k + D\lambda_k + g_0(y_k,z_k,\lambda_k) \\
\lambda_{k+1} &= \lambda_k,
\end{split}
\end{equation}
where all eigenvalues of $B$ are inside the unit circle, all eigenvalues of
$C$ are on the unit circle, and $g_-$ and $g_0$ are nonlinear terms. One can
prove, in much the same way we used for differential equations, the
existence of a local invariant center manifold $y = h(z,\lambda)$. This
manifold has similar properties as in the ODE case: it is locally
attractive, and can be computed by solving approximately the equation
\begin{equation}
\label{bmap4}
h\bigpar{Cz + D\lambda + g_0(h(z,\lambda),z,\lambda)} = 
B h(z,\lambda) + g_-(h(z,\lambda),z,\lambda).
\end{equation}
The dynamics on this manifold is governed by the $n_0$-dimensional map
\begin{equation}
\label{bmap5}
u_{k+1} = C u_k + D\lambda + g_0(h(u_k,\lambda),u_k,\lambda).
\end{equation}
Since $C$ is a real matrix with all eigenvalues on the unit circle, the most
generic cases are the following:
\begin{enum}
\item	one eigenvalue equal to $1$: $n_0=1$ and $C=1$;
\item	one eigenvalue equal to $-1$: $n_0=1$ and $C=-1$;
\item	two complex conjugate eigenvalues of module $1$: $n_0=2$ and $C$
having eigenvalues $\e^{\pm2\pi\icx\th_0}$ with $2\th_0\not\in\Z$.
\end{enum}
The first case is easily dealt with, because extremely similar to the case
of differential equations. Indeed, the dynamics on the one-dimensional
center manifold is governed by the equation
\begin{equation}
\label{bmap6}
u_{k+1} = u_k + G(u_k,\lambda), 
\qquad G(0,0)=0, 
\quad \dpar Gu(0,0) = 0.
\end{equation}
The fixed points are obtained by solving the equation $G(u,\lambda)=0$,
which behaves exactly as the equation $F(u,\lambda)=0$ in
Subsection~\ref{ssec_bifode1}. Let us simply indicate the normal forms of
the most common bifurcations. For the saddle-node bifurcation, we have
\begin{equation}
\label{bmap7}
u_{k+1} = u_k + c_{01}\lambda + c_{20}u_k^2;
\end{equation}
for the transcritical bifurcation, one can reduce the equation to
\begin{equation}
\label{bmap8}
u_{k+1} = u_k + c_{11}\lambda u_k + c_{20}u_k^2;
\end{equation}
and for the pitchfork bifurcation, the normal form is given by
\begin{equation}
\label{bmap9}
u_{k+1} = u_k + c_{11}\lambda u_k + c_{30}u_k^3.
\end{equation}

\begin{exercise}
\label{exo_bifmap}
Find the fixed points of the above three maps and determine their stability.
Show that higher order terms do not affect the behaviour near the origin. 
If the map is a Poincar\'e map associated with a periodic orbit $\Gamma$,
what is the meaning of these bifurcations for the flow?
\end{exercise}


\subsection{Period-Doubling Bifurcation}
\label{ssec_bifmap1}

We turn now to the case $n_0=1$, $C=-1$, and assume that the map is of class
$\cC^3$. The map restricted to the center manifold has the form
\begin{equation}
\label{flip1}
\begin{split}
u_{k+1} &= - u_k + G(u_k,\lambda), \\
G(u,\lambda) &= c_{01}\lambda + c_{02}\lambda^2 + c_{11} u\lambda + c_{20}
u^2 + c_{30} u^3 + \cdots
\end{split}
\end{equation}
We first note that the implicit function theorem can be applied to the
equation $u_{k+1}=u_k$, and yields the existence of a unique equilibrium
branch through the origin. It has the form 
\begin{equation}
\label{flip2}
u = u^\star(\lambda) = \frac{c_{01}}2 \lambda +
\order{\lambda}
\end{equation}
and changes stability as $\lambda$ passes through $0$ if
$c_{11}+c_{20}c_{01}\neq0$. Thus something must happen to nearby orbits. To
understand what is going on, it is useful to determine the second iterates.
A straightforward computation gives
\begin{align}
\nonumber
u_{k+2} &= - u_{k+1} + G(u_{k+1},\lambda) \\
\nonumber
&= u_k - G(u_k,\lambda) + G(-u_k + G(u_k,\lambda),\lambda) \\
&= u_k + c_{01}(c_{11} + c_{01}c_{20})\lambda^2 - 2(c_{11} 
+ c_{01}c_{20})u_k\lambda - 2(c_{30}+c_{20}^2)u_k^3 + \cdots
\label{flip3}
\end{align}
Now let us consider the equation $u_{k+2}=u_k$, the solutions of which
yield orbits of period $2$. Using the method of Newton's polygon, we find
that in addition to the solution $u=u^\star(\lambda)$, there exist solutions
of the form
\begin{equation}
\label{flip4}
u^2 = -\frac{c_{11}+c_{20}c_{01}}{c_{30}+c_{20}^2}\lambda +
\order{\lambda},
\end{equation}
provided $c_{30}+c_{20}^2\neq0$ and $c_{11}+c_{20}c_{01}\neq0$.  In other
words, the map $u_k\mapsto u_{k+2}$ undergoes a pitchfork bifurcation. The
equilibrium branch \eqref{flip4} does not correspond to fixed points of
\eqref{flip1}, but to an orbit of least period $2$. We have thus obtained
the following result:

\begin{theorem}
\label{thm_flip}
Let $F(\cdot,\lambda):\R\to\R$ be a one-parameter family of maps of the
form \eqref{flip1}, i.e., such that $F(\cdot,0)$ admits $0$ as a fixed
point with linearization $-1$. Assume
\begin{equation}
\label{flip5}
\begin{split}
c_{11}+c_{20}c_{01}&\neq 0 \\
c_{30} + c_{20}^2 &\neq 0.
\end{split}
\end{equation}
Then there exists a curve of fixed points $u=u^\star(\lambda)$ passing
through the origin, which changes stability at $\lambda=0$. In addition,
there is a curve of the form \eqref{flip4}, tangent to the $u$-axis,
consisting of points of period $2$. The  orbit of period $2$ is stable if
$c_{30} + c_{20}^2>0$ and unstable if $c_{30} + c_{20}^2<0$.
\end{theorem}

This bifurcation is called a \defwd{period doubling}, \defwd{flip} or
\defwd{subharmonic bifurcation}. It is called \defwd{supercritical} if a
stable cycle of period $2$ is created, and \defwd{subcritical} if an
unstable orbit of period $2$ is destroyed. 

\begin{figure}
 \centerline{\psfig{figure=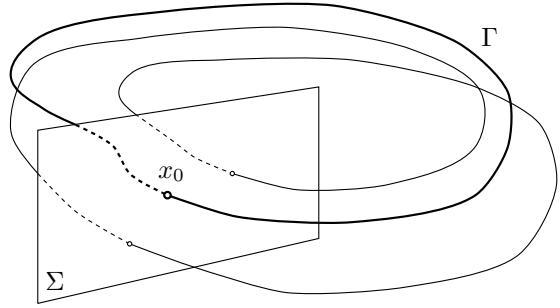,height=40mm}}
 \figtext{
 	\writefig	4.2	0.65	$\Sigma$
 	\writefig	10.0	3.9	$\Gamma$
 	\writefig	5.7	2.15	$x_0$
 }
 \captionspace
 \caption[]{The supercritical flip bifurcation of a Poincar\'e map,
 associated with a periodic $\Gamma$, corresponds to the appearance of a new
 periodic orbit with twice the period.}
 \label{fig_flip}
\end{figure}

\begin{remark}
\label{rem_flip}
If $f$ is a function of class $\cC^3$, its \defwd{Schwartzian derivative} is
defined as
\begin{equation}
\label{flip6}
(Sf)(x) = \frac{f'''(x)}{f'(x)} - \frac32 \biggpar{\frac{f''(x)}{f'(x)}}^2.
\end{equation}
The bifurcation is supercritical if the Schwartzian derivative of the map
at $(0,0)$ with respect to $u$ is negative, and subcritical if this
derivative is positive. 
\end{remark}

The period doubling bifurcation has an interesting consequence if $F$ is
the Poincar\'e map associated with a periodic orbit $\Gamma$. Consider for
instance the supercritical case. For $\lambda<0$, the Poincar\'e map has a
stable fixed point, which means that the periodic orbit $\Gamma$ is stable.
For $\lambda>0$, the fixed point is unstable and there exists a stable
orbit of period two. In phase space, this means that the periodic orbit
$\Gamma$ has become unstable, but a new stable periodic orbit has appeared
(\figref{fig_flip}). If $\Gamma$ has period $T$ for $\lambda=0$, the new
orbit has a period close to $2T$, or half the frequency (which accounts for
the name subharmonic given to the bifurcation). Note that this bifurcation
requires a phase space with dimension at least $3$.


\subsection{Hopf Bifurcation and Invariant Tori}
\label{ssec_bifmap2}

We finally consider what happens when two eigenvalues cross the unit circle
in complex plane. Then we have to study the map
\begin{equation}
\label{tor1}
u_{k+1} = C u_k + G(u_k,\lambda),
\end{equation}
where $u\in\R^2$, and $C$ has eigenvalues $\e^{\pm2\pi\icx\th_0}$ with
$2\th_0\not\in\Z$. We shall assume that $\lambda\in\R$ and $G\in\cC^3$. As
in the case of differential equations, the implicit function theorem shows
the existence of an equilibrium branch $u^\star(\lambda)$ through the
origin. The linearization around this branch has eigenvalues $\mu(\lambda) =
\rho(\lambda)\e^{2\pi\icx\th(\lambda)}$ and $\cc \mu(\lambda)$, with
$\rho(0)=1$ and $\th(0)=\th_0$. An appropriate linear change of variables
casts the system \eqref{tor1} into the form
\begin{equation}
\label{tor2}
z_{k+1} = \mu(\lambda) z_k + g(z_k,\cc z_k,\lambda),
\end{equation}
where $g$ is a nonlinear term, which we will try to simplify by normal form
theory. We can assume that 
\begin{equation}
\label{tor3}
g(z,\cc z,\lambda) = g_2(z,\cc z,\lambda) + g_3(z,\cc z,\lambda) +
\order{\abs{z}^3},
\end{equation}
where $g_2$ is a homogeneous polynomial of degree $2$ in $z$, $\cc z$, and
$g_3$ is of degree $3$. Normal form theory for maps is very similar to
normal form theory for differential equations. In fact, if $h_2(z,\cc
z,\lambda)$ satisfies the \defwd{homological equation}
\begin{equation}
\label{tor4}
h_2(\mu(\lambda)z, \cc\mu(\lambda)\cc z,\lambda) - \mu(\lambda)
h_2(z,\cc z,\lambda) = g_2(z,\cc z,\lambda),
\end{equation} 
then it is easy to see that the transformation $z=\z+h_2(\z,\cc\z,\lambda)$
eliminates terms of order $2$ from \eqref{tor2}. A similar transformation
can be used to eliminate terms of order $3$. Let us now try to solve the
homological equation, assuming 
\begin{equation}
\label{tor5}
h_2(z,\cc z,\lambda) = \sum_{p+q=2} h_{pq}(\lambda)z^p\cc z^q, 
\qquad
g_2(z,\cc z,\lambda) = \sum_{p+q=2} c_{pq}(\lambda)z^p\cc z^q.
\end{equation}
Substitution into \eqref{tor4} shows that 
\begin{equation}
\label{tor6}
h_{pq}(\lambda) = \frac{c_{pq}(\lambda)}{\mu(\lambda)^p\cc\mu(\lambda)^q -
\mu(\lambda)}.
\end{equation}
In particular, for $\lambda=0$, we obtain 
\begin{equation}
\label{tor7}
h_{pq}(0) = \frac{c_{pq}(0)}{\e^{2\pi\icx\th_0}
\bigbrak{\e^{2\pi\icx\th_0(p-q-1)}-1}}.
\end{equation}
This equation can only be solved under the non-resonance condition
\begin{equation}
\label{tor8}
\e^{2\pi\icx\th_0(p-q-1)} \neq 1.
\end{equation}
Since $2\th_0\not\in\Z$, this condition is always satisfied for
$(p,q)=(2,0)$ and $(1,1)$. Thus terms of the form $z^2$ and $z\cc z$ can
always be eliminated. Terms of the form $\cc z^2$, however, can only be
eliminated if $3\th_0$ is not an integer. Similarly, third order terms of
the form $z^3$ and $z\cc z^2$ can always be eliminated, terms of the form
$\cc z^3$ can only be eliminated if $4\th_0$ is not an integer, while the
term $z^2\cc z$ can never be removed. We conclude that if
$\e^{2\pi\icx\th_0}$ is neither a cubic nor a quartic root of unity, there
exists a polynomial change of variables transforming \eqref{tor2} into its
normal form
\begin{equation}
\label{tor9}
\z_{k+1} = \mu(\lambda) \z_k + c(\lambda) \abs{\z_k}^2\z_k +
\order{\abs{\z_k}^3}.
\end{equation}
Using polar coordinates $\z_k = r_k \e^{\icx\ph_k}$, one obtains a map of
the form
\begin{equation}
\label{tor10}
\begin{split}
r_{k+1} &= \rho(\lambda)r_k + \alpha(\lambda) r_k^3 + \order{r_k^3} \\
\ph_{k+1} &= \ph_k + 2\pi\th(\lambda) + \beta(\lambda) r_k^2 + \order{r_k^2}.
\end{split}
\end{equation}
An explicit calculation shows that if
$c(\lambda)=\abs{c(\lambda)}\e^{2\pi\icx\psi(\lambda)}$, then the
coefficients $\alpha$ and $\beta$ are given by
\begin{equation}
\label{tor11}
\alpha(\lambda) = \abs{c(\lambda)} \cos 2\pi(\psi(\lambda)-\th(\lambda)), 
\qquad
\beta(\lambda) = \frac{\abs{c(\lambda)}}{\rho(\lambda)} 
\sin 2\pi(\psi(\lambda)-\th(\lambda)).
\end{equation}
\goodbreak
\noindent
If we neglect the remainders, the map \eqref{tor10} describes a pitchfork
bifurcation for the radial variable $r$. Depending on the signs of
$\rho'(0)$ and $\alpha(0)$, there will be creation or destruction of an
invariant circle. This can be proved to remain true when the
$\order{\cdot}$ terms are present, but the proof is more difficult than in
the case of differential equations. 

\begin{theorem}[Ruelle]
\label{thm_ruelle}
Let $F(\cdot,\lambda):\R^2\to\R^2$ be a one-parameter family of maps,
admitting a smooth curve of fixed points $u^\star(\lambda)$. Assume the
linearization around the fixed points has complex conjugate eigenvalues
$\mu(\lambda)$ and $\cc\mu(\lambda)$ such that
\begin{gather}
\label{tor12a}
\abs{\mu(0)}=1 \quad \text{but $\mu(0)^j\neq 1$ for $j=1,2,3,4$,} \\
\label{tor12b}
\dtot{}{\lambda} \abs{\mu(\lambda)}\Bigevalat{\lambda=0} \neq 0.
\end{gather}
Then there is a smooth change of coordinates transforming the map $F$ into
\eqref{tor10}. If $\alpha(0)\neq 0$, then the map admits, either for small
positive $\lambda$ or for small negative $\lambda$, an invariant curve close
to a circle of radius $\sqrt{\abs\lambda}$. 
\end{theorem} 

\begin{figure}
 \centerline{\psfig{figure=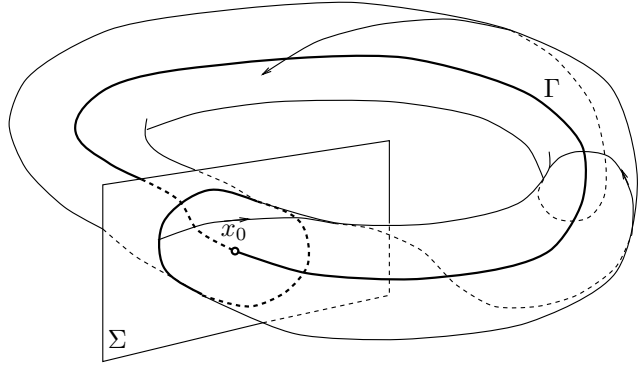,height=48mm}}
 \figtext{
 	\writefig	4.5	0.65	$\Sigma$
 	\writefig	10.3	4.0	$\Gamma$
 	\writefig	6.0	2.15	$x_0$
 }
 \captionspace
 \caption[]{If the Poincar\'e map of a periodic orbit $\Gamma$ undergoes
 supercritical Hopf bifurcation, an attracting invariant torus is created.}
 \label{fig_torus}
\end{figure} 

In the strongly resonant cases $\mu(0)^j= 1$ for $j=1,2,3$ or $4$, the
situation is more complicated, and there is no invariant curve in general.
If such an invariant curve exists, the dynamics on this curve is described
by a \defwd{circle map}. The theory of circle maps is a huge subject in
itself. Roughly speaking, they can be characterized by a \defwd{rotation
number}, measuring the average angle of rotation per iteration. Two cases
can occur:
\begin{itemiz}
\item	if the rotation number is rational, there exists a periodic orbit,
which usually attracts most orbits;
\item	if the rotation number is irrational, all orbits are dense, and
under suitable smoothness assumptions, the map is conjugate to a rotation.
\end{itemiz}
If the bifurcating map is the Poincar\'e map of a periodic orbit $\Gamma$,
the invariant circle will be the intersection of an \defwd{invariant torus}
with the Poincar\'e section $\Sigma$ (\figref{fig_torus}). If the rotation
number is irrational, orbits fill this torus in a dense way. They are
called \defwd{quasiperiodic} (with two frequencies), which means that any
solution $x(t)$ on the invariant torus can be written as 
\begin{equation}
\label{tor13}
x(t) = H(\w t, \w' t), \qquad
H(u+1,v) = H(u,v+1) = H(u,v),
\end{equation}
where $\w/\w'$ is equal to the irrational rotation number. Besides invariant
points, curves, and manifolds, we have thus found a new kind of invariant
set appearing quite commonly in dynamical systems.


\chapter{Introduction to Chaotic Dynamics}
\label{ch_chaos}

The characterization of chaotic dynamics is a large subject, with many
recent developments, and ramifications in several domains of Mathematics.
We will discuss here only a few selected topics, the major objective being
to provide an idea of what we call chaotic motion. 

The definition of chaos may vary from system to system, but usually at least
one of the following elements is present:
\begin{itemiz}
\item	The time dependence of solutions is more complicated than
stationary, periodic or quasiperiodic.
\item	The motion is very sensitive to variations in the initial
conditions: nearby solutions diverge exponentially fast.
\item	The asymptotic motion takes place on a geometrically complicated
object (often a fractal), called a \defwd{strange attractor}.
\item	Chaotic orbits coexist with a (countable) infinity of unstable
periodic orbits; the number of orbits of period less or equal $T$ grows
exponentially fast with $T$.
\item	As time goes by, the images under the flow of any two subsets of
phase space get entangled in a complicated way.
\end{itemiz}
These properties do not necessarily all occur at once: for instance, strange
attractors cannot occur in conservative systems, but are quite typical in
chaotic dissipative systems. Often, there are subtle relations between the
above properties, many of which can be characterized quantitatively (by
Liapunov exponents, topological entropy, \dots).

Our approach in this introduction to chaotic dynamics will be to start with
some very simple examples, where many of these chaotic properties can be
proven to hold. Then we will show how similar properties can be proven to
exist for more realistic systems. 


\section{Symbolic Dynamics}
\label{sec_sd}

Symbolic dynamics is a very useful technique for characterizing dynamical
systems. It consists in associating with every orbit a sequence of
symbols, related to a partition of phase space, and describing in
which order the orbit of $x$ visits the elements of the partition. In this
way, one can sometimes reduce the problem to a combinatorial one, by
transposing the dynamics to the space of allowed symbolic sequences.

We will start by illustrating the method on a very simple one-dimensional
map called the tent map, and later discuss a less trivial two-dimensional
map which is already useful in proving existence of chaotic orbits for a
more general class of systems. 
\goodbreak


\subsection{The Tent Map}
\label{ssec_tnt}

In Section \ref{sec_logistic}, we claimed that the logistic map
\begin{equation}
\label{tnt1}
f_\lambda(y) = \lambda y(1-y)
\end{equation}
is as random as coin tossing when $\lambda=4$. Let us now explain what we
meant. It is easy to check that the transformation
$y=h(x)=\frac12(1-\cos\pi x)$ transforms the logistic map $f_4$ into the map
\begin{equation}
\label{tnt2}
g_2(x) = h^{-1}\circ f_4 \circ h(x) = 
\begin{cases}
2x & \text{for $0\leqs x\leqs\frac12$} \\
2-2x & \text{for $\frac12\leqs x\leqs1$.}
\end{cases}
\end{equation}
This map is called the \defwd{tent map} because of its triangular shape (one
can define more general $g_\mu$ with slope $\mu\neq2$). In the
transformation, we have lost differentiability at $x=\frac12$, but the
piecewise linearity will be useful to classify the orbits. 

First of all, we note that the map \eqref{tnt2} has two unstable fixed
points, at $0$ and $2/3$. We also observe that $1$ is mapped to $0$ and
$1/2$ is mapped to $1$, so there exist \lq\lq transient\rq\rq\ orbits,
ending at $0$ after finitely many iterations. 

In order to understand the other orbits, it turns out to be a good idea to
write $x$ in binary expansion. We will write
\begin{equation}
\label{tnt3}
x = \sum_{i=0}^\infty b_i 2^{-i}, \quad b_i\in\set{0,1}, 
\qquad\Rightarrow\qquad
x = b_0.b_1b_2b_3\dots
\end{equation}
This decomposition is not unique. Indeed, 
\begin{equation}
\label{tnt4}
\sum_{i=1}^\infty 2^{-i} = 1 
\qquad\Rightarrow\qquad
0.1^\infty = 1.0^\infty,
\end{equation}
where the superscript $^\infty$ means that the symbol is repeated
indefinitely. This, however, is the only kind of degeneracy. For our
purposes, the following convention will be useful. We denote by $\cB$ the
set of symbolic sequences $\ub = 0.b_1b_2\dots$ which are {\em not}
terminated by $0^\infty$, except that we include the sequence $0.0^\infty$.
Then $\eqref{tnt3}$ defines a bijection $\ub:[0,1]\to\cB$. In fact, $\cB$ 
is a metric space with the distance
\begin{equation}
\label{tnt5}
d(\ub,\ub') = \sum_{i=1}^\infty (1-\delta_{b_ib'_i}) 2^{-i},
\end{equation}
where $\delta_{ab}=1$ if $a=b$ and $0$ otherwise, and $\ub$ is continuous
in the resulting topology. 

How does the tent map act on binary expansions? We first observe that
$x\leqs1/2$ if and only if $b_1=0$. In this case,
\begin{equation}
\label{tnt6}
g_2(0.0b_2b_3\dots) = 2\sum_{i=2}^\infty b_i 2^{-i} 
= \sum_{i=1}^\infty b_{i+1}2^{-i} = 0.b_2b_3\dots
\end{equation}
Otherwise, we have $x>1/2$, $b_1=1$ and 
\begin{align}
\nonumber
g_2(0.1b_2b_3\dots) &= 2 - 2\Bigbrak{\frac12 + \sum_{i=2}^\infty b_i2^{-i}} 
= 1 - \sum_{i=1}^\infty b_{i+1} 2^{-i} \\
&= \sum_{i=1}^\infty (1-b_{i+1}) 2^{-i} = 0.(1-b_2)(1-b_3)\dots
\label{tnt7}
\end{align}
Hence $g_2$ induces a map $\tau$ from $\cB$ to itself, defined by the following
rules:
\begin{itemiz}
\item	shift all digits of $\ub$ one unit to the left, discarding $b_0$;
\item	if the first digit is $1$, reverse all digits;
\item	replace the sequence $10^\infty$, if present, by $01^\infty$.
\end{itemiz}
We can now analyse the dynamics in $\cB$ instead of $[0,1]$, which is easier
due to the relatively simple form of $\tau$. In fact, $g_2$ has an even
simpler representation. With every $\ub\in\cB$, we associate a sequence
$\ueps=(\eps_1,\eps_2,\dots)$ defined by 
\begin{equation}
\label{tnt8}
\eps_j = (-1)^{b_{j-1}+b_j}, \qquad j=1,2,\dots
\end{equation}
The elements of $\eps$ are $-1$ if adjacent digits of $\ub$ are different,
and $+1$ if they are equal. The set $\Sigma$ of sequences $\ueps$ constructed
in this way consists of $\set{-1,+1}^\N$, from which we exclude sequences
ending with $(+1)^\infty$ and containing an even number of $-1$. The
correspondence $\eqref{tnt8}$ admits an inverse defined by
\begin{equation}
\label{tnt9}
b_0 = 0,\qquad
(-1)^{b_j} = \eps_j(-1)^{b_{j-1}} = \prod_{i=1}^j\eps_i, 
\quad b_j\in\set{0,1}, \quad j\geqs1.
\end{equation}
$\Sigma$ can be endowed with a similar distance as \eqref{tnt5}, and then
\eqref{tnt8} defines a homeomorphism between $\cB$ and $\Sigma$, and, by
composition with $\ub$, a homeomorphism between $[0,1]$ and $\Sigma$. 

Treating separately the cases $b_1=0$ and $b_1=1$, it is easy to see that
the tent map induces a dynamics in $\Sigma$ given by 
\begin{equation}
\label{tnt11}
\sigma: (\eps_1,\eps_2,\eps_3,\dots) \mapsto (\eps_2,\eps_3,\dots)
\end{equation}
which is called the \defwd{shift map}. The sequence $\ueps$ has a very
simple interpretation. Indeed,
\begin{equation}
\label{tnt12}
x\in[0,\tfrac12] 
\quad\Leftrightarrow\quad
b_1=0 
\quad\Leftrightarrow\quad
\eps_1=+1, 
\end{equation}
so that
\begin{equation}
\label{tnt13}
\eps_j=+1 
\quad\Leftrightarrow\quad
\sigma^j(\ueps)_1 = +1
\quad\Leftrightarrow\quad
g_2^j(x)\in[0,\tfrac12]. 
\end{equation}
We have thus shown that the sequence of $\eps_j$ indicates whether the
\nth{j} iterate of $x$ is to the left or to the right of $1/2$, and that
this information is encoded in the binary expansion of $\ub(x)$. The
sequence $\ueps$ associated with $x$ is called its \defwd{itinerary}. 

Let us now examine the different possible orbits.

\begin{itemiz}
\item[1.] {\bf Periodic orbits:} The itinerary of a periodic orbit must be
periodic. Conversely, if an itinerary $\ueps$ is periodic with period
$p\geqs1$, then $\sigma^p(\ueps)=\ueps$, and by bijectivity the
corresponding $x$ is a fixed point of $g_2^p$. Thus $x$ is a point of period
$p$ if and only if its itinerary is of the form $\ueps=A^\infty$, with $A$ a
finite sequence of length $p$. There are exactly $2^p$ such sequences, and
thus $2^p/p$ orbits of period $p$ (the number of orbits of least period $p$
can be a bit smaller). For instance, we obtain again the fixed points of
$g_2$:
\begin{align}
\nonumber
\ueps &= (+1)^\infty &&\Leftrightarrow&
\ub &= 0.0^\infty &&\Leftrightarrow&
x&=0 \\
\ueps &= (-1)^\infty &&\Leftrightarrow&
\ub &= 0.(10)^\infty &&\Leftrightarrow&
x&=\tfrac12 + \tfrac18 + \tfrac1{32} + \dots = \tfrac23. 
\label{tnt14}
\end{align}

\item[2.] {\bf Transient orbits:} Itineraries of the form $\ueps=BA^\infty$
with $B$ a finite sequence of length $q$ correspond to orbits that reach a
periodic orbit after a finite number of steps.
\end{itemiz}

\goodbreak

\begin{exercise}
\label{exo_tent}
Find all periodic orbits of periods up to $3$ of the tent map.
\end{exercise}

Let us call \defwd{eventually periodic} orbits which are either periodic,
or reach a periodic orbit after a finite number of iterations. Their
itineraries are of the form $\ueps=A^\infty$ or $\ueps=BA^\infty$. It is
easy to see that the orbit of $x\in[0,1]$ is eventually periodic if and
only if $x$ is rational.\footnote{If $x$ has an eventually periodic orbit,
its itinerary, and hence its binary expansion become eventually periodic,
which implies that $2^{p+q}x-2^qx\in\Z$ for some $p,q\geqs1$, and thus
$x\in\Q$. Conversely, let $x=n/m\in\Q$. The map $x\mapsto\set{2x}$, where
$\set{\cdot}$ denotes the fractional part, shifts the bits of $\ub(x)$ one
unit to the left. It also maps the set $\set{0,1/m,2/m,\dots,(m-1)/m}$ into
itself. Thus the orbit of $x$ under this map is eventually periodic, and so
are its binary expansion and its itinerary.} $\Q$ being dense in $\R$, the
union of all eventually periodic orbits is dense in $[0,1]$. However, this
set is countable and has zero Lebesgue measure.

\begin{itemiz}
\item[3.] {\bf Chaotic orbits:} All irrational initial conditions, by
contrast, admit itineraries which are not periodic. The corresponding orbits
will typically look quite random. 	
\end{itemiz}
The following properties are direct consequences of the symbolic
representation:
\begin{itemiz}
\item	For every sequence $\ueps\in\Sigma$, there exists an $x\in[0,1]$
with itinerary $\ueps$. We may thus choose the initial condition in such a
way that the orbit passes left and right of $1/2$ in any prescribed
order. 
\item	The dynamics is sensitive to initial conditions. If we only know $x$
with finite precision $\delta$, we are incapable of making any prediction on
its orbit after $n$ iterations, whenever $2^n\geqs1/\delta$.
\item	However, simply by looking whether successive iterates of $x$ lie to
the left or right of $1/2$, we are able to determine the binary
expansion of $x$ (even though $g_2$ is not injective).
\item	There exists a dense orbit. Indeed, choose $x$ in such a way that its
itinerary 
\begin{equation}
\label{tnt15}
\ueps = (\mskip1.5mu\underbrace{\underbrace{+1\vphantom{7pt}},
\underbrace{-1\vphantom{7pt}}}_{\text{period 1}},
\underbrace{\underbrace{+1,+1},\underbrace{+1,-1},
\underbrace{-1,+1},\underbrace{-1,-1}}_{\text{period 2}},
\underbrace{\underbrace{+1,+1,+1},\dots}_{\text{period 3}},\dots)
\end{equation}
contains all possible finite sequences, ordered by increasing length. Given
any $\delta>0$ and $y\in[0,1]$, one can find $n\in\N$ such that
$\abs{g_2^n(x)-y}<\delta$. It suffices to take $n$ in such a way that the
shifted sequence $\sigma^n\ueps$ and the itinerary of $y$ agree for the
first $\intpartplus{\abs{\log\delta}/\log2}$ bits. 
\end{itemiz}

\begin{remark}
\label{rem_measure}
The tent map $g_2$ has other interesting properties, from the point of view
of measure theory. The main property is that its unique invariant measure
which is absolutely continuous with respect to the Lebesgue measure is the
uniform measure. This measure is ergodic, meaning that 
\begin{equation}
\label{tnt16}
\lim_{n\to\infty} \frac1n \sum_{i=1}^n \ph(g_2^n(x)) = \int_0^1 \ph(x)\6x
\end{equation}
for every continuous test function $\ph$, and for Lebesgue-almost all
$x\in[0,1]$. Thus from a probabilistic perspective, almost all orbits will
be uniformly distributed over the interval. 
\end{remark}

Itineraries can be defined for more general maps, just by choosing some
partition of phase space. In general, however, there will be no simple
relation between a point and its itinerary, and the correspondence need not
be one-to-one. Some symbolic sequences may never occur, and several
initial conditions may have the same itinerary, for instance if they are in
the basin of attraction of a stable equilibrium. 

Let us now consider the following variant of the tent map:
\begin{equation}
\label{tnt17}
g_3(x) = 
\begin{cases}
3x & \text{for $x\leqs\frac12$} \\
3-3x & \text{for $x\geqs\frac12$.}
\end{cases}
\end{equation}
Since the interval $[0,1]$ is not invariant, we define $g_3$ on all of $\R$.
We first observe that 
\begin{itemiz}
\item	if $x_0<0$, the orbit of $x_0$ converges to $-\infty$;
\item	if $x_0>1$, then $x_1=g_3(x_0)<0$;
\item	if $x_0\in(1/3,2/3)$, then $x_1=g_3(x_0)>1$ and thus $x_2<0$.
\end{itemiz}
Hence all orbits which leave the interval $[0,1]$ eventually converge to
$-\infty$. One could suspect that all orbits leave the interval after a
certain number of iterations. This is not the case since, for instance,
orbits starting in multiples of $3^{-n}$ reach the fixed point $0$ after
finitely many iterations. But there exists a more subtle nontrivial
invariant set. In order to describe it, we use a ternary (base $3$)
representation of $x\in[0,1]$:
\begin{equation}
\label{tnt18}
x = \sum_{i=0}^\infty b_i 3^{-i}, \quad b_i\in\set{0,1,2}, 
\qquad\Rightarrow\qquad
x = b_0.b_1b_2b_3\dots
\end{equation}
Again, this representation is not unique. We can make it unique by replacing
$10^\infty$ with $02^\infty$ and $12^\infty$ with $20^\infty$ if applicable.
Now we observe that
\begin{align}
\nonumber
g_3(0.0b_2b_3\dots) &= 0.b_2b_3\dots \\
g_3(0.2b_2b_3\dots) &= 0.(2-b_2)(2-b_3)\dots 
\label{tnt19}
\end{align}
Points of the form $0.1b_2b_3\dots$ belong to $(1/3,2/3)$, and we already
know that these leave the interval $[0,1]$. It is thus immediate that the
orbit of $x$ never leaves the interval $[0,1]$ if and only if its ternary
expansion $\ub(x)$ does not contain the symbol $1$. The largest invariant
subset of $[0,1]$ is thus 
\begin{equation}
\label{tnt20}
\Lambda = \bigsetsuch{x\in[0,1]}{b_i(x)\neq 1\;\forall i\geqs 1}.
\end{equation}
$\Lambda$ is obtained by removing from $[0,1]$ the open intervals
$(1/3,2/3)$, $(1/9,2/9)$, $(7/9,8/9)$, $(1/27,2/27)$ and so on
(\figref{fig_Cantor}). The resulting set is called a \defwd{Cantor set}:

\begin{figure}
 \centerline{\psfig{figure=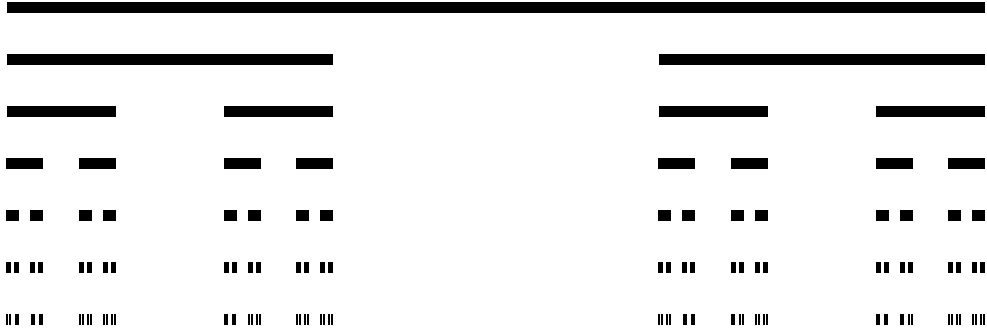,width=130mm,height=45mm}}
 \captionspace
 \caption[]{The first steps of the construction of the Cantor set $\Lambda$.}
 \label{fig_Cantor}
\end{figure}

\begin{definition}
\label{def_Cantor}
A set $\Lambda$ is called a \defwd{Cantor set} if it is closed, its interior
is empty, and all its points are accumulation points.
\end{definition}

\goodbreak

\begin{exercise}
\label{exo_Cantor}
Prove that $\Lambda$ defined in \eqref{tnt20} is a Cantor set. 
\end{exercise}
The Cantor set $\Lambda$ is an example of \defwd{fractal}. 
\begin{definition}
\label{def_boxdim}
Let $\cM$ be a subset of $\R^d$. Assume that for any $\eps>0$, $\cM$ can be
covered by a finite number of hypercubes of side length $\eps$. Let
$N(\eps)$ be the smallest possible number of such cubes. Then the
\defwd{box-counting dimension} of $\cM$ is defined by 
\begin{equation}
\label{tnt21}
D_0(\cM) = \lim_{\eps\to0} \frac{\log N(\eps)}{\log(1/\eps)}.
\end{equation}
\end{definition}
One easily shows that for \lq\lq usual\rq\rq\ sets $\cM\subset\R^d$, such
as a $d$-dimensional hypercube, $D_0=d$. For the Cantor set \eqref{tnt20},
however, we find that $N(\eps)=2^n$ if $3^{-n}\leqs\eps<3^{-n+1}$, and
thus
\begin{equation}
\label{tnt22}
D_0(\Lambda) = \frac{\log2}{\log3}.
\end{equation}
The map $g_3$ restricted to $\Lambda$ is conjugate to the tent map $g_2$.
Indeed, with any point $x=0.b_1b_2b_3\dots$ in $\Lambda$, we can associate
$h(x) = 0.(b_1/2)(b_2/2)(b_3/2)\dots$ in $[0,1]$, and the relations
\eqref{tnt19} and \eqref{tnt6}, \eqref{tnt7} show that 
\begin{equation}
\label{tnt23}
h\circ g_3\evalat{\Lambda} \circ h^{-1} = g_2.
\end{equation}
Thus we can easily compute itineraries of points $x\in\Lambda$, and the
conclusions on the behaviour of orbits of $g_2$ can be transposed to
$g_3\evalat{\Lambda}$. $\Lambda$ is the simplest example of what is called
a \defwd{hyperbolic invariant set}. 

Similar properties can be seen to hold for nonlinear perturbations of the
tent map. Indeed, the ternary representation is a useful tool, but not
essential for the existence of orbits with all possible symbolic
representations. It is, in fact, sufficient that $f$ maps two disjoint
intervals $I_1, I_2\subset[0,1]$ onto $[0,1]$, and maps all other points
outside $[0,1]$. Then the symbolic representation corresponds to the
sequence of intervals $I_1$ or $I_2$ visited by the orbit. 


\goodbreak
\subsection{Homoclinic Tangles and Smale's Horseshoe Map}
\label{ssec_smale}

\begin{figure}
 \centerline{\psfig{figure=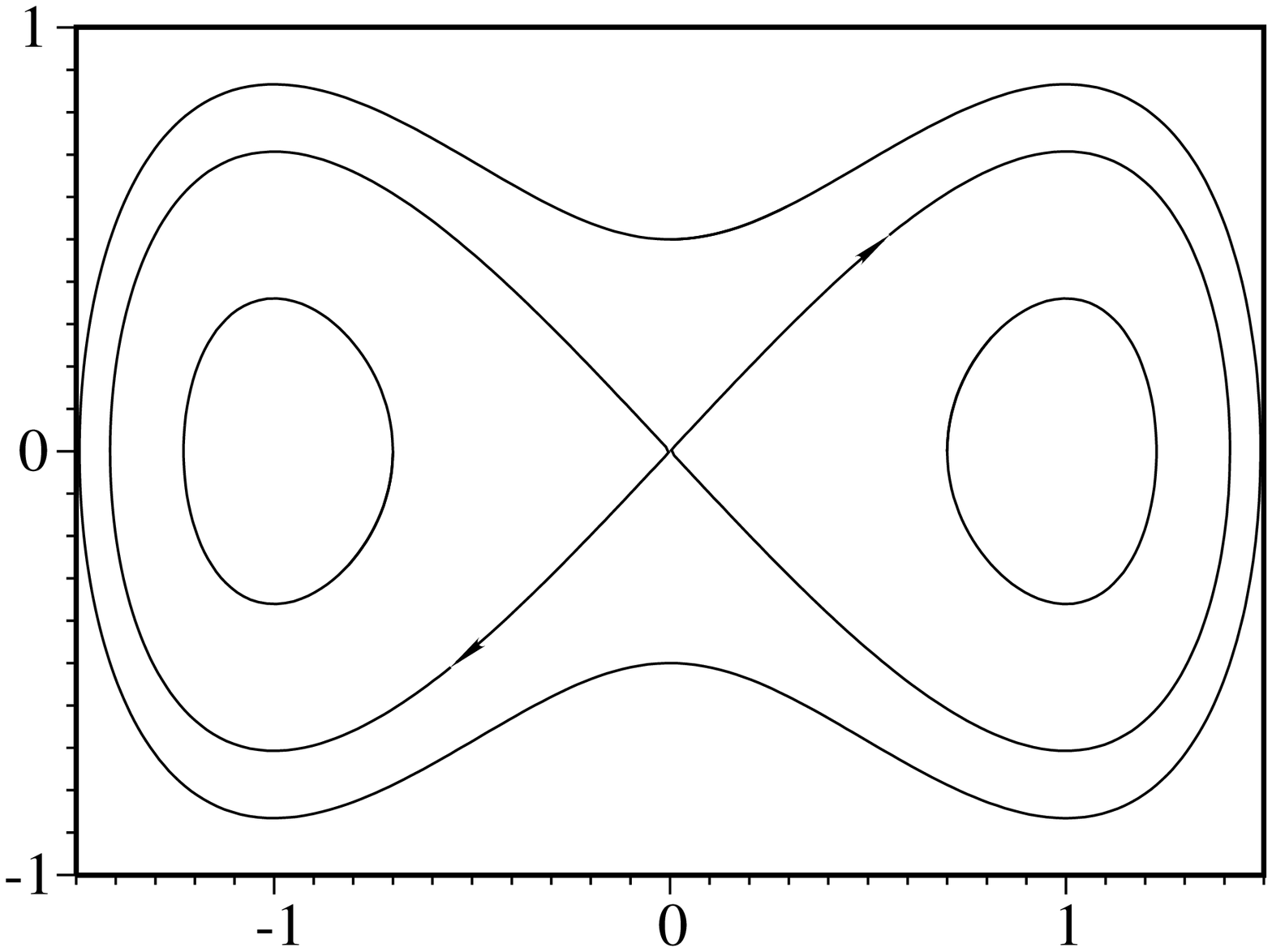,height=50mm}
 \hspace{5mm}
 \psfig{figure=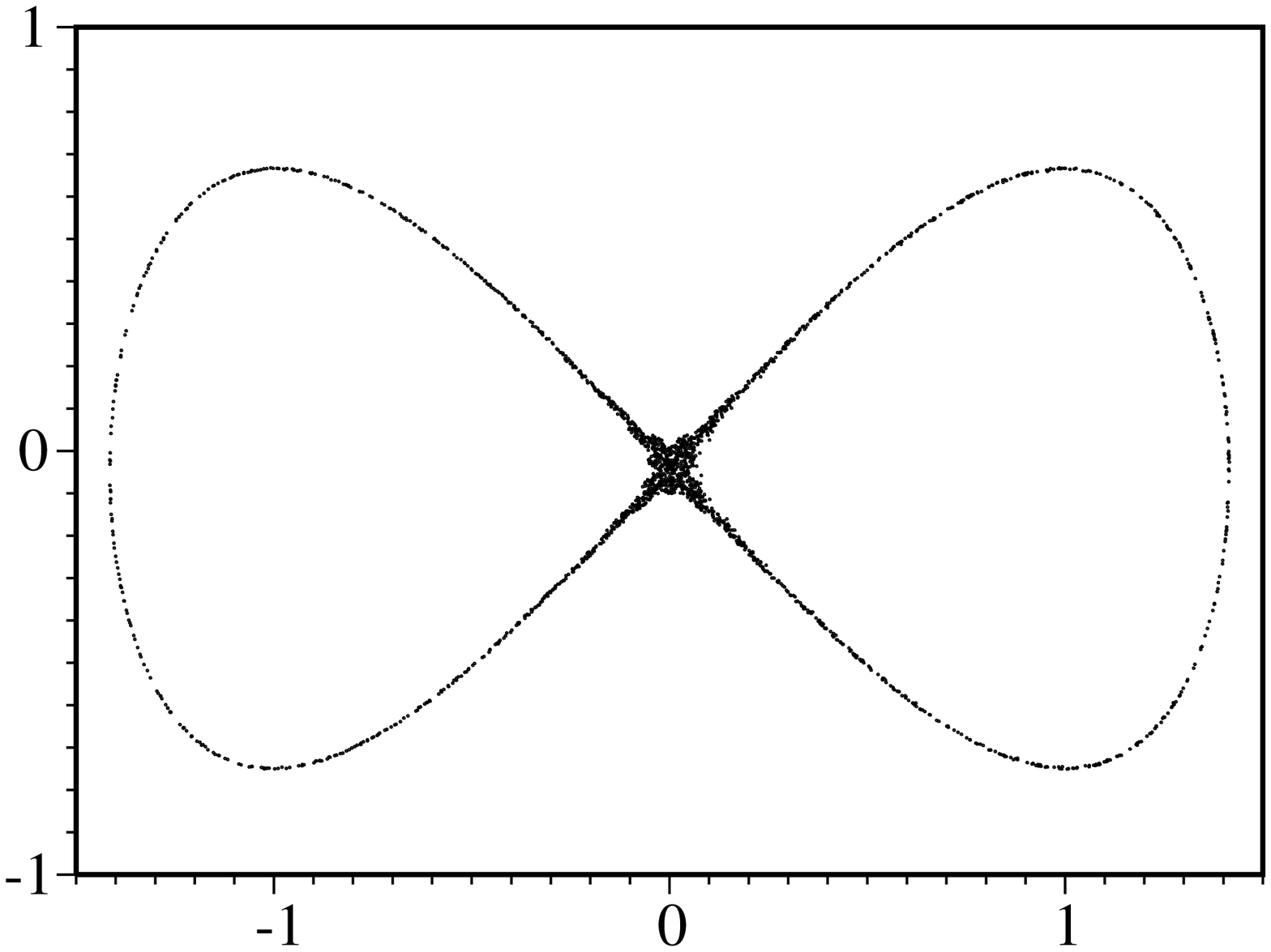,height=50mm}}
 \captionspace
 \caption[]{(a) Orbits of Duffing's equation for $\eps=0$, (b) one orbit of
 the Poincar\'e map obtained for $\eps=0.25$.}
 \label{fig_duffing}
\end{figure}

One-dimensional maps are rather special cases of dynamical systems, but it
turns out that some of their properties can often be transposed to more
\lq\lq realistic\rq\rq\ systems. As a motivating example, let us consider
the equation
\begin{equation}
\label{sh1}
\begin{split}
\dot x_1 &= x_2 \\
\dot x_2 &= x_1 - x_1^3,
\end{split}
\end{equation}
called the (undamped) unforced \defwd{Duffing oscillator}. This is a
Hamiltonian system, with Hamiltonian
\begin{equation}
\label{sh2}
H(x_1,x_2) = \frac12 x_2^2 + \frac14 x_1^4 - \frac12 x_1^2.
\end{equation}
$H$ is a constant of the motion, and thus orbits of \eqref{sh1} belong to
level curves of $H$ (\figref{fig_duffing}a). The point $(0,0)$ is a
hyperbolic equilibrium, with the particularity that its unstable and stable
manifolds $\Wglo{u}$ and $\Wglo{s}$ are identical: they form the curve
$H=0$, and are called \defwd{homoclinic loops}. 

Let us now perturb \eqref{sh1} by a periodic forcing,
\begin{equation}
\label{sh3}
\begin{split}
\dot x_1 &= x_2 \\
\dot x_2 &= x_1 - x_1^3 + \eps\sin(2\pi t). 
\end{split}
\end{equation}
We can introduce $x_3=t$ as additional variable to obtain the autonomous
system
\begin{equation}
\label{sh4}
\begin{split}
\dot x_1 &= x_2 \\
\dot x_2 &= x_1 - x_1^3 + \eps\sin(2\pi x_3) \\
\dot x_3 &= 1. 
\end{split}
\end{equation}
Here $x_3\in\fS^1$ should be considered as a periodic variable. We can thus
use the surface $x_3=0$ as a Poincar\'e section. \figref{fig_duffing}b shows
an orbit of the Poincar\'e map, which has a complicated structure: it seems
to fill a two-dimensional region, in which its dynamics is quite random. 

This phenomenon can be explained at least partially, by analysing the
properties of the Poincar\'e map $P_\eps$. First observe that the system
\eqref{sh4} is conservative. The section being perpendicular to the flow,
the Poincar\'e map is also conservative. When $\eps=0$, $P_0=\ph_1$ is the
time-one flow of \eqref{sh1}, and thus $P_0$ admits the origin as
hyperbolic fixed point. The implicit function theorem implies that for
small $\eps$, $P_\eps$ admits an isolated hyperbolic fixed point
$x^\star(\eps)$ near $(0,0)$. It is unlikely that the stable and unstable
manifolds of $x^\star(\eps)$ still form a loop when $\eps>0$, and if not,
they must intersect transversally at some point $x_0$ because of area
conservation. $x_0$ is called a \defwd{homoclinic point}. A method due to
Melnikov allows to prove that such a transverse intersection indeed exists.

\begin{figure}
 \centerline{
 \psfig{figure=figs/homocline.ps,width=50mm}
 \psfig{figure=figs/homocline2.ps,width=65mm}
 }
 \figtext{
 	\writefig	0.9	5.0	{a}
 	\writefig	7.9	5.0	{b}
 	\writefig	1.9	1.1	{$x^\star$}
 	\writefig	5.5	4.8	{$x_0$}
 	\writefig	4.4	3.6	{$P(x_0)$}
 	\writefig	5.2	2.8	{$P^2(x_0)$}
 	\writefig	1.8	2.3	{$\Wglo{u}$}
 	\writefig	2.9	1.1	{$\Wglo{s}$}
 	\writefig	6.6	5.0	{$A_1$}
 	\writefig	2.9	2.3	{$A_2$}
 	\writefig	7.8	1.8	{$x^\star$}
 }
 \captionspace
 \caption[]{If the stable and unstable manifolds of a conservative map
 intersect transversally, they must intersect infinitely often, forming a
 homoclinic tangle.}
 \label{fig_homoclinic}
\end{figure}

Consider now the successive images $x_n=P_\eps^n(x_0)$. Since $x_0$ belongs
to the stable manifold $\Wglo{s}$ of $x^\star(\eps)$, all $x_n$ must also
belong to $\Wglo{s}$, and they must accumulate at $x^\star$ for
$n\to\infty$. Similarly, since $x_0$ belongs to the unstable manifold
$\Wglo{u}$, all $x_n$ must  belong to $\Wglo{u}$ and accumulate at
$x^\star$ for $n\to-\infty$. Thus $\Wglo{s}$ and $\Wglo{u}$ must intersect
infinitely often. Because the map is area preserving, the area between
$\Wglo{s}$, $\Wglo{u}$, and any consecutive intersection points must be the
same, and thus $\Wglo{u}$ has to oscillate with increasing amplitude when
approaching $x^\star$ (\figref{fig_homoclinic}). $\Wglo{s}$ has a similar
behaviour. This complicated geometrical structure, which was first
described by Poincar\'e, is called the \defwd{homoclinic tangle}. 

\begin{figure}
 \centerline{\psfig{figure=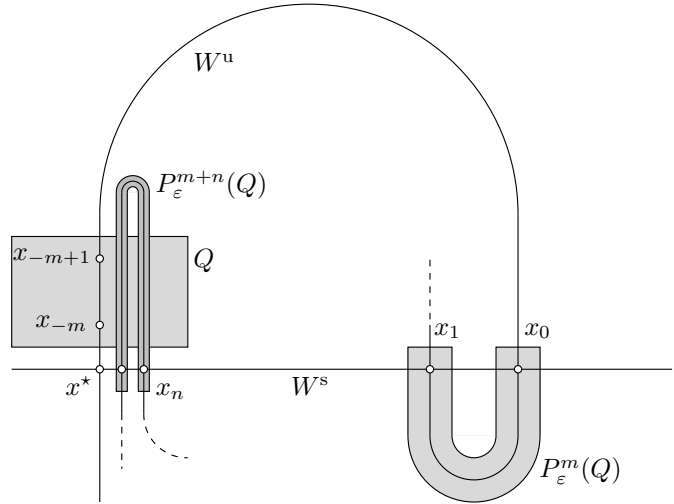,height=70mm}}
 \figtext{
 	\writefig	3.7	2.1	$x^\star$
 	\writefig	3.3	3.0	$x_{-m}$
 	\writefig	3.0	3.9	$x_{-m+1}$
 	\writefig	5.4	3.8	$Q$
 	\writefig	4.9	2.1	$x_n$
 	\writefig	6.7	2.1	$\Wglo{s}$
 	\writefig	8.6	2.9	$x_1$
 	\writefig	9.8	2.9	$x_0$
 	\writefig	4.9	4.8	$P_\eps^{m+n}(Q)$
 	\writefig	10.0	1.0	$P_\eps^m(Q)$
 	\writefig	5.4	6.4	$\Wglo{u}$
 }
 \vspace{1mm}
 \captionspace
 \caption[]{Schematic representation of the homoclinic tangle. The
 homoclinic orbit $\set{x_n}_{n\in\Z}$ is asymptotic to $x^\star$ for
 $n\to\pm\infty$. One can choose $m$ and $n$ in such a way that the image of
 a rectangle $Q$ after $m+n$ iterations intersects $Q$ twice.}
 \vspace{4mm}
 \label{fig_horseshoe1}
\end{figure}

The dynamics near the homoclinic tangle can be described as follows. Let
$Q$ be a small rectangle, containing a piece of unstable manifold
$\Wglo{u}$ and two points $x_{-m}$, $x_{-m+1}$ of the homoclinic orbit
(\figref{fig_horseshoe1}). At least during the first iterations of
$P_\eps$, $Q$ will be stretched in the unstable direction, and contracted
in the stable one. After $m$ iterations, the image $P_\eps^m(Q)$ contains
$x_0$ and $x_1$. Let $n$ be sufficiently large that the points $x_n$ and
$x_{n+1}$ of the homoclinic orbit are close to $x^\star$. Due to area
conservation, one can arrange that the piece of $\Wglo{u}$ between $x_n$
and $x_{n+1}$ crosses $Q$ at least twice, and so does the image
$P_\eps^{m+n}(Q)$. 

\begin{figure}
 \centerline{\psfig{figure=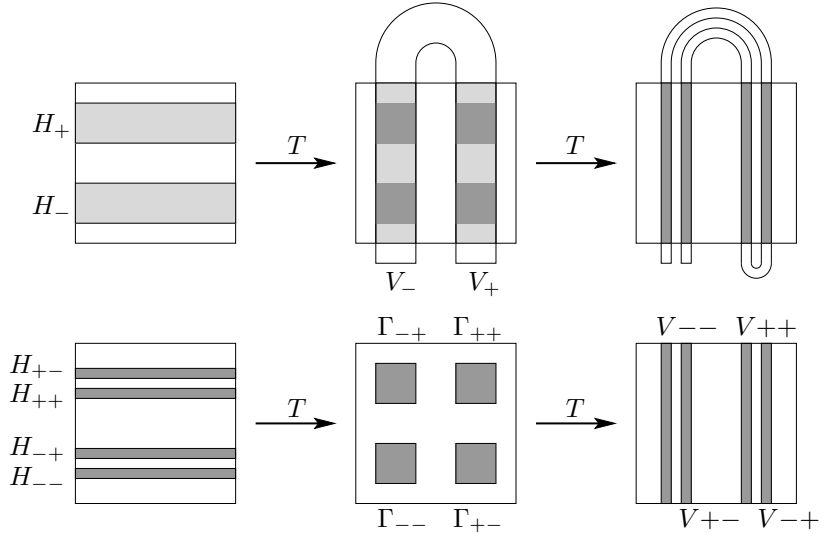,height=70mm}}
 \figtext{
 	\writefig	2.0	5.6	$H_+$
 	\writefig	2.0	4.5	$H_-$
 	\writefig	6.7	3.5	$V_-$
 	\writefig	7.8	3.5	$V_+$
 	\writefig	1.7	2.4	$H_{+-}$
 	\writefig	1.7	2.0	$H_{++}$
 	\writefig	1.7	1.3	$H_{-+}$
 	\writefig	1.7	0.9	$H_{--}$
 	\writefig	6.6	2.9	$\Gamma_{-+}$
 	\writefig	7.6	2.9	$\Gamma_{++}$
 	\writefig	6.6	0.35	$\Gamma_{--}$
 	\writefig	7.6	0.35	$\Gamma_{+-}$
 	\writefig	10.3	2.9	$V{--}$
 	\writefig	11.35	2.9	$V{++}$
 	\writefig	10.6	0.35	$V{+-}$
 	\writefig	11.65	0.35	$V{-+}$
 	\writefig	5.4	5.3	$T$
 	\writefig	5.4	1.8	$T$
 	\writefig	9.1	5.3	$T$
 	\writefig	9.1	1.8	$T$
 }
 \vspace{6mm}
 \caption[]{The horseshoe map $T$ maps the rectangles $H_{\pm}$ of the unit
 square $Q$ to the rectangles $V_{\pm}$. If
 $\Gamma_{\eps_0\eps_1}=V_{\eps_0}\cap H_{\eps_1}$, one also has
 $T(H_{\eps_0\eps_1})=\Gamma_{\eps_0\eps_1}$ and
 $T(\Gamma_{\eps_0\eps_1})=V_{\eps_0\eps_1}$.}
 \label{fig_horseshoe2}
\end{figure}

This behaviour is reproduced qualitatively by the \defwd{Smale horseshoe
map} (\figref{fig_horseshoe2}). This map $T$ takes a square
$[0,1]\times[0,1]$, stretches it vertically by a factor $\mu>2$, and
contracts it horizontally by a factor $\lambda<1/2$. Then it bends the
resulting rectangle in the shape of a horseshoe, and superimposes it with
the initial square. Two horizontal strips $H_+$ and $H_-$, of size
$1\times\mu^{-1}$, are mapped, respectively, to two vertical strips $V_+$
and $V_-$, of size $\lambda\times1$. For simplicity, the map restricted to
$V_+\cup V_-$ is assumed to be linear, but its qualitative features can be
shown to remain unchanged by small nonlinear perturbations. 

The aim is now to construct an invariant set $\Lambda$ of $T$. To this end,
we observe that $T(x)\in Q$ if and only if $x\in H_-\cup H_+$. Since
$T(x)\in V_-\cup V_+$, $T^2(x)\in Q$ if and only if $T(x)$ belongs to one
of the four rectangles $\Gamma_{\eps_0\eps_1} = V_{\eps_0}\cap H_{\eps_1}$,
where $\eps_0, \eps_1 \in\set{-,+}$. These rectangles have size
$\lambda\times\mu^{-1}$. The preimage of each $\Gamma_{\eps_0\eps_1}$ is a
rectangle $H_{\eps_0\eps_1}\subset H_{\eps_0}$, of size $1\times\mu^{-2}$.
The image of $\Gamma_{\eps_0\eps_1}$ is a rectangle
$V_{\eps_0\eps_1}\subset V_{\eps_1}$, of size $\lambda^2\times1$. Thus
$T^2$ maps $H_{\eps_0\eps_1}$ to $V_{\eps_0\eps_1}$. 

More generally, we can define
\begin{equation}
\label{sh5}
H_{\eps_0\eps_1\dots\eps_n} = 
\bigsetsuch{x\in Q}{x\in H_{\eps_0}, T(x)\in H_{\eps_1},\dots,T^n(x)\in
H_{\eps_n}}.
\end{equation}
Note that $H_{\eps_0\eps_1\dots\eps_n}\subset H_{\eps_0\dots\eps_{n-1}}$,
and that
\begin{equation}
\label{sh6}
H_{\eps_0\eps_1\dots\eps_n} = \bigsetsuch{x\in H_{\eps_0}}{T(x)\in
H_{\eps_1\dots\eps_n}}.
\end{equation}
By induction, it is straightforward to show that
$H_{\eps_0\eps_1\dots\eps_n}$ is a rectangle of size $1\times\mu^{-n}$.
Similarly, we introduce
\begin{align}
\nonumber
V_{\eps_{-m}\dots\eps_{-1}} &= 
\bigsetsuch{x\in Q}{T^{-1}(x)\in H_{\eps_{-1}},\dots,T^{-m}(x)\in
H_{\eps_{-m}}} \\
&= \bigsetsuch{x\in V_{\eps_{-1}}}{T^{-1}(x)\in
V_{\eps_{-m}\dots\eps_{-2}}},
\label{sh7}
\end{align}
which is a rectangle of size $\lambda^m\times 1$ contained in
$V_{\eps_{-m+1}\dots\eps_{-1}}$ . It follows that 
\begin{align}
\nonumber
\Gamma_{\eps_{-m}\dots\eps_{-1},\eps_0\eps_1\dots\eps_n} 
&= V_{\eps_{-m}\dots\eps_{-1}}\cap H_{\eps_0\eps_1\dots\eps_n} \\
&= \bigsetsuch{x\in Q}{T^j(x)\in H_{\eps_j},-m\leqs j\leqs n}
\label{sh8}
\end{align}
is a rectangle of size $\lambda^m\times\mu^{-n}$. By construction,
$T^j(x)\in Q$ for $-m\leqs j\leqs n+1$ if and only if
$x\in\Gamma_{\eps_{-m}\dots,\eps_0\dots\eps_n}$ for some sequence 
$\eps_{-m}\dots,\eps_0\dots\eps_n$ of symbols $-1$ and $+1$. 

\goodbreak

Now let $n, m$ go to $\infty$. Let $\Sigma$ be the set of
bi-infinite sequences $\ueps=\dots\eps_{-1},\eps_0\eps_1\dots$ of symbols
$-1$ and $+1$. It is a metric space for the distance
\begin{equation}
\label{sh9}
d(\ueps,\ueps') = \sum_{i=-\infty}^\infty (1-\delta_{\eps_i\eps'_i})
2^{-\abs{i}}.
\end{equation}
From \eqref{sh8} we conclude that the largest subset of $Q$ invariant under
$T$ is 
\begin{equation}
\label{sh10}
\Lambda = \bigcap_{n=-\infty}^{+\infty} T^n(Q) 
= \bigsetsuch{\Gamma_{\ueps}}{\ueps\in\Sigma},
\end{equation}
where $\Gamma_{\ueps}$, defined as the limit of \eqref{sh8} when the finite
sequence $\eps_{-m}\dots,\eps_0\dots\eps_n$ converges to $\ueps$, is a
single point. $\Lambda$ is a Cantor set, obtained by taking the product of
two one-dimensional Cantor sets. The map $\phi:\ueps\mapsto\Gamma_{\ueps}$
is a bijection from $\Sigma$ to $\Lambda$, which is continuous in the
topology defined by \eqref{sh9}. It follows from \eqref{sh8} that $T$ is
conjugated to the shift map
\begin{equation}
\label{sh11}
\sigma = \phi^{-1}\circ T\evalat{\Lambda}\circ\phi: 
\dots\eps_{-2}\eps_{-1},\eps_0\eps_1\ldots \mapsto
\dots\eps_{-1}\eps_0,\eps_1\eps_2\dots
\end{equation}
This conjugacy can be used to describe the various orbits of
$T\evalat{\Lambda}$ in a similar way as we did for the tent map. In
particular, periodic itineraries $\ueps$ correspond to periodic orbits of
$T$, aperiodic itineraries correspond to chaotic orbits. Smale has shown
that these qualitative properties are robust under small perturbations of
$T$:

\begin{figure}
 \centerline{\psfig{figure=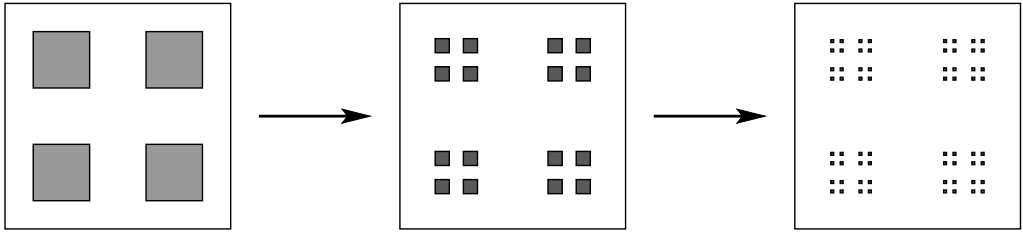,height=35mm}}
 \figtext{
 }
 \captionspace
 \caption[]{Successive approximations of the invariant set $\Lambda$.}
 \label{fig_horseshoe3}
\end{figure}

\begin{theorem}[Smale]
\label{thm_Smale}
The horseshoe map $T$ has an invariant Cantor set $\Lambda$ such that 
\begin{itemiz}
\item	$\Lambda$ contains a countable set of periodic orbits of arbitrarily
long periods;
\item	$\Lambda$ contains an uncountable set of bounded nonperiodic orbits;
\item	$\Lambda$ contains a dense orbit.
\end{itemiz}
Moreover, any map $\widetilde T$ sufficiently close to $T$ in the $\cC^1$
topology has an invariant Cantor set $\widetilde\Lambda$ with $\widetilde
T\evalat{\widetilde \Lambda}$ topologically equivalent to
$T\evalat{\Lambda}$. 
\end{theorem}

The example of the horseshoe map can be generalised to a class of so-called
\defwd{axiom A} systems with similar chaotic properties. These systems can
be used to show the existence of chaotic orbits in a large class of
dynamical systems, including (but not limited to) those with a transverse
homoclinic intersection. 

Note that the examples of hyperbolic invariant sets that we have
encountered have measure zero and are not attracting. Showing the existence
of invariant sets containing chaotic orbits and attracting nearby orbits is
a much more difficult task. 


\newpage
\section{Strange Attractors}
\label{sec_sa}


\subsection{Attracting Sets and Attractors}
\label{ssec_asa}

We consider in this section a dynamical system on $\cD\subset\R^n$, defined
by a flow $\ph_t$. One can include the case of iterated maps by restricting
$t$ to integer values and setting $\ph_t=F^t$. Various subsets of $\cD$
can be associated with the flow.

\begin{definition}\hfill
\label{def_invset}
\begin{itemiz}
\item	A subset $\cS\subset\cD$ is called \defwd{invariant} if
$\ph_t(\cS)=\cS$ for all $t$.
\item	A subset $\cA\subset\cD$ is called an \defwd{attracting set} if
there exists a \nbh\ $\cU$ of $\cA$ such that for all $x\in\cU$,
$\ph_t(x)\in\cU$ for all $t\geqs0$ and $\ph_t(x)\to\cA$ as $t\to\infty$. 
\item	A point $x$ is \defwd{nonwandering} for $\ph_t$ if for every \nbh\
$\cU$ of $x$ and every $T>0$, there exists a $t>T$ such that
$\ph_t(\cU)\cap\cU\neq\emptyset$. The set of all nonwandering points is the
\defwd{nonwandering set} $\Omega$.
\end{itemiz}
\end{definition}

These sets can be partly determined by looking at the asymptotic behaviour
of various orbits of the flow. 

\begin{definition}
\label{def_limitset}
The $\w$-limit set of $x$ for $\ph_t$ is the set of $y\in\cD$ such that
there exists a sequence $t_n\to\infty$ with $\ph_{t_n}(x)\to y$. 
The $\alpha$-limit set of $x$ for $\ph_t$ is the set of $y\in\cD$ such that
there exists a sequence $t_n\to-\infty$ with $\ph_{t_n}(x)\to y$. 
\end{definition}

The $\alpha$- and $\w$-limit sets of any $x$ are invariant sets, and the
$\w$-limit set is included in the nonwandering set. One can show that
asymptotically stable equilibrium points, periodic orbits, and invariant tori
are all $\w$-limit sets of the orbits in their basin of attraction,
nonwandering sets and attracting sets. There is, however, a problem with the
definition of attracting set, as shows the following example.

\begin{example}
\label{ex_attractingset}
Consider the differential equation
\begin{equation}
\label{asa1}
\begin{split}
\dot x_1 &= x_1 - x_1^3 \\
\dot x_2 &= -x_2.
\end{split}
\end{equation}
The $\w$-limit set of $x$ is $(-1,0)$ if $x_1<0$, $(0,0)$ if $x_1=0$ and
$(1,0)$ if $x_1>0$. $(0,0)$ is also the $\alpha$-limit set of all points in
$(-1,1)\times\set{0}$. The nonwandering set is composed of the three 
equilibrium points of the flow, while  $(\pm 1,0)$ are  attracting sets.
However, the segment $[-1,1]\times\set{0}$ is also an attracting set.
\end{example}

One would like to exclude attracting sets such as the segment
$[-1,1]\times\set{0}$ in the example, which contains wandering points. This
is generally solved in the following way.

\begin{definition}
\label{def_attractor}
A closed invariant set $\Lambda$ is \defwd{topologically transitive} if
$\ph_t$ has an orbit which is dense in $\Lambda$. An \defwd{attractor} is a
topologically transitive attracting set. 
\end{definition}

\begin{remark}
\label{rem_chainrecurrence}
Sometimes, one uses a weaker notion of indecomposability than topological
transitivity, based on the notion of \defwd{chain recurrence}. In this
case, one requires that for any pair of points $x,y$ and any $\eps>0$, there
exist points $x_0=x,x_1,\dots,x_n=y$ and times $t_1,\dots,t_n$ such that
$\norm{\ph_{t_j}(x_{j-1})-x_j}\leqs\eps$ for all $j$.
\end{remark}

\begin{remark}
\label{rem_attractor}
A slightly weaker definition of attractor is proposed in \cite{GH}: An
\defwd{attractor} is an indecomposable closed invariant set $\Lambda$ with
the property that, given $\eps>0$, there is a set $\cU$ of positive
Lebesgue measure in the $\eps$-\nbh\ of $\Lambda$ such that, if $x\in\cU$,
the forward orbit of $x$ is contained in $\cU$ and the $\w$-limit set of
$x$ is contained in $\Lambda$.
\end{remark}


\subsection{Sensitive Dependence on Initial Conditions}
\label{ssec_sdi}

A strange attractor is, basically, an attractor on which the dynamics is
chaotic. Obviously, this requires that we define what we mean by \lq\lq
chaotic\rq\rq. In \cite{GH}, for instance, a strange attractor is defined
as an attractor containing a transversal homoclinic orbit. As we saw in
section \ref{ssec_smale}, the existence of such an orbit implies various
chaotic properties (it is not necessary to assume that the system is
conservative). Modern definitions are a bit less specific, they require
that the dynamics be sensitive to initial conditions:

\begin{definition}
\label{def_sdic}
Let $\Lambda$ be a compact set such that $\ph_t(\Lambda)\subset\Lambda$ for
all $t\geqs 0$. 
\begin{itemiz}
\item	The flow $\ph_t$ is said to have \defwd{sensitive dependence
on initial conditions} on $\Lambda$ if there exists $\eps>0$ with the
following property: For any $x\in\Lambda$ and any \nbh\ $\cU$ of $x$, there
exists $y\in\cU$ and $t>0$ such that $\norm{\ph_t(x)-\ph_t(y)}\geqs\eps$. 
\item	$\Lambda$ is a \defwd{strange attractor} if it is an attractor and
$\ph_t$ has sensitive dependence on initial conditions on $\Lambda$. 
\end{itemiz}
\end{definition}

\begin{exercise}
\label{exo_sdic}
Show that the tent map and Smale's horseshoe map have sensitive
dependence on initial conditions.
\end{exercise}

Sensitive dependence on initial conditions requires that for any point $x$,
one can find an arbitrarily close point $y$ such that the orbits of $x$ and
$y$ diverge from each other (the quantity $\eps$ should not depend on
$\norm{x-y}$, though $t$ may depend on it). This is a rather weak property,
and one often requires that the divergence occur at an exponential rate.

\begin{definition}\hfill
\label{def_Liapexp}
\begin{itemiz}
\item
Assume $\ph_t$ is the flow of a differential equation $\dot x=f(x)$. Let
$x\in\cD$ and let $U(t)$ be the principal solution of the equation
linearized around the orbit of $x$, 
\begin{equation}
\label{sdic1}
\dot y = \dpar fx\bigpar{\ph_t(x)} y.
\end{equation}
Oseledec \cite{Oseledec} proved that under quite weak assumptions on $f$,
the limit 
\begin{equation}
\label{sdic2}
L = \lim_{t\to\infty} \frac1{2t} \log\bigpar{\transpose{U(t)}\mskip1.5mu U(t)} 
\end{equation}
exists. The eigenvalues of $L$ are called the \defwd{Liapunov exponents} of
the orbit of $x$. 
\item	If $F$ is an iterated map and $\set{x_k}$ is a given orbit of $F$,
consider the linear equation
\begin{equation}
\label{sdic3}
y_{k+1} = \dpar Fx(x_k) y_k 
\qquad\Rightarrow\qquad
y_k = U_k y_0, \quad
U_k = \dpar Fx(x_{k-1})\dots \dpar Fx(x_0).
\end{equation}
The \defwd{Liapunov exponents} of the orbit $\set{x_k}$ are the eigenvalues
of 
\begin{equation}
\label{sdic4}
L = \lim_{k\to\infty} \frac1{2k} \log\bigpar{\transpose{U_k}\mskip1.5mu U_k}. 
\end{equation}
\end{itemiz}
\end{definition}

Note that $\transpose{U(t)}\mskip1.5mu U(t)$ is a symmetric positive
definite matrix, and hence it is always diagonalizable and has real
eigenvalues.  This definition implies that the solution $y(t)$ of the
equation linearized around the particular solution $\ph_t(x)$ satisfies 
\begin{equation}
\label{sdic5}
\norm{y(t)}^2 = \norm{U(t)y(0)}^2 
= \pscal{y(0)}{\transpose{U(t)}\mskip1.5mu U(t)y(0)} 
\simeq \pscal{y(0)}{\e^{2tL}y(0)}. 
\end{equation}
Let $\lambda_1$ be the largest eigenvalue of $L$. Unless the projection of
$y(0)$ on the eigenspace of $L$ associated with $\lambda_1$ is zero,
$\norm{y(t)}$ will grow asymptotically like $\e^{\lambda_1t}$. We thus say
that the flow has \defwd{exponentially sensitive dependence on initial
conditions} in $\Lambda$ if the largest Liapunov exponent of all orbits in
$\Lambda$ is positive.

\begin{prop}\hfill
\label{prop_sdic}
\begin{itemiz}
\item	If $\ph_t$ is conservative, then the sum of all Liapunov
exponents of any orbit is zero. 
\item	If the flow is dissipative, then this sum is negative.
\item	Assume $\set{\ph_t(x)}$ is an orbit of the differential equation
$\dot x=f(x)$, such that $\norm f$ is bounded below and above by strictly
positive constants on $\set{\ph_t(x)}$. Then this orbit has at least one
Liapunov exponent equal to zero. 
\end{itemiz}
\end{prop}
\begin{proof}
We saw in the proof of Proposition \ref{prop_vol2} that the determinant of
$U(t)$ is constant if the system is conservative, and decreasing if the
system is dissipative. By definition, $\det U(0)=1$ and by uniqueness of
solutions, $\det U(t)>0$ for all $t$. Thus for $t>0$, the product of all
eigenvalues of $U(t)$ is equal to $1$ in the conservative case, and belongs
to $(0,1)$ in the dissipative case. The same is true for
$\det\transpose{U(t)}\mskip1.5mu U(t)$. But for any matrix $B$,
$\det\e^B=\e^{\Tr B}$ because the eigenvalues of $\e^B$ are exponentials of
the eigenvalues of $B$. 

Assume now $x(t)=\ph_t(x)$ is a solution of $\dot x=f(x)$. Then 
\[
\dtot{}{t} \dot x(t) = \dtot {}{t} f(\ph_t(x)) 
= \dpar fx(\ph_t(x)) \dot x(t),
\]
and thus $\dot x(t)=U(t)\dot x(0)$ by definition of $U(t)$. Hence,
\[
\pscal{\dot x(0)}{\transpose{U(t)}\mskip1.5mu U(t) \dot x(0)} 
= \norm{U(t) \dot x(0)}^2 = \norm{\dot x(t)}^2 
= \norm{f(\ph_t(x))}^2.
\]
By the assumption on $\norm{f}$, it follows that as $t\to\infty$, the
function
\[
\pscal{\dot x(0)}{\e^{2tL} \dot x(0)}
\]
is bounded above and below by strictly positive constants. Being symmetric,
$L$ admits an orthonormal set of eigenvectors $u_1,\dots,u_n$, that is, they
satisfy $\pscal{u_i}{u_j} = \delta_{ij}$. If $c_j=\pscal{u_j}{\dot x(0)}$,
then $\dot x(0)=\sum_j c_j u_j$ and thus
\[
\pscal{\dot x(0)}{\e^{2tL} \dot x(0)} = \sum_{j=1}^n c_j^2 \e^{2\lambda_jt}.
\]
Since $\dot x(0)\neq 0$, at least one of the coefficients, say $c_k$, is
different from zero, and thus $\lambda_k$ must be equal to zero.
\end{proof}

It is easy to find systems with positive Liapunov exponents. Consider for
instance the linear system $\dot x = Ax$. In this case, $U(t)=\e^{At}$ and
one shows (using, for instance, the Jordan canonical form of $A$), that the
Liapunov exponents are exactly the real parts of the eigenvalues of $A$.
Thus if $A$ has an eigenvalue with positive real part, all orbits have
exponentially sensitive dependence on initial conditions. 

More generally, if $x^\star$ is a linearly unstable equilibrium point, its
largest Liapunov exponent will be positive. Similarly, if $\Gamma$ is a
periodic orbit, then we have seen in Theorem \ref{thm_floquet} that
$U(t)=P(t)\e^{Bt}$, where $P(t)$ is periodic. The Liapunov exponents are
the real parts of the eigenvalues of $B$, that is, the real parts of the
characteristic exponents. $\e^{Bt}$ has at least one eigenvalue equal to
$1$, corresponding to translations along the periodic orbit (seen in
Proposition \ref{prop_Poincare}). Thus one of the Liapunov exponents is
equal to zero. More generally, orbits on an invariant torus of dimension
$m$ will have at least $m$ Liapunov exponents equal to zero. If  one or
more of the remaining exponents are positive, then the periodic or
quasiperiodic orbits on the torus will depend sensitively on initial
conditions.  However, all the above invariant sets (equilibrium, periodic
orbit or torus) must be unstable in order to have sensitive dependence on
initial conditions, and thus they have zero measure and are not attractors.
Orbits starting near these sets will not necessarily have positive Liapunov
exponents. 

\begin{table}
\begin{center}
\begin{tabular}{|l|c|l|}
\hline
\vrule height 13pt depth 7pt width 0pt
Attractor &  Sign of Liapunov exponents & Asymptotic dynamics \\
\hline
\vrule height 13pt depth 7pt width 0pt
Stable equilibrium & $(-,-,-)$ & Stationary \\
\vrule height 0pt depth 7pt width 0pt
Stable periodic orbit & $(0,-,-)$ & Periodic \\
\vrule height 0pt depth 7pt width 0pt
Attracting torus & $(0,0,-)$ & Quasiperiodic \\
\vrule height 0pt depth 7pt width 0pt
Strange attractor & $(+,0,-)$ & Chaotic \\
\hline
\end{tabular}
\end{center}
\caption[]
{Examples of attractors of a three-dimensional flow.}
\label{t_liap}
\end{table}

By contrast, orbits of a strange attractor must attract all nearby orbits
that do {\em not} belong to the attractor, while they repel nearby orbits
that {\em do} belong to it. In order to attract nearby orbits, the flow
must be dissipative in a \nbh\ of the attractor, so that by
Proposition~\ref{prop_sdic}, the sum of all Liapunov exponents must be
negative. Since orbits on the attractor are bounded and not attracted by an
equilibrium point, one Liapunov exponent is equal to zero. Hence a
two-dimensional flow cannot admit a strange attractor.\footnote{A theorem
due to Poincar\'e and Bendixson states that a non-empty compact $\w$- or
$\alpha$-limit set of a planar flow is either a periodic orbit, or contains
equilibrium points, which also rules out the existence of strange
attractors for two-dimensional flows.} If a three-dimensional flow has a
strange attractor, its Liapunov exponents must satisfy $\lambda_1 >
\lambda_2=0 > \lambda_3$ and $\abs{\lambda_3}>\lambda_1$ (see
\tabref{t_liap}). Because of volume contraction, the attractor must have
zero volume, but it cannot be a surface and have positive Liapunov
exponents. This accounts for the fractal nature of many observed
attractors.

For iterated maps, the situation is less restrictive, since their orbits
need not have one Liapunov exponent equal to zero. Thus two-dimensional
dissipative maps may admit a strange attractor, with Liapunov exponents
satisfying $\abs{\lambda_2} > \lambda_1 > 0 > \lambda_2$. For instance, the
intersection of the strange attractor of a three-dimensional flow with a
Poincar\'e section is also a strange attractor of the associated
two-dimensional Poincar\'e map. One-dimensional iterated maps may
have strange attractors if they are non-invertible. 


\subsection{The H\'enon and Lorenz Attractors}
\label{ssec_hl}

\begin{figure}
 \centerline{\psfig{figure=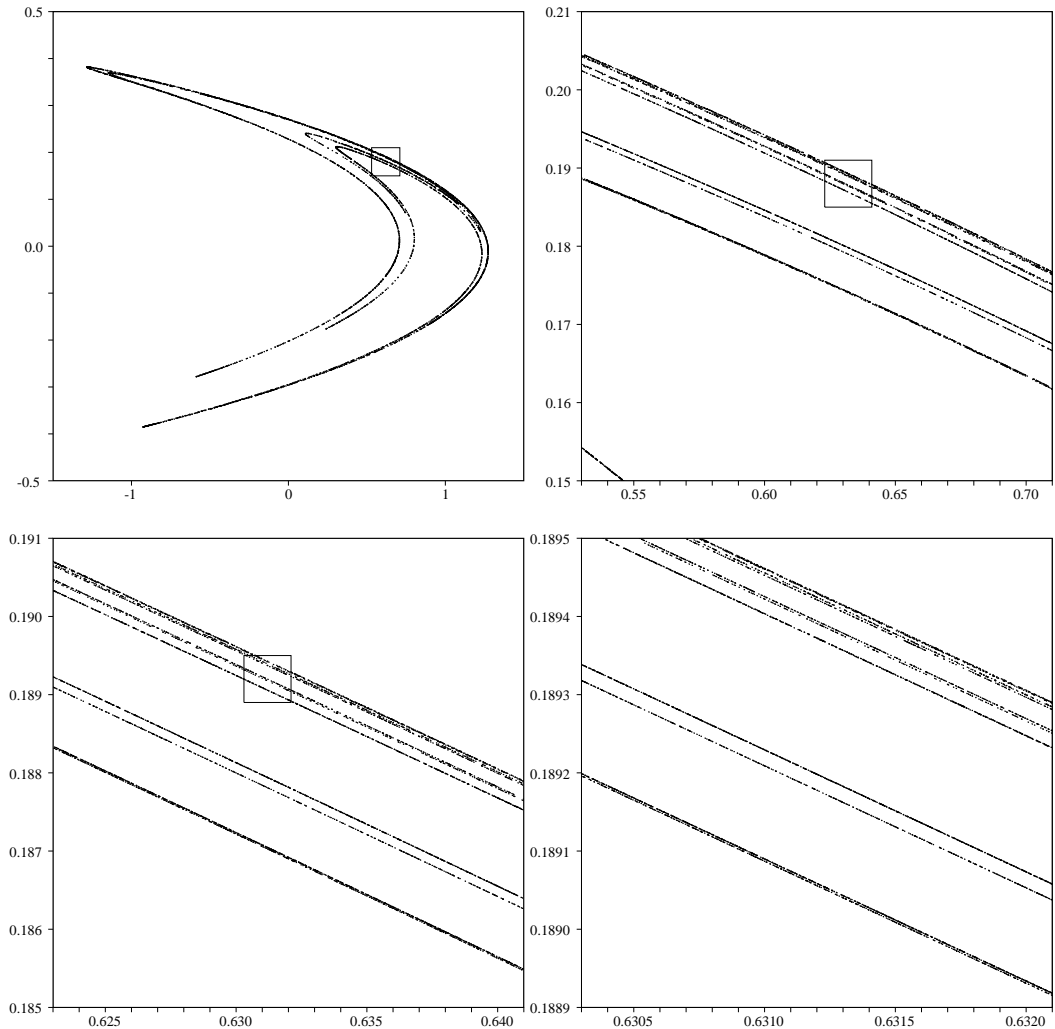,width=140mm}}
 \vspace{4mm}
 \caption[]{H\'enon attractor for $\lambda = 1.4$ and $b = 0.3$. Successive
 magnifications of details, by a factor $10$, show its self-similar
 structure.}
 \label{fig_Henon}
\end{figure}

The first dynamical system for which the existence of a strange attractor
was proved is the H\'enon map
\begin{equation}
\label{hl1}
\begin{split}
x_{k+1} &= 1 - \lambda x_k^2 + y_k \\
y_{k+1} &= b x_k.
\end{split}
\end{equation}
This map does not describe a physical system. It has been introduced as a
two-dimensional generalization of the one-dimensional map $x_{k+1}=1-
\lambda x_k^2$, which is equivalent to the logistic map. When $b=0$,
\eqref{hl1} is reduced to this one-dimensional map. If $b\neq 0$, the
H\'enon map is invertible, and it is dissipative if $\abs{b}<1$. Numerical
simulations indicate that for some parameter values, the H\'enon map has
indeed a strange attractor with a self-similar structure
(\figref{fig_Henon}). 

\begin{theorem}[Benedicks, Carleson \cite{BC}]
\label{thm_Henon}
Let $z^\star=(x^\star,y^\star)$ be the fixed point of \eqref{hl1} with
$x^\star$, $y^\star>0$, and let $\Wglo{u}$ be its unstable manifold. For
all $c<\log2$, there is a set of positive Lebesgue measure of parameters
$(b,\lambda)$ for which 
\begin{enum}
\item	there is an open set $\cU\subset\R^2$, depending on $\lambda$ and
$b$, such that for all $z\in\cU$, 
\begin{equation}
\label{hl2}
T^k(z)\to\overline{\Wglo{u}}\qquad \text{as} \qquad k\to\infty;
\end{equation}
\item	there is a point $z_0\in\Wglo{u}$ such that
\begin{enum}
\item	the positive orbit of $z_0$ is dense in $\Wglo{u}$;
\item	the largest Liapunov exponent of the orbit of $z_0$ is larger than
$c$. 
\end{enum}
\end{enum}
\end{theorem}

Hence the closure of the unstable manifold $\Wglo{u}$ is a strange
attractor. Locally, the strange attractor is smooth in the unstable
direction (with positive Liapunov exponent), and has the structure of a
Cantor set in the transverse, stable direction (with negative Liapunov
exponent). 

As a final illustration, we return to the Lorenz equations
\begin{equation}
\label{hl3}
\begin{split}
\dot x_1 &= \sigma(x_2-x_1) \\
\dot x_2 &= r x_1 - x_2 - x_1x_3 \\
\dot x_3 &= -b x_3 + x_1x_2.
\end{split}
\end{equation}
The existence of a strange attractor for this system has not been proved to
our knowledge, although there is strong numerical evidence that such an
attractor exists for certain parameter values, including in particular
$\sigma=10$, $b=8/3$ and $r=28$. The qualitative properties of dynamics
are nonetheless quite well understood (see for instance \cite{Sparrow}).
The strange attractor seems to appear after a rather subtle sequence of
bifurcations. Let us consider the case $\sigma=10$, $b=8/3$, and take $r$
as bifurcation parameter. 

\begin{figure}
 \centerline{
 \psfig{figure=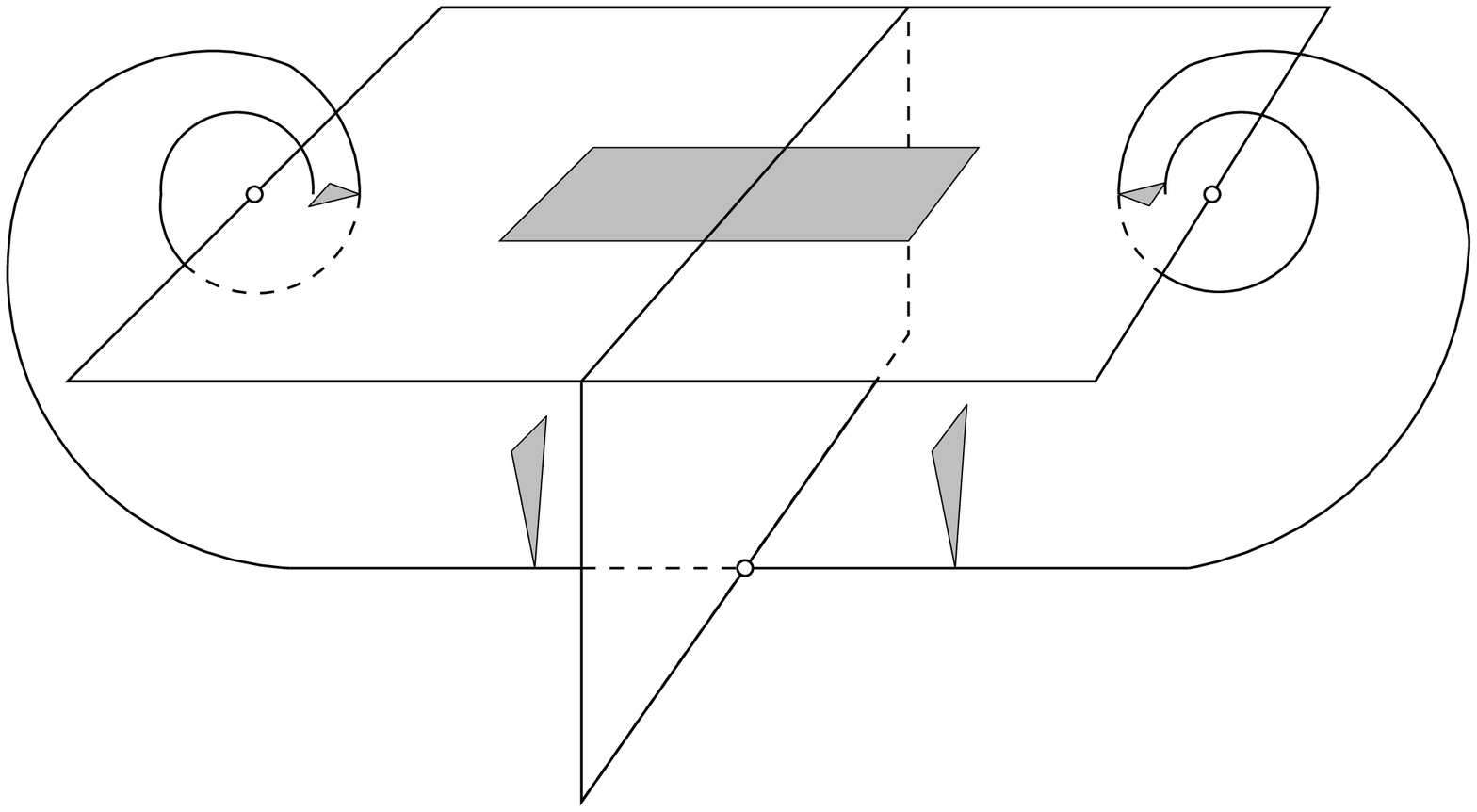,width=72.5mm}
 \hspace{10mm}
 \psfig{figure=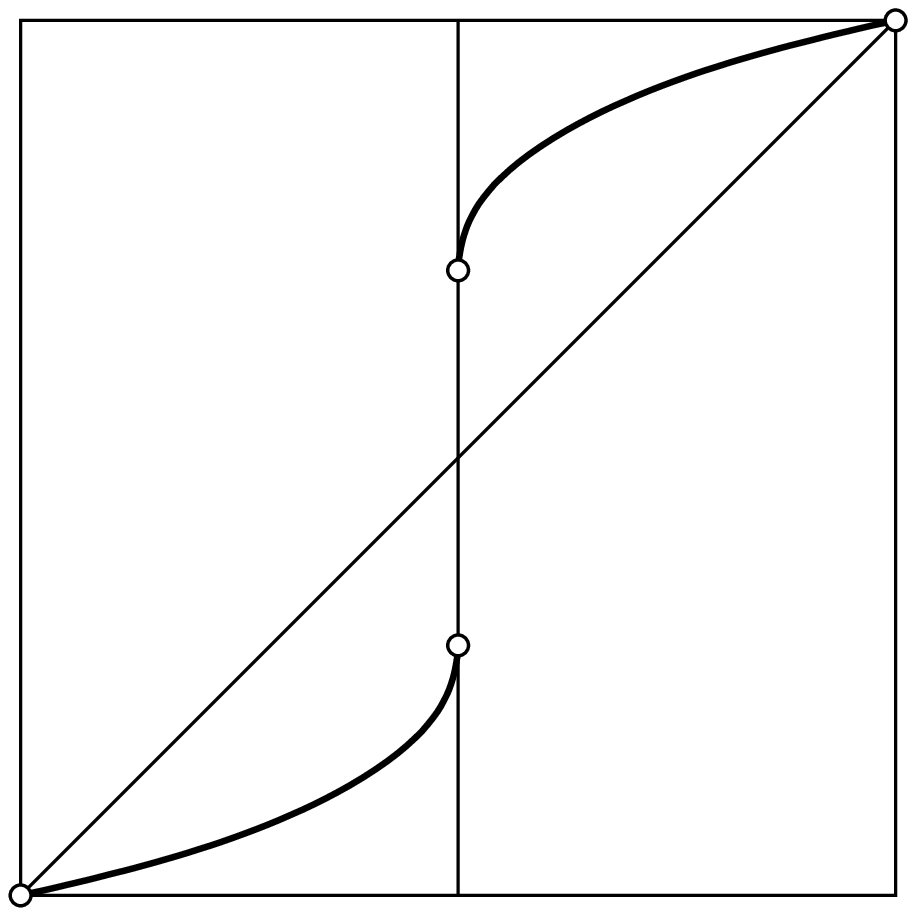,width=45mm}}
 \centerline{
 \psfig{figure=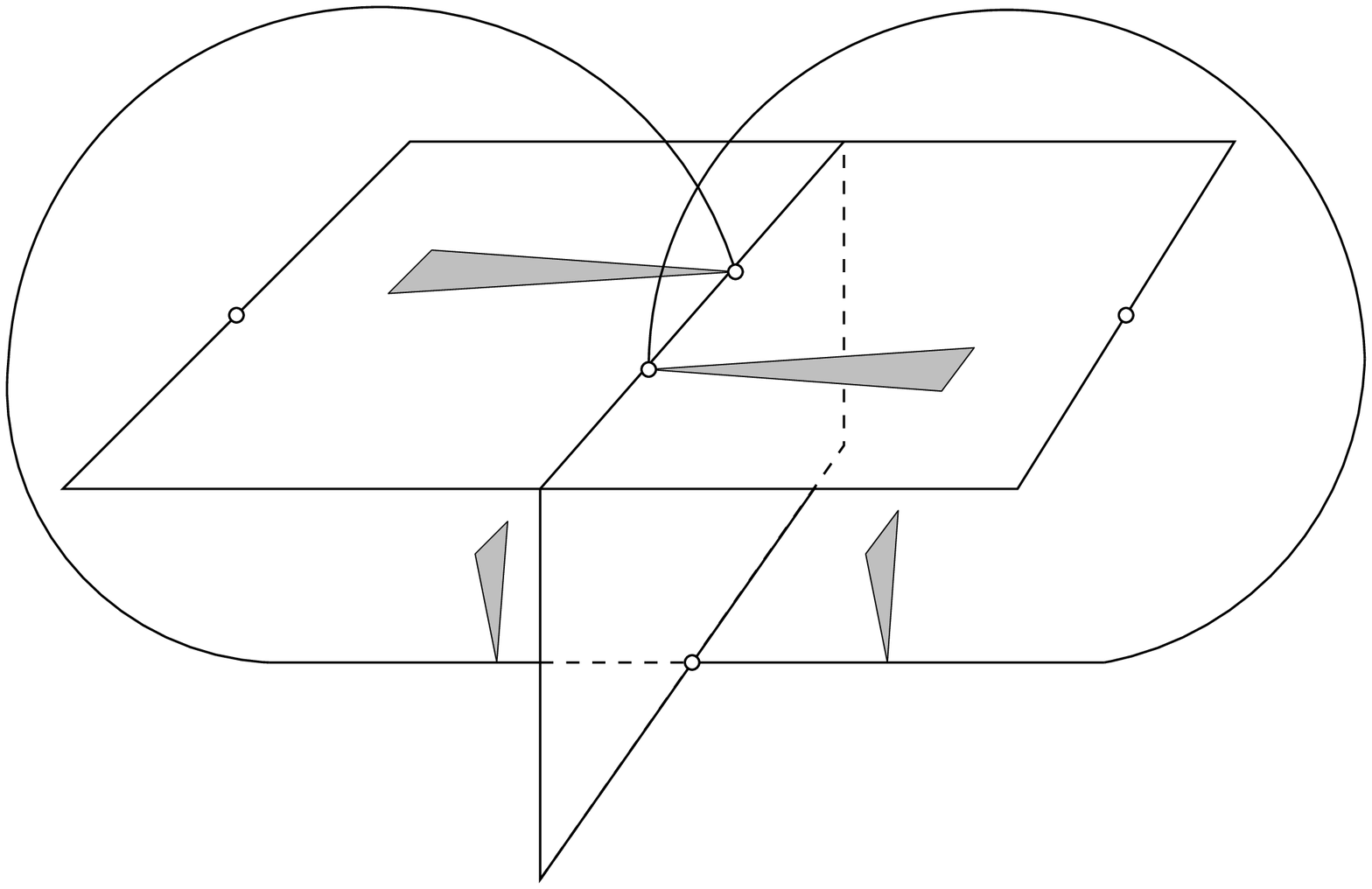,width=72.5mm}
 \hspace{10mm}
 \psfig{figure=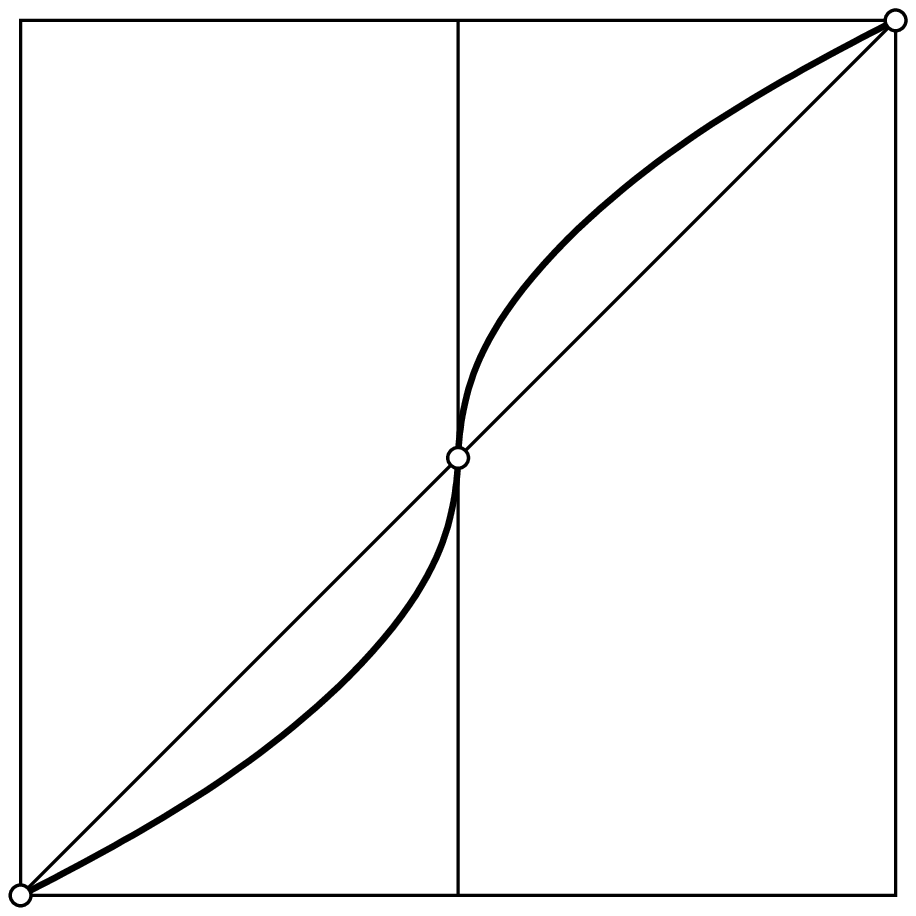,width=45mm}}
 \centerline{
 \psfig{figure=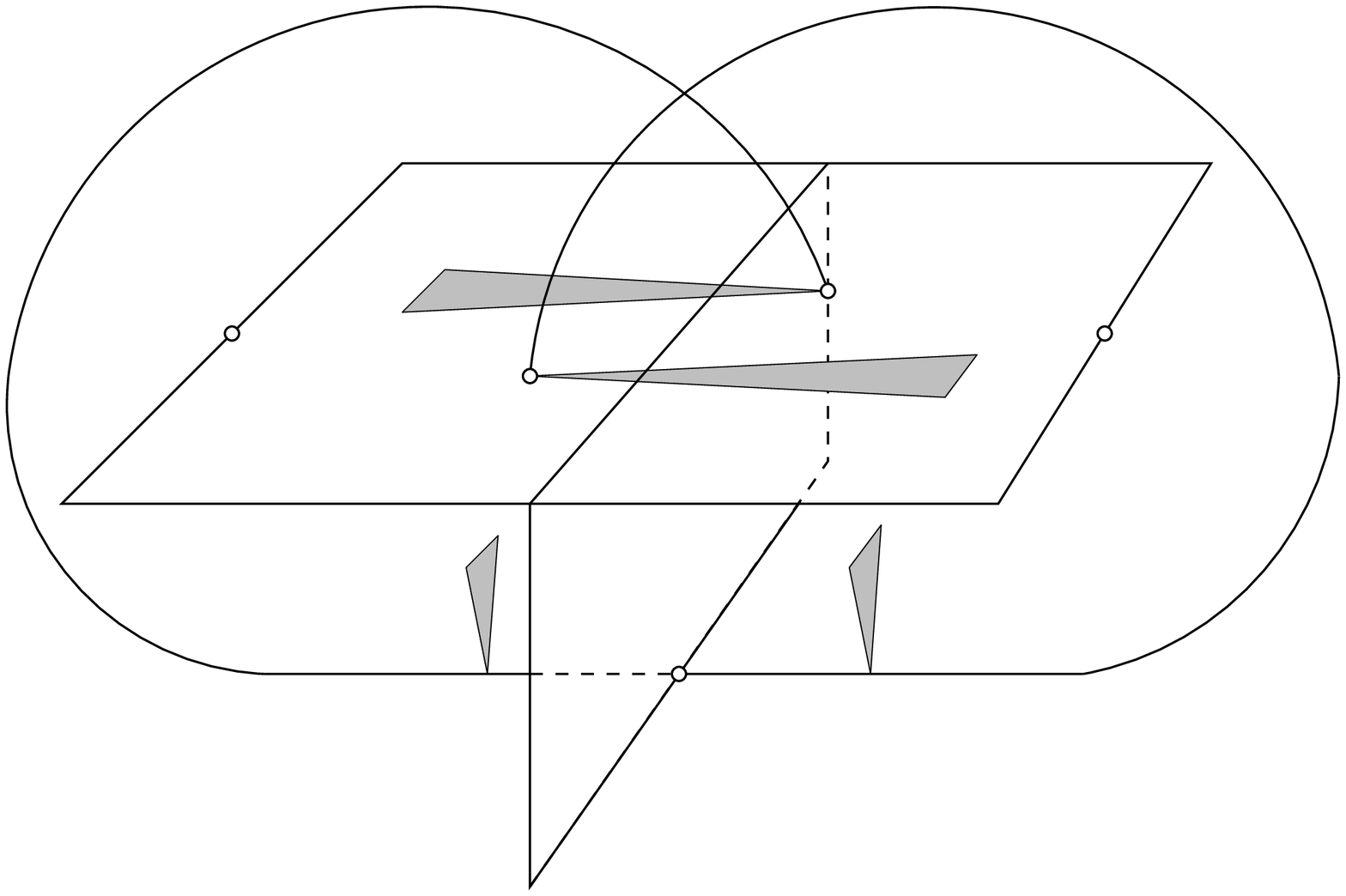,width=72.5mm}
 \hspace{10mm}
 \psfig{figure=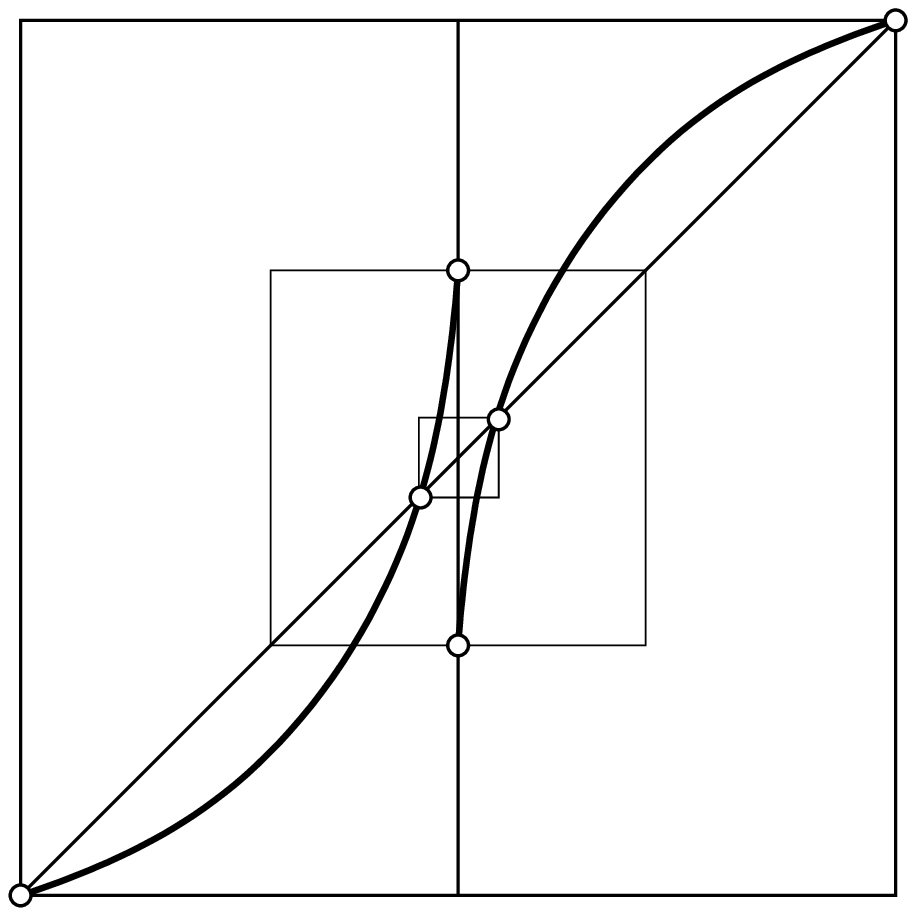,width=45mm}}
 \centerline{
 \psfig{figure=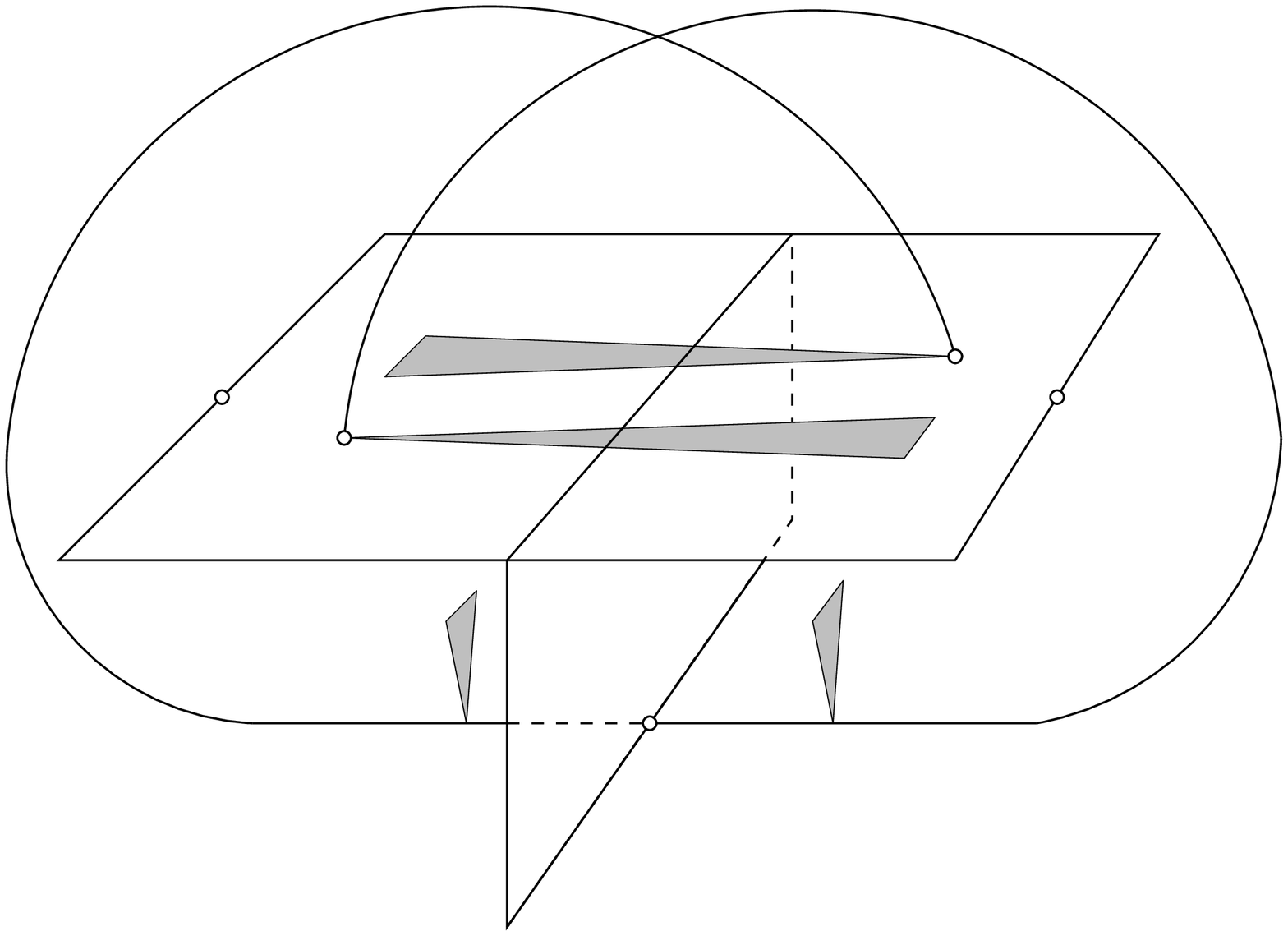,width=72.5mm}
 \hspace{10mm}
 \psfig{figure=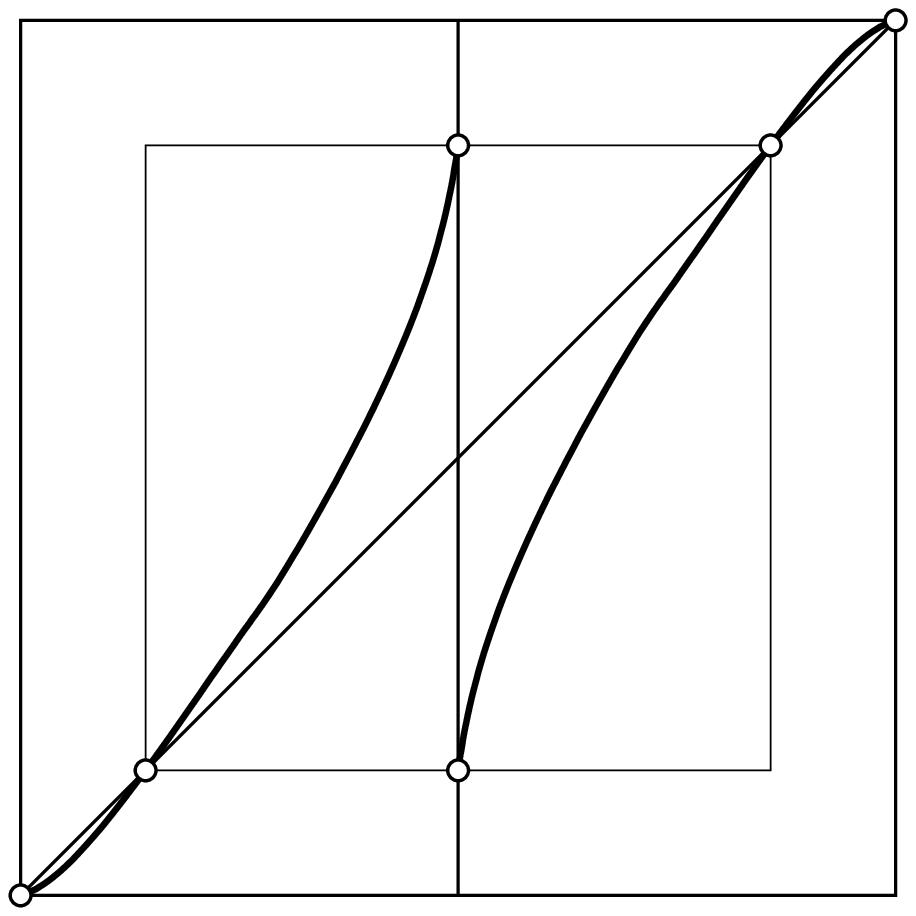,width=45mm}}
 \figtext{
 	\writefig	1.0	5.2	{d}
	\writefig	4.7	1.4	$O$
	\writefig	7.0	1.4	$\Wglo{u}$
	\writefig	3.9	2.1	$\Wglo{s}$
	\writefig	1.9	2.8	$\Sigma$
	\writefig	4.2	2.8	$S$
	\writefig	1.6	3.5	$C_-$
	\writefig	7.0	3.4	$C_+$
	\writefig	9.6	0.4	$c_-$
	\writefig	10.2	0.9	$b_-$
	\writefig	13.3	4.0	$b_+$
	\writefig	13.8	4.6	$c_+$
 	\writefig	1.0	10.55	{c}
	\writefig	4.7	6.75	$O$
	\writefig	7.0	6.75	$\Wglo{u}$
	\writefig	3.9	7.45	$\Wglo{s}$
	\writefig	1.9	8.15	$\Sigma$
	\writefig	4.2	8.15	$S$
	\writefig	1.6	8.85	$C_-$
	\writefig	7.0	8.75	$C_+$
	\writefig	9.6	5.75	$c_-$
	\writefig	11.0	7.85	$b_-$
	\writefig	12.1	8.15	$b_+$
	\writefig	13.8	9.95	$c_+$
 	\writefig	1.0	15.4	{b}
	\writefig	4.7	11.6	$O$
	\writefig	7.0	11.6	$\Wglo{u}$
	\writefig	3.9	12.3	$\Wglo{s}$
	\writefig	1.9	13.0	$\Sigma$
	\writefig	4.2	13.0	$S$
	\writefig	1.6	13.7	$C_-$
	\writefig	7.0	13.6	$C_+$
	\writefig	9.6	10.6	$c_-$
	\writefig	13.8	14.8	$c_+$
 	\writefig	1.0	19.55	{a}
	\writefig	4.7	16.35	$O$
	\writefig	7.0	16.35	$\Wglo{u}$
	\writefig	3.9	17.05	$\Wglo{s}$
	\writefig	1.9	17.75	$\Sigma$
	\writefig	4.2	17.75	$S$
	\writefig	9.6	15.35	$c_-$
	\writefig	13.8	19.55	$c_+$
  }
 \captionspace
 \caption[]{Schematic representation of the stable and unstable manifolds of
 the Lorenz equations (left), and one-dimensional approximation of the
 Poincar\'e map for the coordinate perpendicular to $\Wglo{s}$ (right). (a)
 $1<r<r_1$, (b) homoclinic bifurcation at $r=r_1$, (c) existence of a
 nonwandering repeller for $r_1<r<r_2$, and (d) bifurcation to a strange
 attractor at $r=r_2$.}
 \label{fig_Lorenzbif}
\end{figure}

\begin{enum}
\item	For $0\leqs r\leqs 1$, the origin is globally asymptotically stable,
that is, all orbits converge to the origin (c.f.\
Exercise~\ref{exo_Liapunov}). 

\item	At $r=1$, the origin undergoes a pitchfork bifurcation (c.f.\
Exercise~\ref{exo_cmanif}), and two new stable equilibria
$C_{\pm}=(\pm\sqrt{b(r-1)}, \pm\sqrt{b(r-1)}, r-1)$ appear. These points
become unstable in a Hopf bifurcation at $r=\frac{470}{19}\simeq 24.74$
(c.f.\ Exercise~\ref{exo_slin}). The Hopf bifurcation is subcritical, and
thus corresponds to the destruction of an unstable periodic orbit, that must
have been created somehow for a smaller value of $r$. 

The dynamics for $r>1$ can be described by taking a Poincar\'e section on
the surface $\Sigma: x_3=r-1$ (we only take into account intersections with
$\dot x_3<0$). Consider a rectangle containing the segment $C_-C_+$ and the
intersection $S$ of the two-dimensional stable manifold of the origin with
$\Sigma$ (\figref{fig_Lorenzbif}). $S$ is attracted by the origin, which
has the effect to pinch the rectangle and map it to two pieces of
triangular shape. One vertice of each triangle belongs to a piece of the
one-dimensional unstable manifold $\Wglo{u}$ of the origin. Due to the
dissipation, the angle at this vertice is quite small. In first
approximation, the dynamics can thus be described by a one-dimensional map
for a coordinate transverse to $S$, parametrizing the long side of the
triangles. For $r$ sufficiently small, this map is increasing,
discontinuous at $0$, and admits two stable fixed points $c_\pm$
corresponding to $C_\pm$ (\figref{fig_Lorenzbif}a). 

\item	At $r=r_1\simeq 13.296$, a homoclinic bifurcation occurs
(\figref{fig_Lorenzbif}b): the unstable manifold $\Wglo{u}$ belongs to the
stable manifold $\Wglo{s}$, and thus the sharpest vertices of both
triangles belong to $S$. The one-dimensional approximation of the
Poincar\'e map is continuous, but still monotonously increasing. 

\item	For $r$ slightly larger than $r_1$, each piece of $\Wglo{u}$ hits
$\Sigma$ on the opposite side of $S$ (\figref{fig_Lorenzbif}c). The
one-dimensional map becomes non-invertible, and has two new unstable
equilibria $b_-$ and $b_+$. The map does not leave the interval $[b_-,b_+]$
invariant, but maps two of its subintervals onto $[b_-,b_+]$. The situation
is thus similar to that of the tent map $g_3$: there exists an invariant
Cantor set containing chaotic orbits. The same can be seen to hold for the
two-dimensional Poincar\'e map, which resembles the horseshoe map. Thus
there exists a strange nonwandering set $\Lambda$, but it is repelling and
has zero measure (thus it is difficult to observe numerically). 

\item	As $r$ increases, the fixed points $b_\pm$ move towards the fixed
points $c_\pm$. It appears that for $r\geqs r_2\simeq 24.06$, the interval
$[b_-,b_+]$ is mapped into itself (\figref{fig_Lorenzbif}d). The strange
nonwandering set $\Lambda$ becomes attracting, although it does not yet
attract a full \nbh, since orbits starting in $b_\pm$ may converge to
$c_\pm$. 

\item	The subcritical Hopf bifurcation at $r=r_3\simeq 24.74$ makes the
fixed points $C_\pm$ unstable, so that finally all orbits can converge to
the nonwandering set $\Lambda$, which has become an attractor. 
\end{enum}

The method used here to describe the dynamics near a homoclinic bifurcation
by a Poincar\'e map has been applied to other systems. Transitions to
chaotic behaviour are quite frequently associated with such homoclinic
bifurcations. The Poincar\'e map in a vicinity of the unstable manifold is
described as the composition of an almost linear map, reflecting the motion
near the equilibrium point, and a nonlinear, but usually rather simple map,
reflecting the motion near $\Wglo{u}$ away from the equilibrium. The
composition of two rather innocent-looking maps contains all the ingredients
necessary for the existence of chaotic dynamics. 


\end{document}